\documentclass[12pt]{article}
\usepackage{amsmath,amssymb}
\usepackage[french]{babel}
\newtheorem{theorem}{Th{\'e}or{\`e}me}[subsection]
\newtheorem{lemma}[theorem]{Lemme}
\newtheorem{proposition}[theorem]{Proposition}
\newtheorem{coro}[theorem]{Corollaire}
\begin{document}
\def\smc{\bf}
\def\sk{\section}
\def\ssk{\subsection}
\def\sssk{\subsubsection}
\def\io{\rbrack -1,0\lbrack}
\def\oi{\lbrack 0,1\rbrack}
\def\ooi{\rbrack 0,1\rbrack}
\def\ci{C^{\infty}}
\def\inn{\}_{n\in{\mathbb N}}}
\def\rr{{\mathbb R}}
\def\cc{{\mathbb C}}
\def\hh{{\mathbb H}}
\def\pp{{\mathbb P}}
\def\qq{{\mathbb Q}}
\def\zz{{\mathbb Z}}
\def\aa{{\mathbb A}}
\def\nn{{\mathbb N}}
\def\vv{{\mathbb V}}
\def\preu{\noindent {\sl Preuve~:~}\ }
\def\bo{\partial }
\def\fl{\longrightarrow }
\def\cad{c'est-{\`a}-dire~}
\def\det{{\rm det}}
\def\trace{{\rm trace}}
\def\n{\nu}
\def\f{\phi}
\def\a{\alpha}
\def\b{\beta}
\def\g{\gamma}
\def\k{\kappa}
\def\l{\lambda}
\def\S{\Sigma}
\def\s{\sigma}
\def\g{\gamma}
\def\e{\varepsilon}
\def\D{\Delta}
\def\d{\delta}
\def\L{\Lambda}
\def\G{\Gamma}
\def\r{\rho}
\def\mb{\bar{\mu}}
\def\nb{\bar{\nu}}
\def\begitem{\begin{itemize}}
\def\enditem{\end{itemize}}
\def\iti{\item[(i)]}
\def\itii{\item[(ii)]}
\def\itiii{\item[(iii)]}
\def\itiv{\item[(iv)]}
\def\itv{\item[(v)]}
\def\itvi{\item[(vi)]}
\def\itvii{\item[(vii)]}
\def\ita{\item[(a)]}
\def\itb{\item[(b)]}
\def\itc{\item[(c)]}
\def\itd{\item[(d)]}
\def\ite{\item[(e)]}
\def\itf{\item[(f)]}
\def\itg{\item[(g)]}
\def\ith{\item[(h)]}
\def\tit{\item[-]}
\def\sup{\geq}
\def\inf{\leq}
\newcommand{\eqalign}[1]{
\begin{eqnarray*}
#1
\end{eqnarray*}}

\def\preu{{\sl Preuve:}}

\newcommand{\theo}[2]{\begin{theorem}
\label{#1}
#2
\end{theorem}}

\newcommand{\lem}[2]{
\begin{lemma}
\label{#1}
#2
\end{lemma}}

\newcommand{\cor}[2]{\begin{coro}
\label{#1}
#2
\end{coro}}

\newcommand{\pro}[2]{
\begin{proposition}
\label{#1}
#2
\end{proposition}}

%%%%%%%%%%%%%%%%%%%%%%%%%%%%%%%%%%%%%%%%%%%%%%%%%%%%%%%%%%%%%%%%%%%%%%%%%%%%
\def\ma{\cite{L1}~}
\def\u{ \left \Vert p(u)\right \Vert}
\def\v{\left \Vert p(Ju)\right \Vert}
\def\Sif{S_{\infty}}
\def\Siif{\S_{\infty}}     
\def\nin{n\in{\mathbb N}}
\def\Sinn{\{\S\inn}
\def\Snn{\{S_{n}\inn}\
\def\xnn{\{x_{n}\inn}
\def\ynn{\{y_{n}\inn}
\def\wnn{\{w_{n}\inn}
\def\vnn{\{v_{n}\inn}
\def\snn{\{s_{n}\inn}

\def\dS{\partial S}
\def\dSi{\partial\Sigma}
\def\dgt{\dot\gamma}
\def\Nn{N_{\S_{n}}}
%%%%%%%%%%%%%%%%%%%%%%%%%%%%%%%%%%%%%%%%%%%%%%%%%%%%%%%%%%%%

\def\m{\mu}
\def\di{\partial_\infty M}
\def\mdi{{\cal M}(\di)}
\def\bi{{\cal B}_\infty}
\def\G{\Gamma}
\def\g{\gamma}
\def\bm{{\bar\mu}}
\def\D{\Delta}
\def\d3{\partial_{3}M}
\def\cpi{{\mathbb C}{\mathbb P}^{1}}
\def\b0i{{\cal B}^{0}_\infty}
\def\bbf{{\cal B}^{f}}
\def\bbm{{\cal B}^{m}}
\def\h2{{\mathbb H}^{2}}
\def\auteur{
\vskip 2truecm
\centerline{Fran\c cois Labourie}
\centerline{Topologie et Dynamique}
\centerline{Universit{\'e} Paris-Sud}
\centerline{F-91405 Orsay (Cedex)}}

\hsize=15truecm
\vsize=23truecm
\parindent 1,2cm

\title{Un lemme de Morse \\  pour les surfaces convexes}
\author{Fran{\c c}ois LABOURIE \thanks{l'auteur remercie l'Institut Universitaire de France}}
\maketitle
\section{Introduction}

Le lemme de Morse pour les g{\'e}od{\'e}siques affirme 
que dans une vari{\'e}t{\'e} {\`a} courbure strictement n{\'e}gative et g{\'e}om{\'e}trie born{\'e}e, toute quasi-g{\'e}od{\'e}sique
est {\`a} une distance born{\'e}e d'une g{\'e}od{\'e}sique. 
 
D'un point de vue dynamique, il est une version g{\'e}om{\'e}trique du "Shadowing Lemma". Il a pour 
cons{\'e}quence les trois propri{\'e}t{\'e}s hyperboliques du flot g{\'e}od{\'e}sique d'une vari{\'e}t{\'e} 
compacte {\`a} courbure strictement n{\'e}gative, vu ici un feuilletage de dimension 1 :
\begitem
\iti l'ensemble des feuilles compactes est dense,

\itii une feuille g{\'e}n{\'e}rique est dense,

\itiii les flots g{\'e}od{\'e}siques de deux vari{\'e}t{\'e}s {\`a} courbure n{\'e}gative
proches sont conjugu{\'e}s.
\enditem
\par

Le but de notre d'article est d'{\'e}noncer un lemme de Morse pour les
$k$-surfaces dans les vari{\'e}t{\'e}s de dimension 3 {\`a} courbure n{\'e}gative 
et d'en tirer des cons{\'e}quences analogues.

D{\'e}finissons bri{\`e}vement une $k$-surface, o{\`u} $k\in ]0,1[$, dans une vari{\'e}t{\'e} $M$ de dimension 3 
{\`a} courbure plus petite que -1 : il s'agit d'une surface dont la {\it courbure extrins{\`e}que}, 
\cad le produit des courbures 
principales, vaut $k$. Nos hypoth{\`e}ses sur $k$ entra{\^\i}nent qu'une $k$-surface est localement 
convexe et {\`a} courbure intrins{\`e}que strictement n{\'e}gative. Nous nous int{\'e}resserons principalement
aux $k$-surfaces non compactes ( elles le sont toutes si $M$ est simplement connexe ) et 
compl{\`e}tes dans un certain sens ( voir \ref{ksurf} pour des
pr{\'e}cisions). Du point de vue analytique, il s'agit d'un probl{\`e}me
elliptique, que nous d{\'e}crivons comme de Monge-Amp{\`e}re au sens de
\cite{L1}.

Le lemme de Morse pour les surfaces convexes nous dit en particulier que
toute surface
localement convexe {\`a} courbure extrins{\`e}que plus grande que $k$, ayant de 
bonnes propri{\'e}t{\'e}s {\`a} l'infini, est {\`a} distance born{\'e}e 
d'une $k$-surface. Il r{\'e}soud en particulier une sorte de probl{\`e}me de Plateau pour
les $k$-surfaces. Malheureusement, pour {\'e}noncer correctement ce lemme, il nous faut
introduire des d{\'e}finitions techniques qui sortent du cadre de cette
pr{\'e}sentation (voir \ref{lemmorse}).

Nous allons pr{\'e}senter maintenant deux types d'application de ce lemme de Morse. Nous {\'e}tudierons les 
probl{\`e}mes de Plateau asymptotiques et les propri{\'e}t{\'e}s hyperboliques 
de l'espace des $k$-surfaces. Ces propri{\'e}t{\'e}s hyperboliques, qui font
appara{\^\i}tre l'espace des $k$-surfaces comme une g{\'e}n{\'e}ralisation du flot
g{\'e}od{\'e}sique, sont la motivation essentielle de notre travail.

Pour {\'e}noncer de fa\c con simple nos r{\'e}sultats, nous 
supposerons dans cette introduction que $M$ est le rev{\^e}tement universel d'une vari{\'e}t{\'e} $N$ compacte.

\subsection{Probl{\`e}me de Plateau asymptotique}

Soit $i$ un hom{\'e}omorphisme local d'une surface $S$ dans le bord {\`a} l'infini, $\bo_\infty M$, de $M$. Une {\it solution} 
du probl{\`e}me de Plateau asymptotique d{\'e}fini par $(i,S)$, est une immersion $f$ de $S$ dans 
$M$ dont l'image est une $k$-surface et telle que $i$ est l'application qui {\`a} un point de $S$
associe l'extr{\'e}mit{\'e} {\`a} l'infini de la normale ext{\'e}rieure {\`a} $f(S)$ en ce point.

La terminologie se justifie par le fait que lorsque $i$ est un plongement
du disque ouvert qui s'{\'e}tend {\`a} un plongement du disque ferm{\'e}, alors on a
envie de penser que le bord {\`a} l'infini de la solution coinc{\"\i}de avec le bord
du disque. Elle peut-{\^e}tre aussi source de confusion car en g{\'e}n{\'e}ral une
solution d'un probl{\`e}me asymptotique n'a pas de bord {\`a} l'infini en un
sens raisonnable. 

Nous d{\'e}montrons une s{\'e}rie de r{\'e}sultats sur ce probl{\`e}me de Plateau asymptotique.

\par\noindent{\smc Th{\'e}or{\`e}me A.~~}{\sl Il existe au plus une solution du probl{\`e}me de 
Plateau asymptotique.}\par

Le probl{\`e}me de Plateau asymptotique n'a pas toujours de solutions :   

\par\noindent{\smc Th{\'e}or{\`e}me B.~~}{\sl Si $S$ est $\bo_\infty M$ auquel on a {\^o}t{\'e} 0, 1, ou 2 points et si
$i$ est l'injection canonique, alors le probl{\`e}me $(i,S)$ n'a pas de solutions.}\par

Nous avons les r{\'e}sultats d'existence suivants :

\par\noindent{\smc Th{\'e}or{\`e}me C.~~}{\sl Si $(i,S)$ est un probl{\`e}me de Plateau asymptotique et si
$\bo_\infty M\setminus i(S)$ contient au moins trois points distincts, alors $(i,S)$
admet une solution.}\par

\par\noindent{\smc Th{\'e}or{\`e}me D.~~}{\sl Si $\Gamma$ est un groupe agissant sur $S$, tel que $S/\Gamma$ 
soit une surface compacte de genre plus grand ou {\'e}gal {\`a} 2, si $\rho$ est une repr{\'e}sentation de $\Gamma$ dans 
le groupe des  isom{\'e}tries de $M$, et $i$ v{\'e}rifie
$$
\forall \gamma \in \Gamma,~~~~ i\circ\gamma=\rho(\gamma)\circ i,
$$
alors le probl{\`e}me $(i,S)$ a une solution.}\par

Enfin

\par\noindent{\smc Th{\'e}or{\`e}me E.~~}{\sl Si $(i,U)$ est un probl{\`e}me de Plateau asymptotique et si $S$ 
est un ouvert relativement compact de $U$ alors $(i,S)$ admet une solution.}\par

On peut remarquer {\`a} ce stade une analogie entre les donn{\'e}es de probl{\`e}mes de Plateau 
asymptotiques ayant 
des solutions, et les immersions hyperboliques du disque dans ${\mathbb CP}^1$.

Nous avions d{\'e}montr{\'e} dans \cite{M} le th{\'e}or{\`e}me D dans le cas o{\`u} la vari{\'e}t{\'e} ambiante {\'e}tait 
{\`a} courbure constante.
Toujours dans le cas de la courbure ambiante contante, mais en toutes
dimensions, J.~Spruck et H.~Rosenberg ont d{\'e}montr{\'e} des versions
partielles de nos r{\'e}sultats : 
ils ont montr{\'e} que le 
probl{\`e}me de Plateau d{\'e}fini par $(i,U)$ avait une solution dans le cas o{\`u} $U$ {\'e}tait 
l'int{\'e}rieur du disque ferm{\'e} $D$ et o{\`u}  $i$ s'{\'e}tendait 
en un plongement de $D$ dans $U$ \cite{SR}. Les autres r{\'e}sultats de cet article sont des
cas particuliers, toujours dans le cadre de la courbure ambiante 
constante ( mais en toutes dimensions ),
de \ref{platdisk} : les auteurs ne consid{\`e}rent que les graphes au dessus d'une horosph{\`e}re.
Ceci leur
permet d'{\'e}crire
explicitement  les 
{\'e}quations satisfaites par la fonction dont  une telle surface est le graphe, {\'e}quations qui, 
gr{\^a}ce {\`a} l'hypoth{\`e}se de coubure ambiante constante, sont d'une forme simple {\`a} laquelle peuvent 
s'appliquer les estim{\'e}es {\it a priori}, classiques pour les {\'e}quations de
Monge-Amp{\`e}re. Quant {\`a} nous, nous appliquons les techniques plus flexibles 
de courbes pseudo-holomorphes mises au point dans \cite{L1}. Ceci nous
permet de soritir du cadre de la courbure constant, et de consid{\'e}rer des
surfaces qui ne sont pas des graphes au dessus des horosph{\`e}res.
\subsection{Propri{\'e}t{\'e}s hyperboliques de l'espace des $k$-surfaces}

Nous voulons poursuivre notre analogie avec les g{\'e}od{\'e}siques plus loin. Soit maintenant $N$ une 
vari{\'e}t{\'e} de dimension 3 compacte. Le fibr{\'e} unitaire de $N$ est l'espace des g{\'e}od{\'e}siques 
point{\'e}es, \cad des paires $(\gamma ,x)$ o{\`u} $x$ est un point de la
g{\'e}od{\'e}sique ${\gamma}$. Le flot g{\'e}od{\'e}sique correspond alors au feuilletage
de dimension 1 obtenu en faisant bouger le point le long d'une g{\'e}od{\'e}sique
donn{\'e}e.

D{\'e}finissons donc les $k$-surfaces point{\'e}es, comme les paires $(S,x)$ o{\`u}
$x$ est un point sur la  $k$-surface $S$. Nous attirons encore une fois l'attention 
sur le fait qu'il s'agit de surfaces en g{\'e}n{\'e}ral non compactes et ayant 
m{\^e}me beaucoup de r{\'e}currence. Nous d{\'e}finissons {\'e}galement les tubes point{\'e}s 
comme les paires
$(T,x)$, o{\`u} $T$ est une g{\'e}od{\'e}sique et $x$ est  un vecteur normal unitaire {\`a} $T$. 

Dans \cite{L1}, dans le cadre plus g{\'e}n{\'e}ral des probl{\`e}mes de Monge-Amp{\`e}re, nous avons montr{\'e}
l'espace ${\cal N}$,  constitu{\'e} des $k$-surfaces point{\'e}es et des tubes point{\'e}s, est 
compact,
et poss{\`e}de de plus une structure de lamination d{\'e}finie par la condition : deux points $(F,x)$ et $(G,y)$
sont sur la m{\^e}me feuille si $F=G$.

Par abus le langage, nous appelerons ${\cal N}$ l'espace des $k$-surfaces  de $N$.
 
Nous donnerons plus de pr{\'e}cisions sur la topologie sous-jacente dans \ref{laminations} et renvoyons
surtout {\`a} \cite{L1}.

Le r{\'e}sultat essentiel de notre article est le th{\'e}or{\`e}me suivant qui montre que l'espace ${\cal N}$ est une 
bonne g{\'e}n{\'e}ralisation du flot g{\'e}od{\'e}sique, car il en poss{\`e}de les
propri{\'e}t{\'e}s hyperboliques.

\par\noindent{\smc Th{\'e}or{\`e}me F.~~}{\sl Soit ${\cal N}$ l'espace des $k$-surfaces d'une vari{\'e}t{\'e} compacte $N$
munie d'une m{\'e}trique $g$ {\`a} courbure plus petite que -1, alors
\begitem
\iti pour tout entier $g$, l'ensemble des feuilles compactes de genre plus grand que $g$ 
est dense dans ${\cal N}$,
\itii une feuille g{\'e}n{\'e}rique est dense dans ${\cal N}$,
\itiii enfin ${\cal N}$ est stable dans le sens suivant, si $h$ est une m{\'e}trique 
suffisamment proche de $g$, et si  ${\bar{\cal N}}$ est l'espace des $k$-surfaces
pour la m{\'e}trique $h$, alors il existe un hom{\'e}omorphisme entre  ${\cal N}$ et  ${\bar{\cal N}}$
envoyant feuille sur feuille
\enditem
.}\par

L'espace des $k$-surfaces contient l'espace des tubes point{\'e}s, espace
qui fibre au-dessus du flot g{\'e}od{\'e}sique. En ce sens, l'espace des
$k$-surfaces ``contient le flot g{\'e}od{\'e}sique''. Dans le cas, o{\`u} la
courbure ambiante est constante, il contient {\'e}\-ga\-le\-ment l'espace des
plans totalement g{\'e}od{\'e}siques ${\cal P}$ : les surfaces {\'e}quidistantes des
plans totalement g{\'e}od{\'e}siques sont {\`a} courbure constante. Notre espace
${\cal N}$ de dimension infinie est cependant  beaucoup plus gros.

Signalons ici {\`a} titre de comparaison  
quelques propri{\'e}t{\'e}s de l'espace ${\cal P}$ des plans
totalement g{\'e}od{\'e}siques toujours dans le cas de la courbure ambiante
constante. Tout d'abord, par le th{\'e}or{\`e}me d'ergodicit{\'e} de Moore, une feuille g{\'e}n{\'e}rique est dense. 
En ce qui concerne les feuilles compactes, il est connu que pour certaines vari{\'e}t{\'e}s
arithm{\'e}tiques ${\cal P}$ ne contient aucune feuilles compactes, pour
d'autres il en contient un ensemble dense. 

Cet espace ${\cal P}$ des plans totalement g{\'e}od{\'e}siques, n'a bien s{\^u}r pas beaucoup d'int{\'e}r{\^e}t quand la courbure n'est
pas constante, mais il poss{\`e}de une sorte de propri{\'e}t{\'e} de ``stabilit{\'e}''
d{\'e}montr{\'e}e par  M.~Gromov 
\cite{F} : si une m{\'e}trique $h$ sur
une vari{\'e}t{\'e} $N$ est suffisamment proche d'une m{\'e}trique $g$ {\`a} courbure
constante, alors il existe une application continue de ${\cal P}$, l'espace des
plans totalement g{\'e}od{\'e}siques pour $g$,  dans $N$ telle que les images
des
feuilles soient des surfaces minimales. 

La situation d{\'e}crite par le th{\'e}or{\`e}me F am{\`e}ne {\`a} se poser un certain
nombre de questions :

(i) on a une abondance de mesures tranverses invariantes, mais y en
a-t-il 
une qui charge tous les
ouverts ?

(ii) Quelle est la statistique des feuilles compactes ? On peut associer
{\`a} chaque $k$-surface l'int{\'e}grale de sa courbure moyenne, pour de bonnes
raisons nous appelerons ceci {\it l'aire} de la $k$-surface. De mani{\`e}re
naturelle, l'aire des tubes est alors la longueur de la g{\'e}od{\'e}sique
sous-jacente ( {\`a} $2\pi$ pr{\`e}s ). On peut alors
montrer ( ce n'est pas fait dans cet article ) qu'il n'y a qu'un nombre
fini $N(A)$ de feuilles compactes d'aire born{\'e}e par $A$. Notons $N(A,g)$
le nombre de celles-ci qui sont de genre $g$. On a envie de poser
$$
ent_g (N)=\liminf_{h\leq g,~ A\rightarrow \infty} \frac{log (N(A,g))}{A}.
$$

Le nombre $ent_0$ est l'entropie du flot g{\'e}od{\'e}sique, quelles sont les
valeurs de $ent_g$ ?

(iii) Les feuilles compactes sont elles {\'e}quidistribu{\'e}es ?

(iv) On peut {\'e}galement se poser une question qui est l'analogue de la
question de conjugaison pour les flots g{\'e}od{\'e}siques : pour deux m{\'e}triques
proches, supposons que l'hom{\'e}omorphisme de conjugaison puisse {\^e}tre
choisi conforme sur chacune des feuilles, les m{\'e}triques sont-elles
isom{\'e}triques ?

\subsection{Structure de l'article}
\noindent{\bf 2}- {\sl D{\'e}finitions, {\'e}nonc{\'e} du lemme de Morse.} Nous {\'e}non{\c c}ons 
le lemme de Morse, ainsi que les d{\'e}finitions n{\'e}cessaires. Cette section contient {\'e}galement
des d{\'e}\-fi\-ni\-tions utilis{\'e}es de mani{\`e}re r{\'e}currente dans l'article ainsi que quelques
propri{\'e}t{\'e}s pr{\'e}\-li\-mi\-nai\-res.

\noindent {\bf 3} - {\sl Variations infinit{\'e}simales et d{\'e}formations.} Nous y d{\'e}montrons 
le lemme \ref{} qui d{\'e}crit les variations infinit{\'e}simales de $k$-surfaces {\`a} bord.

\noindent{\bf 4}- {\sl Th{\'e}or{\`e}me de compacit{\'e}.} Nous y {\'e}tudions les limites de $k$-surfaces.
Cette section d{\'e}bute par un rappel des r{\'e}sultats de compacit{\'e} sur les probl{\`e}mes de Monge-Amp{\`e}re
d{\'e}montr{\'e}s dans \cite{L1}.

\noindent{\bf 5} - {\sl Probl{\`e}me de Plateau pour les disques.} Nous y d{\'e}montrons la 
proposition \ref{platdisk} qui est une version faible du lemme de Morse. 

\noindent{\bf 6}- {\sl D{\'e}monstration du lemme de Morse pour les surfaces convexes.}
Nous d{\'e}\-mont\-rons l'existence et l'unicit{\'e}, en utilisant une m{\'e}thode de d{\'e}formation.

\noindent{\bf 7}- {\sl Probl{\`e}mes asymptotiques.}
Nous y d{\'e}montrons les th{\'e}or{\`e}mes {\'e}nonc{\'e}s dans l'introduction. Nous traitons le th{\'e}or{\`e}me A
en \ref{uniasym} ; les th{\'e}or{\`e}mes C,D et E sont vus en \ref{picard},
\ref{equib} et \ref{pasequi} respectivement, et enfin B est d{\'e}montr{\'e} en \ref{inex}.

\noindent{\bf 8}- {\sl Espace des $k$-surfaces.} A partir de cette section, nous nous 
int{\'e}ressons {\`a} l'espace des $k$-surfaces (sans bord) d'une vari{\'e}t{\'e} compacte. Cet espace est introduit en
\ref{}. 

\noindent{\bf 9}- {\sl Densit{\'e} des feuilles p{\'e}riodiques.} Il s'agit de de F-(i). 
La d{\'e}monstration est la plus d{\'e}licate de cet article.

\noindent{\bf 10}-{\sl G{\'e}n{\'e}ricit{\'e} des feuilles denses.} Nous y d{\'e}montrons F-(ii).

\noindent{\bf 11}-{\sl Stabilit{\'e}.} Cette section contient la preuve de F-(iii).
%%%%%%%%%%%%%%%%%%%%%%%%%%%%%%%%%%%%%%%%
\section{D{\'e}finitions}

Dans cette section, nous allons pr{\'e}senter les d{\'e}finitions utilis{\'e}es dans
l'{\'e}\-non\-c{\'e} du lemme de Morse. Nous donnerons {\'e}galement quelques
d{\'e}finitions connexes utilis{\'e}es dans
les preuves, ainsi que quelques r{\'e}sultats pr{\'e}liminaires.

Enon\c cons tout le suite le lemme de Morse qui sera d{\'e}montr{\'e} en
\ref{morse}
\par

\noindent{\smc Lemme de Morse. }\label{lemmorse}{\sl Soit $M$ une vari{\'e}t{\'e}
d'Hadamard {\`a} g{\'e}om{\'e}trie born{\'e}e et {\`a}
courbure strictement plus petite que $-c<0$. Soit  $S$  une surface
localement convexe, {\'e}ventuellement {\`a} bord, {\`a} courbure plus grande que $c$, {\`a}
g{\'e}om{\'e}trie born{\'e}e et qui n'est ni horosph{\'e}rique {\`a} l'infini, ni
tubulaire, ni compacte sans bord.
 Alors, pour tout $k\in \rbrack 0, c\lbrack$,
il existe une unique $k$-surface lentille
pour $S$.

De plus, si  $S$ n'est pas tubulaire {\`a} l'infini, cette
$k$-surface a sa courbure moyenne uniform{\'e}ment born{\'e}e et en particulier
n'est pas d{\'e}g{\'e}n{\'e}r{\'e}e}\par

Les termes de ce lemme vont {\^e}tre expliqu{\'e}s dans cette section.

Nous parlerons des $k$-surfaces d{\'e}g{\'e}n{\'e}r{\'e}es et non d{\'e}g{\'e}n{\'e}r{\'e}es dans  le
paragraphe
\ref{ksurf}. Nous d{\'e}finirons les surfaces lentilles en \ref{deflent}
et nous aurons besoin de \ref{defbout} pour cela. Le paragraphe sur la
g{\'e}om{\'e}trie born{\'e}e \ref{geoborn} est n{\'e}cessaire pour donner la d{\'e}finition
\ref{deftubhor} de tubulaire,
horosph{\'e}rique {\`a} l'infini {\it etc ...}. 

Cette section contient {\'e}galement des r{\'e}sultats souvent  utilis{\'e}s cet article.

\subsection{$k$-surfaces}

Si $S$ est une surface convexe immerg{\'e}e dans une vari{\'e}t{\'e} de dimension 3
 $M$,
 nous noterons $n(S)$ son {\it
relev{\'e}
de Gauss}, \cad la surface immerg{\'e}e dans le fibr{\'e} unitaire $UM$ de $M$
 constitu{\'e}e des vecteurs normaux ext{\'e}rieurs {\`a} $S$.

Une {\sl $k$-surface} est une surface ({\'e}ventuellement  {\`a} bord) immerg{\'e}e dans
$M$ dont le produit des courbures principales vaut $k$, et telle que la
m{\'e}trique induite de $n(S)$ soit compl{\`e}te. Attention, ceci n'entra{\^\i}ne pas
{\it a priori}  que la m{\'e}trique induite de celle de $M$ est compl{\`e}te. Si
tel est le cas , la surface est {\it non d{\'e}g{\'e}n{\'e}r{\'e}e}, et {\it d{\'e}g{\'e}n{\'e}r{\'e}e}
dans le cas contraire. Si la courbure moyenne est born{\'e}e, alors la
$k$-surface est non d{\'e}g{\'e}n{\'e}r{\'e}e.

\subsection{Normal {\'e}tendu, surfaces lentilles, bout}\label{ksurf}

\subsubsection{Normal {\'e}tendu, bout}\label{defbout}
Soit $S$ une surface localement convexe immerg{\'e}e {\`a} bord $\bo S$ dans
$M$. Soit $n$ son champ de vecteur normal ext{\'e}rieur et $n_{\bo}$ le
champ de vecteur normal int{\'e}rieur {\`a} $\bo S$ dans $S$. Posons
$$
N_{\bo}^{+}=\{u\in (T{\bo S})^{\perp}/~~\langle u,n_{\bo}\rangle\leq 0 \}.
$$
Le {\it normal {\'e}tendu}
de $S$ est le sous-ensemble $N_{S}$ de $UM$, d{\'e}fini par  
$$
N_{S}=n(S)\cup N_{\bo}^{+}.
$$  

Dans la suite, nous identifierons
souvent abusivement $n(S)$ ( vu comme sous-ensemble de $UM$ ) et $S$. Il est facile de voir que $N_{S}$ est
une sous-vari{\'e}t{\'e} $C^{0}$-immerg{\'e}e dans $UM$ dont le bord $\bo N_{S}$ est
l'ensemble de vecteurs normaux int{\'e}rieurs {\`a} $S$ le long de $\bo S$.

Le {\sl bout} $B$ de $S$ sera la vari{\'e}t{\'e} de dimension 3, hom{\'e}omorphe {\`a}
$N_{S}\times ]0,+\infty[$ et muni de la m{\'e}trique induite par
l'application 
$$
\left\{ \begin{array}{c}
N_{S}\times ]0,+\infty[\to M \\ (t,n)\mapsto exp(tn)
\end{array}\right.
$$
D'un point de vue m{\'e}trique, on compl{\'e}te le bout en ajoutant $S$, on
parlera alors d'un bout complet.

\subsubsection{Graphe {\'e}tendu, surface lentille, champ focal, pied, fonctions inverses}{\label{deflent}}
Une surface $\S$ localement convexe immerg{\'e}e {\`a} bord $\bo\S$  est un {\it graphe {\'e}tendu}
au-dessus de $S$ s'il existe 
\begitem
\item[-] un ouvert $U$ de  $int(N_{S})=N_{S}\setminus\bo N_{S}$,  v{\'e}rifiant 
$$
{\bar{n(S)}}\subset U \subset{\bar U}\subset int( N_{S})
$$ 
o{\`u} ${\bar A}$
d{\'e}signe l'adh{\'e}rence de $A$ dans $N_{S}$, 
\item[-] une fonction $f$ continue d{\'e}finie sur
${\bar U}$,
strictement positive sur $U$,
\enditem
tels que
\begitem
\iti $\forall y\in {\bar U}\setminus U,~~f(y)=0$,
\itii $\S =\{ \exp(f(u)u)~~ /~~u\in U\}$,
\itiii le segment g{\'e}od{\'e}sique $\exp([f(u),+\infty[u)$ est ext{\'e}rieur {\`a}
$\S$.
\enditem
\par
Enfin,  si $\S$ est  le graphe {\'e}tendu pour une fonction $f$  au-dessus
de $S$, nous dirons que  $\S$ est un graphe {\'e}tendu {\it born{\'e}} au-dessus
de $S$ si $f$ est born{\'e}e et de mani{\`e}re sym{\'e}trique  que
$S$ est {\it lentille} pour $\S$, et nous appellerons $f$ la {\it fonction associ{\'e}e}.  

Le {\it champ  de vecteur
 focal} pour un graphe {\'e}tendu $\S$ au-dessus de $S$ sera le champ de
 vecteur
$$
U~:~exp(f(u)u)\mapsto \frac{d}{dt}\vert_{t=f(u)}exp(tu).
$$ 

Enfin si $\S$ est un graphe {\'e}tendu au-dessus de $S$, le {\it pied} d'un
point $y=exp(f(u)u)$ de $\S$ sera le point $x=\pi (u)$ de $S$, o{\`u} $\pi$ est la
projection du fibr{\'e} unitaire de $M$ sur $M$.

\subsubsection{Remarques}{\label{remdeflent}}
\begitem
\iti L'ouvert $U$ de la d{\'e}finition d'un graphe
{\'e}tendu
contenant l'adh{\'e}rence de $n(S)$ et $f$ {\'e}tant strictement positive,
$\nu(s)\not= n(s)$ pour tous les points $s$ de $\bo S$, o{\`u} $\nu$ est le
champ de vecteur normal ext{\'e}rieur {\`a} $\S$.

\itii Soit  $\S$ est un graphe {\'e}tendu au-dessus de
$S$. Si $\nu$ est le
champ de vecteur normal ext{\'e}rieur {\`a} $\S$ alors $\nu (\bo S)\subset int
(N_{S})$. En effet, d'apr{\`e}s \ref{deflent} (iii), le vecteur $n_\bo$
pointe vers l'int{\'e}rieur de $\S$, c'est-{\`a}-dire,
si $s$ est un point de $\bo S =\bo \S$, nous avons 
$$
\langle
n_{\bo}(s),\nu(s)\rangle\leq 0.
$$
Ensuite, $\nu(s)$ {\'e}tant
perpendiculaire {\`a} $T(\bo S)$ appartient au
plan engendr{\'e} par $n_{\bo }(s)$ et $n(s)$. En particulier, si
$\nu (s) \notin int(N_{S})$, nous avons n{\'e}cessairement $\nu (S)= -n(S)$ ,
mais ceci est impossible : deux surfaces strictement convexes ayant un
point d'intersection et en ce point deux normales oppos{\'e}es sont telles
que, au moins localement, leur intersection est r{\'e}duite {\`a} ce point.

\itiii Un graphe {\'e}tendu $\S$  au-dessus de $S$, se plonge naturellement dans le
bout complet de $S$ de telle sorte que $\partial\S\subset\partial S$.
\enditem
\subsubsection{Fonction inverse}
Introduisons une deni{\`e}re notion utile. Si $S$ est une surface lentille
pour $\S$ de fonction associ{\'e}e $\l$, la {\sl fonction inverse} est la fonction $\mu$ d{\'e}finie sur $\S$
par
$$
\mu (exp(\l(u)u))=\l(u).
$$

Nous avons alors le
\lem{fonc_inv}{Toute fonction inverse est 2-lipschitzienne}
\preu notons $\mu$ la fonction inverse.
Soient $u$ et $v$ deux {\'e}l{\'e}ments de $N_{S}$. En notant
$x=\pi (u)$ et $y=\pi (v)$, o{\`u} $\pi$ est la projection de $UM$ sur $M$,
$d_{\S}$ la distance riemanienne intrins{\'e}que de $\S$,
nous avons par convexit{\'e} locale de $S$ et gr{\^a}ce {\`a} notre hypoth{\`e}se de
courbure n{\'e}gative
$$
d(x,y)\leq d_{\S}(exp(\l(u)u),exp(\l(v)v).
$$
Enfin comme
$$
|\l(u)-\l(v)|\leq d(x,y) + d_{\S}(exp(\l(u)u),exp(\l(v)v),
$$

Nous en d{\'e}duisons bien le r{\'e}sultat.$\diamond$

\subsubsection{D{\'e}formations}
Gr{\^a}ce {\`a} la remarque \ref{remdeflent} (i), nous avons le lemme {\'e}vident de d{\'e}formation suivant

\lem{defolent}{Soit $S_{t}$, $t\in \oi$ une famille de surfaces immerg{\'e}es
localement con\-ve\-voi morse.dvi
xes, compactes {\`a} bord. Soit $\S_{t}$ une autre famille
de surfaces compactes, localement con\-ve\-xes, immerg{\'e}es et telle que $\bo
S_{t}=\bo \S_{t}$. Alors, si $\S_{0}$ est lentille pour $S_{0}$ alors
$\S_{t}$ est lentille pour $S_{t}$ pour tout $t$ dans un voisinage de 0.}

\subsection{G{\'e}om{\'e}trie born{\'e}e}\label{geoborn}

Rappelons quelques d{\'e}finitions.

\subsubsection{Convergence de vari{\'e}t{\'e}s point{\'e}es}
Une {\it vari{\'e}t{\'e} point{\'e}e} est une paire $(M,x)$ o{\`u} $M$ est une vari{\'e}t{\'e} et $x$
un point de $M$. Nous dirons que la suite de vari{\'e}t{\'e}s point{\'e}es $\{(
M_{n},x_{n})\inn$ {\'e}quip{\'e}es des m{\'e}triques $g_{n}$
{\it converge $\ci$ sur tout compact vers la vari{\'e}t{\'e} riemanienne
$(M_{\infty},x_\infty )$ } {\'e}quip{\'e}e de la m{\'e}trique $g_{\infty}$ s'il existe une
suite d'applications $\{f_{n}\inn$ (pas n{\'e}\-ces\-sai\-re\-ment continues)  d{\'e}finies de $M_{\infty}$ dans $M_{n}$,
envoyant $x_{\infty}$ dans $x_{n}$ telle
que pour tout $R$, alors
\begitem 
\iti pour $n$ suffisamment grand, la restriction  de $f_{n}$ {\`a} la boule
de rayon $R$ et de centre $x_{\infty}$ est $\ci$ et injective
\itii $f_{n}^{*}g_{n}$ converge $\ci$ sur tout compact vers $g_{\infty}$.
\enditem
\par

Nous dirons qu'une vari{\'e}t{\'e} $M$ est {\`a} {\sl g{\'e}om{\'e}trie born{\'e}e} si pour
toute suite de points $\{x_{n}\inn$ de $M$, la suite de vari{\'e}t{\'e}s
point{\'e}es $\{(M,x_{n})\inn$ poss{\`e}de une sous-suite convergente. 

Par abus de langage, nous dirons qu'une vari{\'e}t{\'e} $M$ d'Hadamard est {\`a}
{\sl g{\'e}o\-m{\'e}\-trie born{\'e}e} si quelque soit la suite de points $\{x_{n}\inn$ de $M$, la suite de vari{\'e}t{\'e}s
point{\'e}es $\{(M,x_{n})\inn$ poss{\`e}de une sous-suite convergente vers une
vari{\'e}t{\'e} {\`a} courbure strictement n{\'e}gative.

\subsubsection{Remarques}\label{poly}
\begitem
\iti Une vari{\'e}t{\'e} de Hadamard ayant un groupe discret cocompact d'i\-so\-m{\'e}\-tries
est {\`a} g{\'e}om{\'e}trie born{\'e}e.

\itii Si une vari{\'e}t{\'e} de Hadamard {\`a} courbure n{\'e}gative
est {\`a} g{\'e}o\-m{\'e}trie bor\-n{\'e}e, alors la
croissance des horosph{\`e}res est polynomiale. Remarquons en effet tout
d'abord que la courbure reste coinc{\'e}e entre deux bornes
strictement n{\'e}gatives. Le flot g{\'e}od{\'e}sique de $M$ va donc {\^e}tre d'Anosov. De
plus, par g{\'e}om{\'e}trie born{\'e}e le volume des boules de rayon 1 sur les
horosph{\`e}res va {\^e}tre uniform{\'e}ment born{\'e}. L'argument classique montrant
la croissance polynomiale des vari{\'e}t{\'e}s stables d'un flot d'Anosov nous
fournit le r{\'e}sultat.
\enditem

\subsubsection{Convergence de sous-vari{\'e}t{\'e}s immerg{\'e}es}

Nous nous int{\'e}resserons aux {\it vari{\'e}t{\'e}s immerg{\'e}es point{\'e}es}, c'est-{\`a}-dire aux quadruplets de la forme $Q=(N,x,f,M)$ o{\`u} $x$ est un point d'une
vari{\'e}t{\'e} $N$, $f$ une immersion de $N$ dans une vari{\'e}t{\'e} riemanienne  $M$. 
Nous dirons qu'une suite de vari{\'e}t{\'e}s immerg{\'e}es $\{Q_{n}=(N_{n},x_{n},f_{n},M_{n})\inn$ {\sl
converge $\ci$
sur tout compact}
vers une vari{\'e}t{\'e} immerg{\'e}e
$Q_{\infty}=(N_{\infty},x_{\infty},f_{\infty},M_{\infty})$ si 
\begitem
\iti La suite de vari{\'e}t{\'e}s point{\'e}es  $\{(M_{n},f_{n}(x_{n}))\inn$ converge
vers la vari{\'e}t{\'e} point{\'e}e
$(M_{\infty},f_{\infty}(x_{\infty}))$

\itii La suite de vari{\'e}t{\'e}s point{\'e}es $\{(N_{n},x_{n})\inn$, o{\`u} $N_{n}$ est
munie de la m{\'e}t\-ri\-que induite par $f_{n}$, converge vers  $\{(N_{\infty},x_{\infty})\}$, o{\`u} $N_{\infty}$ est
munie de la m{\'e}trique induite par $f_{\infty}$.

\itii $\{f_{n}\inn$ converge vers $f_{\infty}$ sur tout compact au sens
o{\`u} on l'imagine. 
\enditem
\par
Enfin, une vari{\'e}t{\'e} immerg{\'e}e $(N,x,f,M)$ sera dite {\sl {\`a} g{\'e}om{\'e}trie
born{\'e}e} si pour toute suite $\{x_{n}\inn$ la suite
$\{(N,x_{n},f,M)\inn$ admet une sous-suite convergente. 

Plus
g{\'e}n{\'e}ralement nous dirons qu'une suite de vari{\'e}t{\'e}s immerg{\'e}es 
$$
\{(N_{n},x_{n},f_{n},M_{n})\inn
$$
est {\it {\`a}
g{\'e}om{\'e}trie born{\'e}e} si pour toute suite de points $\{y_{n}\in N_{n}\inn$,
la suite 
$$
\{(N_{n},y_{n},f_{n},M_{n})\inn
$$ 
poss{\`e}de une sous-suite
convergente.

Lorsque la vari{\'e}t{\'e} ambiante $M$ est sous-entendue, nous abr{\'e}vierons souvent le quadruplet $(N,x,f,M)$
d{\'e}crivant une vari{\'e}t{\'e} immerg{\'e}e, en une paire $(N,x)$ et en confondant
ainsi de
mani{\`e}re abusive $N$ et son image.

\subsection{Surfaces horosph{\'e}riques, surface tubulaires}\label{deftubhor}
Dans ce paragraphe, $M$ sera toujours une vari{\'e}t{\'e} simplement connexe de
dimension $3$ {\`a} courbure strictement n{\'e}gative et {\`a} g{\'e}om{\'e}trie born{\'e}e.
\subsubsection{Surfaces horosph{\'e}riques, pseudo-horosph{\`e}res}
Par d{\'e}finition, une surface {\it horosph{\'e}rique} sera une surface convexe
compl{\`e}te plong{\'e}e $\S$, pour laquelle il existe une fonction de Busemann
$h$ v{\'e}rifiant les deux conditions suivantes
\begitem
\iti $h$ est major{\'e}e sur $\S$,
\itii le gradient de $h$ est partout transverse {\`a} $\S$.
\enditem
Une suite de surfaces immerg{\'e}es $\{(S_{n},x_{n})\inn$ sera dite {\it de
type ho\-ros\-ph{\'e}\-ri\-que} s'il existe une suite de points $\{y_{n}\in
S_{n}\inn$
telle que $\{(S_{n},y_{n})\inn$ poss{\`e}de une sous-suite qui converge vers
une surface horopsh{\'e}rique.

Une
surface immerg{\'e}e  $(\S, x)$ sera dite {\it horosph{\'e}rique {\`a} l'infini}
s'il existe une suite de points $\{x_{n}\inn$ de $\S$ telle que
$\{(\S,x_{n})\inn$ converge vers une surface horosph{\'e}rique.

Nous dirons  qu'une surface convexe compl{\`e}te plong{\'e}e $S$ dans $M$ est
une {\sl pseudo-horosph{\`e}re} si l'application de Gauss-Minkowski d{\'e}finie
de $S$ dans $\bo_{\infty}M$ par, 
$$
x\mapsto exp(+\infty n(x))
$$ 
o{\`u} $n$ est
le champ de vecteur normal ext{\'e}rieur, est une bijection sur
$\bo_{\infty}M$ priv{\'e} d'un point. 

Une surface horosph{\'e}rique est une
pseudo-horosph{\`e}re et r{\'e}ciproquement nous avons le 
\lem{psudo}{Dans une vari{\'e}t{\'e} {\`a} g{\'e}om{\'e}trie born{\'e}e une pseudo-horosph{\`e}re
{\`a} g{\'e}o\-m{\'e}\-trie born{\'e}e est horosph{\'e}rique {\`a} l'infini} 
\preu soit $y\in\bo_{\infty}M$ le point {\'e}vit{\'e} par l'application de Gauss-Minkowski de la pseudo-horosph{\`e}re $S$. Soit $B$ l'ensemble convexe bord{\'e}
par $S$.  Construisons une suite de convexes $\{C_{n}\inn$ ayant les
propri{\'e}t{\'e}s suivantes
\begitem
\iti $\bo_{\infty} C_{n}$ est un voisinage de $y$,
\itii $\{C_{n}\inn$ converge au sens de Haussdorff dans le compactifi{\'e}
de $M$ vers $y$.
\enditem

Soit maintenant $d_{n}$ la fonction distance {\`a} $C_{n}\cap B$. D'apr{\`e}s
(i) sur $S$,
$d_{n}$ atteint son maximum en un point $x_{n}$. D'apr{\`e}s (ii), 
$d_{n}(x_{n})$ tend vers $+\infty$.

Consid{\'e}rons la suite de surfaces immerg{\'e}es $\{(S,x_{n})\inn$ immerg{\'e}es
dans la suite $\{(M,x_{n})\inn$. Par nos hypoth{\`e}ses,  $\{(S,x_{n})\inn$ va converger
vers une surface  $(S_{\infty},x_{\infty})$ immerg{\'e}e dans une vari{\'e}t{\'e}
$(M_{\infty},x_{\infty})$. Enfin la suite de fonctions
$d_{n}-d_{n}(x_{n})$ va converger vers une fonction de Buseman sur
$M_{\infty}$, major{\'e}e sur $S$ et dont le gradient est dirig{\'e} vers
l'ext{\'e}rieur de $S_{\infty}$. La surface $S_{\infty}$ est donc de type
horosph{\'e}rique, et $S$ est donc bien horosph{\'e}rique {\`a} l'infini.$\diamond$

\subsubsection{Surfaces tubulaires, tubes}
Par d{\'e}finition, une surface {\it tubulaire} sera une surface localement convexe
compl{\`e}te plong{\'e}e $\S$, pour laquelle il existe une g{\'e}od{\'e}sique telle que
si $d$ est la fonction distance {\`a} cette g{\'e}od{\'e}sique,
$d$ v{\'e}rifie les deux conditions suivantes
\begitem
\iti $d$ est born{\'e}e sur $\S$,
\itii le gradient $u$ de $d$ est dirig{\'e} vers l'ext{\'e}rieur  de $\S$,
c'est-{\`a}-dire $\langle u,n\rangle \geq 0$, o{\`u} $n$ est le champ de vecteur
normal ext{\'e}rieur.
\enditem
Il est int{\'e}ressant de remarquer que la surface {\'e}tant convexe, et les
lignes de gradient de $d$ {\'e}tant des g{\'e}od{\'e}siques, la condition (iii) est
{\'e}quivalente {\`a} la condition o{\`u} on impose  $\langle u,n\rangle > 0$.

Une suite de surface immerg{\'e}e sera dite {\it de type tubulaire} si elle
poss{\`e}de une sous-suite qui converge vers une surface tubulaire.

Une
surface immerg{\'e}e  $(\S, x)$ sera dite {\it tubulaire {\`a} l'infini}
s'il existe une suite de points $\{x_{n}\inn$ de $\S$ telle que
$\{(\S,x_{n})\inn$ converge vers une surface tubulaire.

Enfin une d{\'e}finition connexe nous sera souvent utile par la suite :
le {\it tube d'une g{\'e}od{\'e}sique $\gamma$} est d{\'e}fini par
$$
N(\gamma)=\{u\in UM / \langle u,\dgt \rangle  =0\}.
$$

\subsubsection{Remarque}
Il est utile de remarquer qu'{\^e}tre de type horosph{\'e}rique, tubulaire,
ho\-ro\-sph{\'e}\-ri\-que ou tubulaire  {\`a} l'infini est une propri{\'e}t{\'e} qui ne d{\'e}pend que de
$\{S_{n}\inn$ et $S$ r{\'e}ciproquement et ne fait pas intervenir les points
choisis sur ces surfaces

\subsection{Principe du maximum g{\'e}om{\'e}trique et applications}

Nous avons l'{\'e}nonc{\'e} {\'e}vident suivant que nous appelerons du nom pompeux
de {\it principe du maximum g{\'e}om{\'e}trique}.

\lem{max}{Soit $S_{1}$ et $S_{2}$ deux surfaces convexes tangentes en un
point $x$ et telle que $S_{1}$ est {\`a}
l'int{\'e}rieur de $S_{2}$ au voisinage de $x$, alors $k_{1}(x)\geq
k_{2}(x)$, o{\`u} $k_{i}$ d{\'e}signe la courbure extrins{\`e}que  de $S_{i}$.}

Nous allons en tirer quelques cons{\'e}quences. Pour cela supposons que la
courbure de $M$ est plus petite que $-c<0$ et soit $k$ tel
$0<k<c$. Nous avons 

\pro{lentboul}{ Si $S$ est une $k$-surface compacte telle que $\partial S$ est
inclus dans une boule $B$, alors $S$ toute enti{\`e}re est incluse dans
cette boule.}
\preu en effet avec nos hypoth{\`e}ses, toute sph{\`e}re a une courbure plus
grande que $c$.

Rappelons en rapidement la d{\'e}monstration en utilisant \ref{varinfcal} :  si nous notons $\l^{r}_{i}$ les valeurs
propres de l'op{\'e}rateur deuxi{\`e}me forme fondamentale de la sph{\`e}re de rayon
$r$, $S^{r}$, associ{\'e}es
aux vecteurs propres $e_{i}^{r}$, ainsi que $k^{r}_{i}$ les courbures
sectionelles des plans perpendiculaires {\`a} $S^{r}$ passant par $e^{i}$,
d'apr{\`e}s \ref{sign}, la courbure $\k_{r}$ de la sph{\`e}re de rayon
$r$ v{\'e}rifie l'{\'e}quation
$$
\frac{d\k_{r}}{dr}=-\k_{r}( \l^{r}_{2}(1+\frac{k^{r}_{1}}{\k_{r}})+\l^{r}_{1}(1+\frac{k^{r}_{2}}{\k_{r}})) 
$$
La propri{\'e}t{\'e} annonc{\'e}e d{\'e}coule imm{\'e}diatement de l'in{\'e}galit{\'e} :
$$
{\frac{d\k_{r}}{dr}}\geq (\k_{r}-c)(-\l^{r}_{1}-\l^{r}_{2}).
$$

Nous pouvons conclure par l'absurde maintenant. Si $S$ n'est pas incluse dans
$B$,  on peut construire une  sph{\`e}re {\`a} laquelle $S$ est tangente
int{\'e}rieurement ce qui contredit notre principe du maximum. $\diamond$

Le m{\^e}me raisonnement montre

\pro{lenthoros}{Une $k$-surface ne peut-{\^e}tre horosph{\'e}rique.}

\subsection{Domination}

Soit $S$ une surface lentille pour $\S$, de fonction associ{\'e}e $f$ 
d{\'e}finie sur l'ad\-h{\'e}\-ren\-ce d'un ouvert $U$ du normal {\'e}tendu $N_{S}$ {\`a} $S$. Nous dirons
que {\sl $S$ domine la surface $S_{1}$},  s'il existe un ouvert $V$ de $U$
tel que ${\bar V}\subset U$, o{\`u} ${\bar V}$ est l'adh{\'e}rence de $V$ dans
$N_{S}$, une fonction positive  ou nulle  $g$ d{\'e}finie sur ${\bar V}$,
inf{\'e}rieure ou {\'e}gale {\`a}
$f$ sur $V$ et {\'e}gale   {\`a}
$f$ sur ${\bar V}\setminus V$, telle que $S_{1}$ est le graphe de $g$
sur $V$. Dans le cas o{\`u} la fonction $g$ est strictement positive, nous
dirons que   {\sl $S$ domine strictement la surface $S_{1}$}

Un lemme facile et utile est le
\lem{domi}{Soient $S$ une surface lentille pour $\S$, $\S_{t}$
une famille  d'ouverts {\`a} bord de $\S$ pour $t\in [0,1]$ et $S_{t}$ une
famille de surfaces lentilles pour $\S_{t}$, alors
\begitem
\iti l'ensemble des $t$ tels que $S$ domine strictement $S_{t}$ est ouvert,  
\itii si $\forall t\in ]0,1]$, $S$ domine strictement $S_{t}$ alors,
soit $S$ domine strictement
$S_{0}$, soit $S$ est tangente int{\'e}rieurement {\`a} $S_{0}$ en un point
int{\'e}rieur {\`a} $S$.
\enditem}

\section{Variations infinit{\'e}simales et d{\'e}formations}

Le but essentiel de cette section est le

\lem{defoloc}{Soit $k$ un r{\'e}el positif tel que   $M$ soit {\`a} courbure
strictement plus petite que $-k$. Soit $\S$ une k-surface compacte {\`a}
bord $\dSi$ et soit $c_{t}$, une famille continue de d{\'e}formations  de $\dSi$
telle que $c_{0}=\dSi$, soit $k(t)$ une fonction $C^{\infty}$ de $t$ telle
que $k(0)=k$. Il existe alors une famille unique
de $k(t)$-surfaces immerg{\'e}es $\S^{t}$, d{\'e}finie au voisinage de 0,  telle que
$\S^{0}=\S$ et  $\bo\S^{t}=c^{t}$.

De plus, s'il existe une famille immerg{\'e}e de surfaces $S^{t}$ v{\'e}rifiant  $c_{t}=\bo S^{t}$, telle
que $\S$ est lentille pour $S^{0}$, alors $\S^{t}$ est lentille pour
$S^{t}$ au voisinage de 0. 

Enfin, si $S^{t}\subset S^{s}$ pour $t\geq s$ et $k(t)$ est 
croissante alors $\S$ domine 
$\S^{t}$, toujours pour $t$ au voisinage de $0$. Si de plus $S^t
\not=S^{0}$ ou $k(t)$ est strictement croissante, alors $\S$ domine strictement
$\S^t$ pour $t$ au voisinage de z{\'e}ro. 
}

Le coeur de la d{\'e}monstration de cette proposition est une proposition
de d{\'e}formation infinit{\'e}simale que nous allons maintenant pr{\'e}senter.

Soit  $S\subset M$ une surface
compacte  immerg{\'e}e {\`a} bord. Notons $C_{0}^{\infty}(S)$ l'espace des fonctions $\ci$
d{\'e}finies sur $S$ et nulles sur le bord.

A toute $f\in\ci (S)$, on peut associer une
{\it variation de surfaces} d{\'e}finie
pour $t$ suffisamment petit par les immersions $s^{f}_{t}$
\eqalign{
S&\fl M\cr 
x&\longmapsto \exp(tfn(x))=s_{t}^{f}(x),
}

o{\`u} $n$ d{\'e}signe le champ de vecteur normal ext{\'e}rieur {\`a} $S$.

D{\'e}signons  par $k^{f}_{t}(x)$ la courbure extrins{\`e}que {\`a} la surface
$s^{f}_{t}(S)$ au point $s_{t}^{f}(x)$,
l'op{\'e}rateur $L$ {\it de variation infinit{\'e}simale de courbure
extrins{\`e}que} est alors
\eqalign{
C^{\infty}_{0}(S)&\fl C^{\infty}(S)\cr
f& \longmapsto  L(f)=\frac{d}{dt}\mid_{t=0}{k^{f}_{t}}.
}

Nous montrerons

\pro{varinf}{Soit $k$ un r{\'e}el positif. Si $M$ est {\`a} courbure  strictement plus petite que $-k$
et $S$ est {\`a} courbure extrins{\`e}que strictement comprise entre 0 et $k$,
alors l'op{\'e}rateur $L$ est elliptique et inversible}

Dans cette section, nous expliciterons tout d'abord l'op{\'e}rateur $L$, puis
d{\'e}\-mon\-tre\-rons \ref{varinf} en utilisant le principe du maximum. Enfin
nous montrerons  \ref{defoloc}
\subsection{Explicitation de l'op{\'e}rateur $L$}

Nous allons calculer explicitement l'op{\'e}rateur $L$ pour une surface
quelconque $S$. Notons pour cela $\n$ la connexion de Levi-Civita de $M$,
$R$ son tenseur de courbure, $n$ le vecteur normal {\`a} $S$, $W$ l'endomorphisme de $TS$ d{\'e}fini par
$W(u)=R(n,u)n$,  $\k$ la courbure extrins{\`e}que de $S$, $B$
l'op{\'e}rateur deuxi{\`e}me forme fondamentale de $S$ d{\'e}fini par $B(u)=\n_{u}n$
et enfin ${\rm Hess}(f)$ la hessienne de $f$. 

Montrons alors
\pro{varinfcal}{Nous avons 
$$
L(f)=\k ( -\trace  ({\rm Hess}(f)\circ B^{-1})+f \trace  (W\circ B^{-1})- f\trace  (B))
$$
}

\preu
donnons tout d'abord le cadre de ce calcul. Il nous faut consid{\'e}rer la famille d'immersions
$s^{f}_{t}$ comme une application $s^{f}$ de $S\times\mathbb R$ dans
$M$. Nous identifions $S$ {\`a} $S\times\{0\}$. Sur le fibr{\'e} $E$ induit de $TM$ par $s^{f}$, nous noterons par abus
de langage $\n$ la connexion induite de la connexion de Levi-Civita de $M$. Nous pouvons maintenant voir $Ts^{f}=F$, comme une section de
$TS^{*}\otimes E$, et nous avons bien s{\^u}r  $d^{\n}F=0$ et $F(\frac{\partial}{\partial t})=fn$ le long de $S$.

Tout champ de vecteur $u$ sur $S$ donne canoniquement naissance {\`a} un
champ de vecteur not{\'e} {\'e}galement $u$ sur $S\times \mathbb R$ qui commute avec le champ
$\frac{\partial}{\partial t}$. Consid{\'e}rons {\'e}galement $n_{t}$ le champ de vecteur normal
{\`a} $S_{t}=s_{t}^{f}(S)$, et $n$ la section de $E$ qui s'en d{\'e}duit.

Par abus de notation, si
$v$ est une section de ce fibr{\'e} $E$ nous noterons 
$$
\frac{d}{dt}\mid_{t=0}{v}(x) =(\n_{\frac{\partial}{\partial t}}v )(x,0)
$$

Consid{\'e}rons  $A$,  la section de  $TS^{*}\otimes E$  d{\'e}finie par
$$
A(u)=\n_{u}n.
$$

Un premier calcul donne
$$
\frac{d}{dt}\mid_{t=0}{F(u)}=\n_{u}(fn)=df(u).n +fA(u).
$$
En notant $\langle ,\rangle$ la m{\'e}trique de $M$ et celle qui s'en d{\'e}duit
sur $E$, nous avons
$$
0=\frac{d}{dt}\mid_{t=0}{\langle n,{F(u)}\rangle}=\langle \frac{d}{dt}\mid_{t=0}{n},F(u)\rangle+\langle n,\frac{d}{dt}\mid_{t=0}{F(u)}\rangle.
$$
Nous en tirons facilement 
$$
\frac{d}{dt}\mid_{t=0}{}=F(-\n f).
$$
Nous obtenons ainsi
\eqalign{
\frac{d}{dt}\mid_{t=0}{A(u)}&=\n_{\frac{\partial}{\bo t}}\n_{u}n\cr
&=R(\frac{\partial}{\bo t},u)n +\n_{u}\n_{\frac{\partial}{\bo t}}n\cr
&=fR(n,u)n -\n_{u}\n f.
}

Consid{\'e}rons $g$ la m{\'e}trique induite sur $TS$ par $F$, c'est-{\`a}-dire
d{\'e}finie par 
$$
g(u,u)=\langle F(u),F(u)\rangle
$$ 
et $B$ l'endomorphisme de $TS$ d{\'e}fini par
$$
g(B(u),u)=\langle A(u),u\rangle.
$$

Rappelons que nous voulons calculer
$$
L(f)=\frac{d}{dt}\mid_{t=0}{~\det (B)}.
$$

Pour cela utilisons le fait que 
$$
\frac{d}{dt}\mid_{t=0}{g(u,u)}=2f\langle A(u),F(u)\rangle=2fg(B(u),u).
$$
En d{\'e}rivant l'equation $g(B(u),u)=\langle A(u),u\rangle$, nous obtenons

\eqalign{
&g(\frac{d}{dt}\mid_{t=0}{B(u)},u)=\langle \frac{d}{dt}\mid_{t=0}{A(u)},u\rangle-fg(B(u),B(u))\\
&=f\langle R(n,u)n,u\rangle-\langle \n_{u}\n f,u\rangle-fg(B(u),B(u));
}

Et donc 
$$
\frac{d}{dt}\mid_{t=0}{B}=fR(n,u)n-{\rm Hess}(f)-fB^{2}.
$$
La formule classique 
$$
\frac{d}{dt}\mid_{t=0}{~{\rm log}(\det (B))}=\trace  (\frac{d}{dt}\mid_{t=0}{B}\circ B^{-1}),
$$
nous donne la proposition.$\diamond$

Citons un corollaire utile des formules d{\'e}montr{\'e}es dans ce paragraphe

\cor{equic}{ Supposons que $M$ ait une courbure plus petite que $-c$. Les
sph{\`e}res et les horosph{\`e}res de $M$ ont alors  une courbure plus
grande que $c$. De plus, si $S$ est une surface convexe, pour tout
$k<c$, il existe $R$ ind{\'e}pendant de $S$ tel que la surface $S_{R}=\exp(Rn(S))$ ait ses
courbures principales plus grandes que $k^{1/2}$ et en particulier sa 
courbure plus grande que $k$}
 
\subsection{D{\'e}monstration de la proposition \ref{varinf}}

Si $\k$ est strictement positif, $B$ est d{\'e}fini positif. En particulier,
$L$ est elliptique d'indice nul.

Pour conclure, il nous suffit de d{\'e}montrer que $L$ est injectif. La
proposition suivante et une application standard du principe du maximum
permettent de conclure :

\pro{sign}{Le terme de degr{\'e} z{\'e}ro de la formule de \ref{varinfcal} est
strictement positif :
$$
J= \trace  (W\circ B^{-1})- \trace  (B)>0.
$$
}
\preu utilisons $(e_{1},e_{2})$ une base de vecteurs propres de $B$ associ{\'e}e
aux valeurs propres $\l_{1}$ et $\l_{2}$. Soit $k_{i}$ la courbure du
plan engendr{\'e} par $n$ et $e_{i}$.

Nous obtenons
\eqalign{
J=&-\frac{k_{1}}{\l_{1}}-\frac{k_{2}}{\l_{2}}-\l_{1}-\l_{2}\\
=1 -\l_{2}(1+\frac{k_{1}}{\k})-\l_{1}(1+\frac{k_{2}}{\k})>0.
}

$\diamond$

\subsection{D{\'e}monstration du lemme \ref{defoloc}} Pour d{\'e}montrer la premi{\`e}re partie de cette proposition, nous
allons pro\-c{\'e}\-der par {\'e}tapes.

Nous r{\'e}solvons tout d'abord le probl{\`e}me infinit{\'e}simal. Une variation infinit{\'e}simale de
${\S}^{t}$, est 
un champ de vecteur $\zeta$ le long de ${\S}^{t}$ que nous pouvons {\'e}crire
sous la forme
$$
\zeta=fn+u,
$$
o{\`u} $~u\in TS$,  $~f\in \ci ({\S}^{t})$, et  v{\'e}rifiant la condition au
bord 
$$
\forall x\in\bo \S,~~~\zeta (x)=\frac{d}{dt}{{c_{t}(x)}}.
$$

Ici, nous voyons abusivement $c_{t}$ comme une famille d'immersions de
$\dSi$ dans $M$.

Pour une surface de courbure extrins{\`e}que $\k$, la variation infinit{\'e}simale de courbure
extrins{\`e}que associ{\'e}e {\`a} une telle  variation $\zeta$ est
$$
L(f)+d\k (u),
$$
o{\`u} $L$ est l'op{\'e}rateur de \ref{varinf}. 

Pour r{\'e}soudre infinit{\'e}simalement notre probl{\`e}me, nous devons donc
montrer qu'il existe une unique fonction $f$ telle que
$$
L(f)=0,
$$
avec la condition au bord
$$\
~\forall x\in\bo S,~f(x)=\langle \frac{d}{dt}{c_{t}(x)},n\rangle .
$$

Ceci d{\'e}coule  de \ref{varinf}.

Nous pouvons maintenant utiliser le th{\'e}or{\`e}me d'inversion locale pour
les op{\'e}rateurs elliptiques pour
r{\'e}soudre notre probl{\`e}me localement, \cad cons\-trui\-re au voisinage de $0$
une unique famille $t\mapsto\S^{t}$, continue en $t$, de $k$-surfaces
v{\'e}rifiant $\S^{0}=\S$ et $\bo\S^{t}=c_{t}$.

Mettons nous maintenant dans le cadre de la deuxi{\`e}me partie de
\ref{defoloc}. Notre lemme \ref{defolent} nous permet de montrer
que $\S^{t}$ est lentille pour $S^{t}$

Enfin, $S^{t}\subset S^{0}$ pour $t\geq 0$, la proposition
\ref{sign}, assure que la fonction $f$ d{\'e}crivant la variation
infinit{\'e}simale de $\S$ est positive, puisque dans le cas d'une surface lentille nous
avons $\langle \frac{d}{dt}{c^{t}(x)},n\rangle  >0$. Ainsi $\S$ domine $\S^{t}$
pour $t$ positif et suffisamment petit.

%%%%%%%%%%%%%%%%%%%%%%%%%%%%%%%%%%%%%%%%%%%%%%%%%%%%%%%%%%%%%%%%

\section{Th{\'e}or{\`e}me de compacit{\'e}}
Nous allons dans cette section {\'e}noncer et d{\'e}montrer  le th{\'e}or{\`e}me de compacit{\'e} qui nous
sera utile par la suite.

\theo{compalent}{ Soit $M$ une vari{\'e}t{\'e} d'Hadamard {\`a} g{\'e}om{\'e}trie born{\'e}e et {\`a}
courbure plus petite que $-c < 0$, soit $\{(S_{n},x_{n})\inn$ une suite de surfaces convexes
immerg{\'e}es {\`a} courbure extrins{\`e}que plus grande que $c$,  convergeant vers $(S_{\infty},x_{\infty})$. Nous supposerons de plus cette suite est {\`a} g{\'e}om{\'e}trie
born{\'e}e et n'est pas de type
ho\-ros\-ph{\'e}\-ri\-que. 

Soit enfin $\Sinn$ une suite de $k_{n}$-surfaces lentilles pour $\Snn$, o{\`u}
$k_{n}$ converge vers $k\in ]0,c[$, et notons  $y_{n}\in\S_{n}$ le pied de
$x_{n}$.

Nous avons alors les trois possibilit{\'e}s suivantes
\begitem
\iti $S_{\infty}$ n'est pas  tubulaire {\`a} l'infini et alors, apr{\`e}s
extraction d'une sous-suite  $\{(\S_{n},y_{n})\inn$ converge vers une
$k$-surface lentille non d{\'e}g{\'e}n{\'e}r{\'e}e pour $S_{\infty}$,

\itii $S_{\infty}$ n'est pas  tubulaire, tout en {\'e}tant
tubulaire {\`a} l'infini, alors apr{\`e}s extraction d'une sous-suite,  $\{(\S_{n},y_{n})\inn$ converge vers une
$k$-surface lentille {\'e}ventuellement  d{\'e}g{\'e}n{\'e}r{\'e}e pour $S_{\infty}$,

\itiii $S_{\infty}$ est tubulaire pour une g{\'e}od{\'e}sique $\gamma$
 et alors $\{(n(\S_{n}),n(y_{n}))\inn$ converge apr{\`e}s
 extraction d'une sous-suite vers le tube de $\gamma$.
\enditem
}

Nous allons tout d'abord  rappeler les r{\'e}sultats principaux de \ma et leurs
cons{\'e}quences sur ce que nous appelerons le probl{\`e}me de Dirichlet pour les
$k$-surfaces. Ensuite, nous appliquerons ces r{\'e}sultats pour d{\'e}montrer
\ref{compalent}.

\subsection{Rappels sur les probl{\`e}mes de Monge-Amp{\`e}re}

\subsubsection{D{\'e}finitions} Dans \ma, nous avions {\'e}tudi{\'e} une classe de probl{\`e}mes que nous avions
appel{\'e}e de Monge-Amp{\`e}re et dont les $k$-surfaces forment une classe
d'exemple. De plus, nous y avions d{\'e}fini une classe de probl{\`e}me {\`a} bord
convexe correspondant, dans notre contexte,  {\`a} la situation suivante.

Un {\it  probl{\`e}me de Dirichlet} pour les $k$-surfaces est la donn{\'e}e d'un
quadruplet $(M,S,c,x)$, o{\`u} $c$ est une courbe plong{\'e}e compl{\'e}te, trac{\'e}e sur une
surface plong{\'e}e $S$ localement convexe, incluse dans une vari{\'e}t{\'e}
riemanienne $M$. Le point $x$ sert en particulier {\`a}
donner un sens {\`a} la notion de convergence pour une suite de probl{\`e}mes de
Dirichlet.

Une {\it solution du probl{\`e}me} de Dirichlet est une $k$-surface $\S$
compl{\`e}te connexe immerg{\'e}e passant par $x$ telle que  $\bo\S \subset c$ et $\S$ est int{\'e}rieure {\`a} $S$ le long de
$c$.

Une {\it solution d{\'e}g{\'e}n{\'e}r{\'e}e} est une $k$-surface $\S$ dont le bord $\bo\S$
est complet et inclus dans  $c$, telle $\S$ est int{\'e}rieure {\`a} $S$ le long
de $c$ et telle qu'enfin $n(\S )\subset UM$ soit  compl{\`e}te, sans que
$\S$ le soit. Remarquons que dans ce cas, toute surface {\'e}quidistante de
$\S$ est compl{\`e}te.

Pour des raisons techniques, nous sommes oblig{\'e}s d'introduire la
d{\'e}finition suivante :  si $U$ est un ouvert de $M$ contenant $x$, le probl{\`e}me $(U,S\cap
U, c\cap U, x)$  sera appel{\'e} {\it probl{\`e}me
restreint} {\`a} $U$  et la {\it restriction {\`a} $U$ d'une solution $\S$} d{\'e}finie
sur
$M$ est la composante connexe de $\S\cap U$ contenant $x$. 

\subsubsection{Remarques}

\begitem

\iti un tube est une surface rideau au sens de \ma ;

\itii une $k$-surface lentille $\S$ pour $S$ est un cas particulier de solution du
prob\-l{\`e}\-me de Dirichlet d{\'e}fini par $S$ et $\dS$ ;

\itiii Si la courbure moyenne d'une solution est born{\'e}e, alors la
solution est non d{\'e}g{\'e}n{\'e}r{\'e}e ;

\itiv Si $S$ est une $k$-surface d{\'e}g{\'e}n{\'e}r{\'e}e, alors toute surface
{\'e}quidistante est com\-pl{\`e}\-te.

\itv enfin, il se peut que $S$, $c$ ou $\partial \S$ soient vide.
\enditem
\subsubsection{Compacit{\'e}} Pour {\'e}noncer le r{\'e}sultat principal de \ma, il nous faut introduire une
notation : pour tout $\e >0$ et probl{\`e}me de Dirichlet $\D$,  $~\D^{\e}$
d{\'e}signera 
 le probl{\`e}me restreint {\`a} la boule ouverte de centre $x$ de rayon
$\e$. De m{\^e}me, $\S^{\e}$ d{\'e}signera la restriction d'une solution $\S$ {\`a}
cette boule.

Nous avons alors le 

\theo{comrap}{ Soit $\{ \D_{n}=(M_{n},S_{n},c_{n},x_{n})\inn$ une suite de
probl{\`e}mes de Dirichlet convergeant $\ci$ sur tout compact vers un probl{\`e}me
de Dirichlet 
$$
\D_{\infty}=(M_{\infty},S_{\infty},c_{\infty},x_{\infty})
$$ 
et soit $\{\S_{n}\inn$
une suite de $k_{n}$-solutions, {\'e}ventuellement d{\'e}\-g{\'e}\-n{\'e}\-r{\'e}e de
$\{\D_{n}\inn$ o{\`u} $\{k_{n}\inn$ converge vers $k\in ]0,c[$. 
Il existe alors $\epsilon$, tel que   nous ayions
l'alternative suivante, apr{\`e}s extraction d'une sous-suite :

(a) soit $\{\S_{n}^{\e}\inn$ converge vers une $k$-solution, {\'e}ventuellement
d{\'e}gen{\'e}r{\'e}e, du probl{\`e}me de Dirichlet $\D_{\infty}^{\e}$ ;

(b) soit $\{n(\S_{n}^{\e})\inn$ converge vers un tube. Dans ce cas, {\`a} partir
d'un certain rang $\partial \S_{n}^{\e}=\emptyset$. }

Enfin, nous tirons de ce th{\'e}or{\`e}me les corollaires suivant, particulier {\`a}
la courbure n{\'e}gative 

\cor{corcomrap}{Soit $\{ \D_{n}=(M_{n},S_{n},c_{n},x_{n})\inn$ une suite de
probl{\`e}mes de Dirichlet convergeant $\ci$ sur tout compact vers un probl{\`e}me
de Dirichlet 
$$
\D_{\infty}=(M_{\infty},S_{\infty},c_{\infty},x_{\infty})
$$ 
et soit $\{\S_{n}\inn$
une suite de $k_{n}$-solutions, {\'e}ventuellement d{\'e}\-g{\'e}\-n{\'e}\-r{\'e}e de
$\{\D_{n}\inn$, o{\`u} $\{k_{n}\inn$ converge vers $k\in ]0,c[$. Alors la suite de surfaces
$$
W_{n}=\{exp(n(s)),~~s\in \S_{n}\}
$$ 
point{\'e}e en $exp(n(x_{n}))$ converge $\ci$ apr{\`e}s extraction d'une
sous suite.}
\preu ceci provient de ce que la restriction de l'exponentielle {\`a} un
tube est une immersion {\`a} valeur dans $M$.$\diamond$

\cor{graphe}{Soit $\{(f_{n}, S)\inn$ une suite de surfaces immerg{\'e}es
localement con\-ve\-xes
convergeant $\ci$ sur tout compact vers une surface immerg{\'e}e
$(f_{0},S)$ et soit $\{l_{n}\inn$ une suite de fonctions positives d{\'e}finies
sur $S$ dont les graphes sont des $k$-surfaces. On suppose que pour
tout $y\in S$, la suite $\{l_{n}(y)\inn$ est born{\'e}e, alors apr{\`e}s
extraction d'une sous-suite, la suite de
fonctions $\{l_{n}\inn$ converge $\ci$ sur tout compact de $S$ vers une
fonction dont le graphe est une $k$-surface}

\preu il suffit pour cela d'appliquer le th{\'e}or{\`e}me de compacit{\'e}
\ref{comrap} pour la suite de probl{\`e}mes de Dirichlet ({\`a} bord vide)
d{\'e}finie dans les bouts des surfaces $f_{n}(S)$.$\diamond$
\subsection{Une premi{\`e}re majoration} Nous nous donnons donc {\`a} partir de maintenant une suite $\Sinn$ de
$k_{n}$-surfaces lentilles pour $\Snn$. On suppose de plus que $\Snn$ converge
vers $S_{\infty}$ qui est {\`a} g{\'e}om{\'e}trie born{\'e}e. 

 Soit  $\{\l\inn$ la suite de
fonctions associ{\'e}es d{\'e}finies sur $\Nn$ 
l'ensemble des vecteurs normaux {\'e}tendus {\`a} $\S_{n}$.

Nous voulons tout d'abord montrer

\pro{premaj}{ Si la suite $\Snn$ n'est pas de type horosph{\'e}rique, alors la suite de
fonctions $\{\l\inn$ est born{\'e}e }

\preu raisonnons par l'aburde et supposons que  $\L_{n}$, le maximum de
$\l_{n}$, tende vers l'infini. Il existe alors  deux suites
de points  $\{w_{n}\in S_{n}\inn$ et $\{v_{n}\in\Nn\inn$, o{\`u} $\Nn$ d{\'e}signe
l'ensemble des vecteurs normaux {\'e}tendus {\`a} $\S_{n}$ telles  que 
\begitem
\iti $w_{n}=exp(\l_{n}(v_{n})v_{n})$ ;

\itii la suite $\{\l_{n}(v_{n})-\L_{n})\inn$ tende vers z{\'e}ro.
\enditem
Notons, pour tout $R$, $S^{R}_{n}$, la composante connexe de
l'intersection de $S_{n}$ avec la boule ouverte de $M$ de centre $w_{n}$ et de rayon $R$. Notons $U_{n}$ le champ
de vecteur focal le long de $S_{n}$.

La  suite $\Snn$ {\'e}tant {\`a} g{\'e}om{\'e}trie born{\'e}e, en extrayant au besoin une
sous-suite, nous pouvons supposer que $\{S^{R}_{n},w_{n}\inn$
converge.

Les fonction inverses $\mu_{n}$ d{\'e}finies sur $S_{n}$ {\'e}tant
2-lipschtziennes par le lemme \ref{fonc_inv}, nous pouvons extraire une sous-suite telle que la suite de fonctions $\{\mu_{n}-\L_{n}\inn$  converge uniform{\'e}ment vers une
fonction 2-lipschtzienne n{\'e}gative $g$. En particulier $\{\mu_{n}\inn$ converge
uniform{\'e}ment vers l'infini
sur $S^{R}_{n}$.

Nous en tirons deux cons{\'e}quences :

(a) pour tout $L$, il existe $n_{0}$, tel que si $n\geq n_{0}$ alors
$d_{i}(y_{n},\bo S^{R}_{n})\geq L$ ( o{\`u} $d_{i}$ d{\'e}signe la distance
riemanienne  dans $S_{n}$ ) 
puisque la restriction de $\mu_{n}$ au bord de $\S_{n}$ est nulle 

(b) {\`a} partir d'un certain rang,  pour tout $0<t<2R$ et $z,~w\in
S^{R}_{n}$, $d(exp (-tU_{n}(z),-tU_{n}(w))\leq d(z,w)$.

Rappelons que dans une vari{\'e}t{\'e} {\`a} courbure n{\'e}gative, une hypersurface
sans bord localement convexe et compl{\`e}te respectivement {\`a} une boule, est
plong{\'e}e et borde un convexe.

En utilisant (a) et ce rappel, nous en d{\'e}duisons 
que $\{S_{n},v_{n}\inn$
converge sur tout compact vers une hypersurface plong{\'e}e globalement
convexe $(H,v)$.

Gr{\^a}ce {\`a} (b), en extrayant une sous-suite de telle sorte que
$\{U_{n}(v_{n})\inn$ converge vers un vecteur $U$, nous en d{\'e}duisons que $\{U_{n}\inn$
converge sur tout compact vers le champ de gradient de la fonction de
Buseman $f$ associ{\'e}e {\`a} $exp(-\infty U)$ qui va se retrouver transverse {\`a}
$H$.
 
Par construction $f\vert_{H} -f(v)=g-g(v)$, et nous en d{\'e}duisons que
$f\mid_{H}-f(v)\leq g(v)=0$, et $H$ est de donc de type horosph{\'e}rique, ce qui est
interdit par hypoth{\`e}se.$\diamond$

\subsection{Surfaces tubulaires et tubes} Introduisons une d{\'e}finition interm{\'e}diaire.
 
Soit  $\wnn$ une suite convergente de
points de $\Snn$. Notons alors $y_{n}\in\S_{n}$ le point pied associ{\'e} {\`a}
$w_{n}=exp(\l_{n}(u_{n})u_{n})$, ({\it cf} \ref{deflent}) et $\S_{n}^{\e}$ la composante connexe de
$\S_{n}\cap B(y_{n},\e)$, o{\`u} $\e$ est fourni par le th{\'e}or{\`e}me \ref{comrap}de
compacit{\'e} sur les probl{\`e}mes de Monge-Amp{\`e}re. Nous dirons
que la suite $\wnn$ est {\it critique}, si apr{\`e}s extraction d'une sous
suite, $\{n(\S_{n}^{\e}),u_{n})\inn$ converge vers un tube, \cad si nous sommes
dans le deuxi{\`e}me cas de l'alternative de \ref{comrap}, nous dirons
de plus dans ce cas, par abus de langage,  que la suite critique {\it
converge} vers $N(\g )$. Remarquons que dans ce cas, d'apr{\`e}s le th{\'e}or{\`e}me
cit{\'e}, {\`a} parir d'un certain rang $y_{n}$ est un point int{\'e}rieur {\`a} $S_{n}$ et $u_{n}$ est le vecteur
normal {\`a} $y_{n}$.

Nous d{\'e}montrerons deux propositions, nous avons tout d'abord

\pro{tub1}{Si $\Snn$ n'est pas de type horosph{\'e}rique et s'il existe une
suite critique, alors la surface
$S^{\infty}$ est tubulaire}

\preu raisonnons par l'absurde et supposons nous donn{\'e}e une suite
critique $\wnn$. Soit $y_{n}$ la suite de points pieds associ{\'e}s dans $\Sinn$,
et posons $u_{n}=n(y_{n})\in n(\S_{n})$ Le th{\'e}or{\`e}me \ref{comrap} assure
alors que $\{n(\S_{n}),u_{n}\inn$ converge en tant que vari{\'e}t{\'e} immerg{\'e}e
point{\'e}e vers un rev{\^e}tement  ${\bar N(\g )}$ d'un tube point{\'e} $N(\g
)$. Comme $\Snn$ elle m{\^e}me converge, la suite de fonctions inverse
$\{\mu_{n}\inn$ va elle aussi converger sur tout compact vers la fonction
fonction  $d_{\g} $ distance {\`a} cette g{\'e}od{\'e}sique qui sera born{\'e}e par
d'apr{\`e}s \ref{premaj}. Par construction, le champ de gradient de
$d_{\g}$ sera dirig{\'e} ver l'int{\'e}rieur de $S^{\infty}$.  Par d{\'e}finition, $S^{\infty}$
est tubulaire. $\diamond$

R{\'e}ciproquement, nous avons

\pro{tub2}{ Si $\Snn$ est {\`a} g{\'e}om{\'e}trie born{\'e}e et n'est pas de type
ho\-ros\-ph{\'e}\-ri\-que, si $S_{\infty}$ est tubulaire pour $\g$ alors toute suite
convergente est critique et converge vers $N(\g )$}

\preu  soit $\g$ la g{\'e}od{\'e}sique associ{\'e}e {\`a}  $S_{\infty}$ et $\mu$ la
fonction distance {\`a} cette g{\'e}od{\'e}sique. Raisonnons par l'absurde et soit
$\wnn$ une suite convergente de points non critique. Notons comme
d'habitude $\ynn$, $y_{n}\in \S_{n}$ la suite de points pieds
associ{\'e}es. En utilisant le th{\'e}or{\`e}me \ref{comrap}, nous pouvons donc
extraire une sous-suite telle que $\{(\S_{n},y_{n})\inn$ converge vers une
$k$-surface, $\S_{\infty}$, {\'e}ventuellement d{\'e}g{\'e}n{\'e}r{\'e}e.

D'apr{\`e}s la proposition \ref{premaj} et le fait que $S_{\infty}$ est
tubulaire nous en d{\'e}duisons que $\mu$ est born{\'e}e sur $\S_{\infty}$. Soit $\L$ sa borne sup{\'e}rieure.  

Construisons maintenant une suite de points $\snn$,
$s_{n}\in\S_{\infty}$ telle que $\mu (s_{n})$ tende vers $\L$. La
vari{\'e}t{\'e} $M$ {\'e}tant {\`a} g{\'e}om{\'e}trie born{\'e}e ainsi que la suite $\Snn$, nous pouvons extraire une sous
suite telle que 
$$
\{\Delta_{n}= (M,s_{n}, \bo S_{n},\bo S_{n})\inn
$$ 
converge vers
$\Delta_{0}=(M_{0},s_{0},\bo S_{0},S_{O})$.

Nous voulons maintenant appliquer notre th{\'e}or{\`e}me de compacit{\'e}
\ref{comrap} {\`a} la suite $\{(\S_{\infty},s_{n})\inn$ solutions des probl{\`e}mes
de Dirichlet $\Delta_{n}$.

Nous avons donc {\it a priori} deux possibilit{\'e}s, 
\begitem
\item[(a)] soit
$\{(\S_{\infty},s_{n})\inn$ converge vers une $k$-surface, {\'e}ventuellement
d{\'e}\-g{\'e}\-n{\'e}\-r{\'e}e, 

\item[(b)] soit $\{(n(\S_{\infty}),n(s_{n}))\inn$ converge vers un tube
$N({\bar\g})$.
\enditem

Eliminons (b) : puisque $\mu$ est born{\'e}e sur
$\S_{\infty}$, n{\'e}\-ces\-sai\-re\-ment ${\bar\g}$ et $\g$ sont {\`a} distance
born{\'e}e. Ces deux g{\'e}od{\'e}siques ont alors confondues. Autrement dit, $\L$
est nulle ce qui entra{\^\i}ne que la restriction de $\mu$ est nulle sur $\S_{\infty}$. Ceci est impossible
pour une $k$-surface.

Ainsi nous ne pouvons {\^e}tre que dans le cadre de (a) et  $\{(\S_{\infty},s_{n})\inn$ converge vers une $k$-surface, {\'e}ventuellement
d{\'e}\-g{\'e}\-n{\'e}\-r{\'e}e, $({\bar\S}_{\infty},s_{\infty})$. La fonction $\mu$ restreinte
{\`a} ${\bar\S}_{\infty}$ atteint son
maximum en $s_{\infty}$ par construction. Ceci signifie que
${\bar\S}_{\infty}$ est in\-t{\'e}\-rieu\-re\-ment tangente {\`a} la surface convexe
$G_{\L}$ {\`a} distance constante $\L$ de $\g$. Ceci est impossible {\`a} cause du principe
du maximum g{\'e}om{\'e}trique puisque $G_{\L}$ est {\`a} courbure extrins{\`e}que
strictement plus grande que $-c$.

Nous venons de montrer que toute suite convergente est critique. Par
ailleurs, si une suite critique converge vers $N({\bar\g})$, nous avons
$\bar\g =\g$  puisque
$\mu$ est born{\'e}e,. $\diamond$

\subsubsection{D{\'e}monstration du th{\'e}or{\`e}me \ref{compalent}}
Soit donc $\Snn$ une suite de surfaces {\`a} courbure sup{\'e}rieure {\`a} $c$,
point{\'e}e en $\xnn$, et convergeant vers $S_{\infty}$ point{\'e}e en
$x_{\infty}$. Soit enfin $\Sinn$ une suite de $k$-surfaces lentilles pour
$\Snn$ et $\ynn$ la suite
des pieds de $\xnn$. On suppose  que la suite $\Snn$ n'est pas de type
horosph{\'e}rique et est {\`a} g{\'e}om{\'e}trie born{\'e}e.

Si $S_{\infty}$ est tubulaire, la proposition \ref{tub2} est exactement
(iii) de \ref{compalent}.

Supposons maintenant que $S_{\infty}$ n'est pas tubulaire.
D'apr{\`e}s la proposition \ref{tub1}, aucune suite
convergente  de points de $\Snn$ n'est critique. Le th{\'e}or{\`e}me \ref{comrap} nous
permet donc d'assurer que la suite $\Sinn$ point{\'e}e en $\ynn$, converge
vers une solution {\'e}ventuellement d{\'e}g{\'e}n{\'e}r{\'e}e  $\S_{\infty}$ du probl{\`e}me de Dirichlet pour
$S_{\infty}$. Il nous reste {\`a} montrer que $\S_{\infty}$ est lentille
pour $S_{\infty}$.

 Notons comme pr{\'e}c{\'e}demment $\{\l_{n}\inn$ la suite de
fonctions associ{\'e}es convergeant vers $\l_{\infty}$, pour conclure il
suffit de  d{\'e}montrer :
\begitem
\ita $\l_{\infty}$ est strictement positive sur $S_{\infty}$;
\itb si $u$ est le vecteur normal int{\'e}rieur {\`a} $\bo S_{\infty}$
dans $S_{\infty}$ et si $n$ est le champ de vecteur normal ext{\'e}rieur {\`a}
$\S_{\infty}$ alors $\langle u,n\rangle >0$ le long de  $\bo
S_{\infty}$.
\enditem
La preuve de (a) est simple : nous savons d{\'e}j{\`a} que $\l_{\infty}$ est
positive ou nulle puisque les $\l_{n}$ le sont. Si maintenant,
$\l_{\infty}$ est nulle en un point int{\'e}rieur {\`a} $S_{\infty}$ nous en
d{\'e}duisons que $\S_{\infty}$ est tangente int{\'e}rieurement {\`a} $S_{\infty}$ en
ce point ce qui est impossible {\`a} cause des conditions de courbure et du
principe du maximum g{\'e}om{\'e}trique \ref{max}.

Pour (b), remarquons que, par passage {\`a} la limite  $\langle u,n\rangle \geq 0$ le long de  $\bo
S_{\infty}$. Si en un point $z$,  $\langle u,n\rangle = 0$, alors les deux
surfaces sont tangentes en ce point $z$. Elles ne peuvent {\^e}tre tangentes
ext{\'e}rieurement : si deux surfaces convexes sont tangentes ext{\'e}rieurement
en un point, alors au moins au voisinage de ce point leur intersection
est r{\'e}duite {\`a} ce point. Rappelons maintenant  que chaque $S_{n}$
est un graphe au dessus d'un ensemble ferm{\'e}  du normal {\'e}tendu {\`a}
$\S_{n}$. Ceci va {\^e}tre {\'e}galement vrai {\`a} la limite ; puisqu'enfin en $z$
les surfaces sont tangentes,  $\S_{\infty}$ est tangente
int{\'e}rieurement {\`a} $S_{\infty}$, ce qui est impossible {\`a} cause
des hypoth{\`e}ses de courbure et du principe du maximum g{\'e}om{\'e}trique.

Enfin, \ref{tub1} et \ref{tub2} montre que $S_{\infty}$ est tubulaire
{\`a} l'infini, si et seulement si il existe une suite $\wnn$ telle que
$\{M,w_{n}\inn$ converge et $\{\S_{\infty},w_{n}\inn$ est critique. En
particulier, si $S_{\infty}$ n'est pas tubulaire
{\`a} l'infini alors la courbure moyenne des $\S_{n}$ est uniform{\'e}ment
born{\'e}e. $\diamond$

\section{Probl{\`e}me de Plateau pour les disques}

La d{\'e}monstration de notre lemme de Morse va se faire par {\'e}tapes et la
premi{\`e}re d'entre elles est la r{\'e}solution du probl{\`e}me de Plateau pour
les disques.

Nous allons en fait d{\'e}montrer un r{\'e}sultat plus pr{\'e}cis qui est le but de
cette section :

\pro{platdisk}{ Soit $S$ un disque compact immerg{\'e} {\`a} courbure
$>c$, dans une vari{\'e}t{\'e} d'Hadamard {\`a} courbure plus petite que $-c$.
Alors, pour
tout $k\in ]0,c[$,  il existe une unique $k$-surface
lentille $\S$ pour $S$.

De plus, si $c>k_{1}\geq k>0$, si ${
\S}_{1}$ est la $k_{1}$-surface lentille pour un disque $S_{1}\subset S$,  alors $\S$ domine $\S_{1}$ au sens de
 \ref{}.}

Nous allons utiliser une m{\'e}thode de d{\'e}formation.

\subsection{D{\'e}formation contractante} Soit $k\in \rbrack 0,c\lbrack$.
 Notons $f_{1}$ l'immersion  du disque $S$ dans $M$. Identifions d'une mani{\`e}re ou d'une
 autre ce disque {\`a} la boule unit{\'e} de ${\mathbb R}^{2}$. Soit alors $f_{t}$,
 ~$t\in\ooi$, la famille d'immersions de $S$ dans $M$ d{\'e}finie par $f_{t}(x)=f_{1}(tx)$.

Le but de  ce paragraphe est de d{\'e}montrer la

\pro{defodisk}{Soit $\S$ une $k$-surface lentille pour $f_{1}$, et
$k(t)$ une fonction croissante {\`a} valeurs dans $]0,c[$ telle que
$k(0)=k$. Il
existe alors une unique famille continue de $k(t)$-surfaces lentilles $\S_{t}$
pour $f_{t}(S)$, d{\'e}finie pour tout $t\in\ooi$ et v{\'e}rifiant $\S_{1}=\S$. Enfin, quand $t$ tend vers
0, $\S_{t}$ tend vers $f_{1}(0)$ au sens de la topologie de Haussdorff
pour les compacts}

\preu  nous utilisons la proposition \ref{defoloc} pour construire au
voisinage de 1 une famille de d{\'e}formations de $\S_{1}$. D'apr{\`e}s la
deuxi{\`e}me partie de cette assertion, nous obtenons de plus que si $t>s$
alors $\S_{s}$ domine $\S_{t}$.
 
Supposons nous donn{\'e}e une d{\'e}formation de  $\S_{1}$ pour $t$ dans un intervalle
$\rbrack a,1\rbrack$.
Pour conclure la d{\'e}monstration de la premi{\`e}re partie dela proposition
\ref{defodisk}, il nous reste {\`a} montrer que pour toute
suite $\{a_{n}\inn$ tendant vers $a$ par valeurs sup{\'e}rieures,
alors la suite $\{\S_{a_{n}}\inn$ converge vers une surface lentille
pour $f(S_{a})$.

Le corollaire de notre th{\'e}or{\`e}me de compacit{\'e} \ref{compalent} nous permet d'affirmer
qu'apr{\`e}s extraction {\'e}ventuelle d'une sous suite,  $\{\S_{a_{n}}\inn$
converge vers une $k$-surface. 

Enfin si $\{a_{n}\inn$ et  $\{{\bar a}_{n}\inn$ sont deux
suites convergeant vers $a$ par valeurs sup{\'e}rieures, telles que les suites
$\{\S_{a_{n}}\inn$ et $\{\S_{{\bar a}_{n}}\inn$ convergent respectivement
vers $\S_{a}$ et $\S_{{\bar a}}$, nous obtenons que $\S_{a}$ domine
$\S_{{\bar a}}$ et r{\'e}ciproquement. Ainsi $\S_{a}=\S_{{\bar a}}$, ce qui
ach{\`e}ve de d{\'e}montrer que si $t$ tend vers $a$ par valeurs inf{\'e}rieures,
$\S_{t}$ converge.

Ce dernier point termine la d{\'e}monstration de la premi{\`e}re partie de notre
proposition.

La deuxi{\`e}me partie de la poposition d{\'e}coule imm{\'e}diatement de  la
proposition \ref{lentboul} qui affirme : si $S$ est une surface
compacte dont le bord $\bo S$ est incluse dans une boule $B$, alors
toute $k$-surface $\S$, pour  $k\leq c$, solution du probl{\`e}me de
Dirichlet pour $S$
est incluse dans la boule $B$.$\diamond$

\subsection{Unicit{\'e}} Nous pouvons maintenant montrer 

\pro{unidisk}{Soit $S$ un disque compact immerg{\'e} {\`a} courbure  $>c$,  $k_{1}$ et $k$ tels que $0<k\leq
k_{1}<c$. Si $\S$ est une $k$-surface lentille pour $S$, et $\S_{1}$ une
$k_{1}$-surface lentille pour $S_{1}\subset S$, alors $\S$ domine
$\S_{1}$.

En particulier, il existe au plus une
$k$-surface
lentille pour $S$. }

\preu consid{\'e}rons comme dans le paragraphe pr{\'e}c{\'e}dent, $f_{1}$  l'immersion
de notre disque $S_{1}$ et $f_{t}$  la famille contractante d'immersions
associ{\'e}e, soit de plus $k(t)$ une fonction strictement d{\'e}croissante {\`a}
valeurs dasn $]0,c[$, telle que $k_{1}=k(1)$.

Soit 
${\S}_{1}$ comme dans l'{\'e}nonc{\'e};  la proposition \ref{defodisk} nous
permet de lui associer une famille  de $k_{1}$-surfaces lentilles ${\S}_{t}$
pour $f_{t}$. Pour montrer que $\S$ domine $\S_{1}$, il suffit 
de montrer que $\S$ domine ${\S}_{t}$ pour tout $t$.

En vertu de notre lemme de domination \ref{domi} et du principe du
maximum g{\'e}om{\'e}trique \ref{max}, il suffit en fait de
montrer qu'il existe $t$, tel que $\S$ domine strictement $\S_{t}$. Mais ceci
d{\'e}coule de la derni{\`e}re partie de la proposition
\ref{defodisk} : clairement, puisque ${ \S}_{t}$ converge vers un
point dans la topologie de Haussdorf quand $t$ tend vers 0, $\S$
domine ${\S}_{t}$, pour $t$ suffisamment petit.$\diamond$

\subsection{Existence} Pour conclure la d{\'e}monstration de \ref{platdisk}, il
nous faut montrer l'existence d'une surface lentille pour notre disque
immerg{\'e} $S$, nous allons pour cela construire une nouvelle famille de
d{\'e}formations de ce disque. 

Ce disque {\'e}tant compact, il est inclus dans une boule $B$. Nous pouvons
alors pertuber la m{\'e}trique de $M$  dans un petit ouvert $U$ n'intersectant
pas cette  boule $B$ de fa\c con {\`a} ce que la nouvelle m{\'e}trique soit toujours {\`a}
courbure plus petite que $-c$, et que de plus  elle soit {\`a} courbure
constante dans une boule $B_{0}\subset U$.

Il suffit de montrer de montrer l'existence d'une surface lentille pour
la m{\'e}\-tri\-que perturb{\'e}e puisque d'apr{\`e}s \ref{lentboul}, cette surface
lentille sera incluse dans $B$ et sera donc lentille pour la m{\'e}trique originelle.

Il est maintenant facile mais technique de montrer le lemme

\lem{}{Soit $S$ un  disque immerg{\'e} {\`a} courbure extrins{\`e}que plus grande que
$c_{1}$, il existe alors une d{\'e}formation continue $S_{t}$,~ $S_{0}=S$ de ce disque
$t\in \oi$,
telle que
\begitem
\iti pour tout $t$, la courbure extrins{\`e}que de $f_{t}$ est plus grande
que $c_{1}$ ;
\itii $S_{1}\subset B_{0}$ est {\`a} courbure constante, et son bord est
inclus dans un plan totalement g{\'e}od{\'e}sique de $B_{0}$.
\enditem}
\preu nous allons donner une esquisse de la construction de cette d{\'e}formation. Nous
utilisons tout d'abord la d{\'e}formation contractante pour d{\'e}former notre
disque en en tout petit disque inclus dans le bord d'une surface
convexe. Ce petit disque va alors {\^e}tre un graphe au dessus d'un disque
topologique sur une sph{\`e}re de petit rayon. Ensuite nous d{\'e}pla\c cons
cette sph{\`e}re de telle sorte quelle soit incluse dans $B_{0}$ et le tour
est jou{\'e}.$\diamond$

Utilisons cette d{\'e}formation pour construire une $k$-surface lentille
pour $S$.

En utilisant les {\'e}quidistantes aux plans totalement g{\'e}od{\'e}siques dans
l'es\-pa\-ce hyperbolique, nous construisons une $k$-surface lentille pour
$S_{1}$.

Notre lemme \ref{defoloc} nous permet de construire une famille
continue de k-surfaces lentilles $\S_{t}$ pour $S_{t}$ sur un intervalle maximal
$\rbrack a,1\rbrack$. Le th{\'e}or{\`e}me de compacit{\'e} \ref{compalent} et
l'unicit{\'e} \ref{unidisk}, entra{\^\i}nent que $\S_{t}$ converge lorsque $t$ tend
vers $a$ ; ainsi $a=0$ et nous avons termin{\'e} notre d{\'e}monstration.

\section{D{\'e}monstration du lemme de Morse pour les surfaces convexes}
 
Dans toute cette section, $S$ d{\'e}signera une surface localement convexe, 
{\'e}\-ven\-tuel\-le\-ment {\`a} bord, {\`a} courbure plus grande que $c$, {\`a}
g{\'e}om{\'e}trie born{\'e}ee et qui n'est ni horosph{\'e}rique {\`a} l'infini, ni
tubulaire, ni compacte sans bord. 

Le but de cette section est  le lemme suivant, que nous appelons lemme
de Morse pour les surfaces convexes.

\lem{morse}{ Soit $M$ une vari{\'e}t{\'e} d'Hadamard {\`a} g{\'e}om{\'e}trie born{\'e}e et {\`a}
courbure strictement plus petite que $-c<0$. Soit  $S$  une surface
localement convexe, {\'e}ventuellement {\`a} bord, {\`a} courbure plus grande que $c$, {\`a}
g{\'e}om{\'e}trie born{\'e}ee et qui n'est ni horosph{\'e}rique {\`a} l'infini, ni
tubulaire, ni compacte sans bord.
 Alors, pour tout $k\in \rbrack 0, c\lbrack$,
il existe une unique $k$-surface lentille, {\'e}ventuellement d{\'e}g{\'e}n{\'e}r{\'e}e,
pour $S$.

De plus, si  $S$ n'est pas tubulaire {\`a} l'infini, cette
$k$-surface est non d{\'e}g{\'e}n{\'e}r{\'e}e.}

Il suffit de d{\'e}montrer ce r{\'e}sultat quand $S$ est un disque topologique : dans le cas
g{\'e}n{\'e}ral, il suffit en effet de passer au rev{\^e}tement universel et
d'utiliser l'unicit{\'e} pour conclure.

Nous supposerons donc dans toute la suite que $S$ est un disque. Notre
lemme de Morse est une cons{\'e}quence imm{\'e}diate des deux propositions
suivantes \ref{exislen} et \ref{unilen}, o{\`u} nous d{\'e}montrons
successivement l'existence et l'unicit{\'e} d'une $k$-surface lentille pour $S$.

\subsection{Existence} Nous allons montrer la proposition :

\pro{exislen}{Il existe une $k$-surface lentille $\S$ pour $S$, telle que pour
toute $k$-surface lentille ${\bar\S}$ pour $S$, ${\bar\S}$ domine
$\S$. De plus si $S$ n'est pas tubulaire {\`a} l'infini, $\S$ est non d{\'e}g{\'e}n{\'e}r{\'e}e.}

\preu consid{\'e}rons $S_{t}$ une famille continue de disques inclus dans $S$ d{\'e}finie
pour $t\in ]0,1]$ et
v{\'e}rifiant
\begitem
\iti $S_{1}=S$ ;
\itii $S_{t}$ est compact {\`a} bord pout $t\not= 1$;
\itiii $S_{t}\subset S_{s}$ si $t\leq s$ ;
\itiv $S_{t}$ converge vers un point $x_{0}$ quand $t$ tend vers 0.
\enditem
La proposition \ref{platdisk} nous permet donc de construire une
famille $\S_{t}$ de $k$-surfaces lentilles pour $S_{t}$. par ailleurs, notre
th{\'e}or{\`e}me de compacit{\'e} \ref{comrap}, nous permet d'extraire une suite
$\{t_{n}\inn$ tendant vers 1, telle que $\{\S_{t_{n}}, y_{n}\inn$ converge vers
une $k$-surface lentille point{\'e}e $\S$ quand $n$ tend vers l'infini. Ici, $\ynn$
d{\'e}signe bien s{\^u}r la suite de points pieds de $x_{0}$ dans $\S_{t_{n}}$.

Si maintenant ${\bar S}$ est une autre surface lentille, il est clair que
${\bar \S}$ domine $\S_{t}$ pour $t$ petit. Ensuite notre lemme \ref{domi}
et le principe du maximum \ref{max} assure que $\bar \S$ domine strictement $\S_{t}$ domine $t$ pour tout $t$. Ainsi
${\bar \S}$
domine-t-elle {\'e}galement $\S$. $\diamond$

\subsection{Unicit{\'e}} L'unicit{\'e} d'une $k$-surface lentille d{\'e}coule imm{\'e}diatement
de la proposition pr{\'e}c{\'e}dente \ref{exislen} et de
\pro{unilen}{Si $\S_{1}$ et $\S_{2}$ sont deux $k$-surfaces lentilles pour $S$ telles que
$\S_{1}$ domine $\S_{2}$ alors $\S_{1}=\S_{2}$}
\preu  soit $N_{\S_{1}}$ le fibr{\'e} normal {\'e}tendu. Soit $F$ l'immersion
naturelle de $M_{1}= N_{\S_{1}}\times \rr^{+}$ dans $M$ qui {\`a} $(n,t)$ associe
$exp(tn)$. Notons $\mu$ la fonction de $M_{1}$ qui a $(n,t)$ associe
$t$, $G_{\l}$ la surface localement convexe de niveau $\l$ pour la
fonction $\mu$, et enfin $\pi$ la projection de $M_{1}$ dans $\S_{1}$.

Par d{\'e}finition, puisque $\S_{1}$ domine $\S_{2}$, nous pouvons trouver une
immersion $i$ de $\S_{2}$ dans $M_{1}$, telle que $F\circ i$ est
l'immersion inititale de $\S_{2}$.  

Confondons $\S_{2}$ et son image dans $M_{1}$. Notre but est de montrer
que $\mu$ restreinte {\`a} $\S_{2}$ est nulle. 

Supposons donc le contraire et soit $\wnn$ une suite de points de
$S_{2}$, telle que $\{\mu (w_{n})\inn$ converge vers la borne sup{\'e}rieure
$\L \not= 0$ de $\mu$ sur $\S_{2}$.

Extrayons de $\wnn$ une sous-suite telle que $\{M,w_{n}\inn$ converge
vers une vari{\'e}t{\'e} $(M_{0},w_{\infty})$,
et utilisons le th{\'e}or{\`e}me de compacit{\'e} \ref{comrap} pour extraire 
extraire de la suite $\wnn$ une sous-suite telle que les deux suites
$\{\S_{1},z_{n})\inn$ et $\{\S_{2},w_{n}\inn$ convergent. Ici
$z_{n}=\pi\circ i (w_{n})$ o{\`u} $\pi$ est la projection naturelle de
$M_{1}$ dans $\S_{1}$.

En ce qui concerne $\{\S_{1},z_{n}\inn$, nous avons comme toujours
deux possibilit{\'e}s : 
\begitem
\item[-] soit   $\{\S_{1},z_{n}\inn$ converge vers une
$k$-surface,

\item[-] soit  $\{n(\S_{1}),n(z_{n}))\inn$ converge vers un
tube. 
\enditem
N{\'e}anmoins, dans ces deux cas, la suite de surfaces convexes
point{\'e}es
$$
\{G_{\L},exp(\L n(z_{n}))\inn
$$ 
converge vers une surface
strictement convexe
lisse {\`a} courbure strictement plus grande que $k$. 

Examinons maintenant  $\{\S_{2},w_{n}\inn$, nous avons aussi deux cas.

(a) Soit  $\{\S_{2},w_{n}\inn$ converge vers une $k$-surface ; le
principe du maximum g{\'e}om{\'e}trique \ref{max} exclut cette possibilit{\'e} puisque cette $k$-surface serait alors tangente
int{\'e}rieurement d'apr{\`e}s notre construction {\`a} $G_{\L}$.

(b) Soit $\{n(\S_{2}),n(w_{n})\inn$ converge vers un tube $N(\g)$, ceci
est {\`a} nouveau exclus puisque alors la g{\'e}od{\'e}sique $\g$ serait tangente
int{\'e}rieurement {\`a} $G_{\L}$ qui est strictement convexe.

Ainsi $\L =0$, et nous avons termin{\'e} notre d{\'e}monstration.$\diamond$
%%%%%%%%%%%%%%%%%%%%%%%%%%%%%%%%%10%%%%%%%%%%%%%%%%%%%%%%%%%%%%%%%%%
\section{Probl{\`e}mes asymptotiques}

Rappelons le probl{\`e}me de Plateau asymptotique : nous nous donnons une
immersion $i$ d'une surface $S$ dans $\bo_{\infty}M$ le bord {\`a} l'infini de $M$,
 et nous cherchons une immersion $f$ de $S$ dans $M$,
telle que 
\begitem
\iti $f(S)$ soit une $k$-surface {\'e}ventuellement d{\'e}g{\'e}n{\'e}r{\'e}e 
\itii si $N$ est l'application de Gauss-Minkowski de $f(S)$ dans
$\bo_{\infty}M$ qui a un point associe le point {\`a} l'infini de la normale
ext{\'e}rieure asymptote, alors 
$$
i=N\circ f.
$$
\enditem
Le couple $(i,S)$ sera la {\it donn{\'e}e} de ce probl{\`e}me de Plateau asymptotique
et $(f,S)$, ou plus simplement $f$, sera une {\it  solution} de ce probl{\`e}me.

Nous allons r{\'e}soudre ce probl{\`e}me tout d'abord dans le cas o{\`u} $S$ est le
disque ouvert et o{\`u} $f$ s'{\'e}tend en un plongement du disque ferm{\'e}. Nous
utiliserons ce r{\'e}sultat pour montrer l'unicit{\'e} de la solution du
probl{\`e}me de Plateau, puis l'existence dans les cas qui nous int{\'e}ressent.

\subsection{Disques plong{\'e}s} Nous voulons montrer le
\lem{diskplon}{ Soit $f$ un plongement du disque ferm{\'e} ${\bar D}$ dans
le bord {\`a} l'infini $\bo_{\infty}M$, alors il existe une unique $k$-surface
plong{\'e}e solution du probl{\`e}me de Plateau {\`a} l'infini d{\'e}fini par $f$
restreinte au disque ouvert $D$}

\preu il suit des constructions de \cite{An}, le fait, d{\'e}sormais
classique, qu'il existe une surface convexe $S$, plong{\'e}e hom{\'e}omorphe au
disque, bordant un convexe $C$  telle que l'application de Gauss-Minkowski soit un
hom{\'e}omorphisme de $S$ sur $f(D)$ ; on consid{\`e}re pour cela  $S_{0}$ le bord de
l'enveloppe convexe de $\bo_{\infty}M\setminus f(D)$. En rempla{\c c}ant au
besoin  cette surface par la surface $S_{\varepsilon_{0}}$ {\`a} distance
constante $\varepsilon_{0}$ o{\`u} $\varepsilon_{0}$ ne d{\'e}pend que de $k$ et pas de $S_{0}$, on peut de plus
supposer qu'elle a une courbure strictement plus grande que $k$ ({\it
cf.} corollaire \ref{equic}).

La $k$-surface $\S_{k}$ lentille pour $S_{\varepsilon}$ produite par le
lemme de Morse \ref{morse} est alors une solution du probl{\`e}me asymptotique. Il
s'agit maintenant de montrer l'unicit{\'e} de cette solution que nous
appelerons pour le moment {\it canonique}. Soit donc
${\bar \S}_{k}$ une autre solution, nous voulons montrer que ${\bar \S}_{k}=\S_{k}$.

Pour cela, nous allons raisonner en deux temps.

(i) Les trois surfaces ${\bar \S}_{k}$, $\S_{k}$ et $S_{0}$ bordent des convexes que
nous noterons   ${\bar O}_{k}$, $O_{k}$ et $U_{0}$ respectivement. Nous
voulons montrer que  ${\bar \S}_{k}$ se trouve {\`a} l'ext{\'e}rieur de
$\S_{k}$, c'est-{\`a}-dire que $O_{k}\subset {\bar O}_{k}$.

Bien s{\^u}r, nous avons
$$
U_{0}\subset {\bar O}_{k},~~~~U_{0}\subset {O}_{k}.
$$

Par construction la surface  $\S_{k}$ est coinc{\'e}e entre
$S_{\varepsilon}$ (c'est une surface lentille) et $S_{0}$ (elle borde un
convexe). En particulier lorsque $k$ tend vers 0,  $\S_{k}$ tend vers
$S_{0}$. Notons $F$ la surface convexe {\`a} distance
$\varepsilon_{0}$ de   ${\bar O}_{k}$ ; comme $U_{0}\subset {\bar
O}_{k}$,  
$F$ se trouve {\`a} l'ext{\'e}rieur de $O_{k}$. 

Soit  alors
$t\mapsto k_{t}$ une fonction continue strictement d{\'e}croissante de
$[0,1]$ dans $[k,k_{0}]$ o{\`u} $k_{0}$ est choisi de telle sorte que la
courbure de $F$ soit strictement plus grande que $k_{0}$. 
Par ailleurs soit  $B_{t}$, $t\in ]0,1]$ une suite exhaustive continue de compacts {\`a}
bord lisse de  $F$, $k_{t}$  et $F_{t}$ la suite de $k(t)$-surfaces
solutions du probl{\`e}me de Plateau pour $B_{t}$. Par le principe de
maximum g{\'e}om{\'e}trique et un argument de d{\'e}formation, $F_{t}$ se trouve
toujours {\`a} l'ext{\'e}rieur de $O_{k}$. Enfin, lorsque $t$ tend vers $1$,
nous avons vu que  $F_{t}$  converge vers la solution du probl{\`e}me de
Plateau pour $F$, c'est-{\`a}-dire ${\bar \S}_{k}$. 

Ainsi  ${\bar \S}_{k}$ se trouve {\`a} l'ext{\'e}rieur de $\S_{k}$. 

(ii) Par ailleurs, soit  $D^{r}$ le disque de rayon $r$ inclus dans $D$, et
$k(r)$ une fonction de $r$ strictement d{\'e}croissante.
notons alors $S^{r}_{\varepsilon}$ {\`a} distance
constante $\varepsilon$ du bord de l'enveloppe convexe de $\bo_{\infty}M\setminus
f(D^{r})$ la surface et  $\S^{r}_{k(r)}$ la $k(r)$-solution canonique du probl{\`e}me de Plateau
asymptotique pour  $D^{r}$. Les  familles $\{S^{r}_{\varepsilon}\}_{r\in\oi}$ et
$\{S^{r}_{0}\}_{r\in\oi}$ d{\'e}pendent  continuement
de $r$, ainsi que les surfaces $\{\S_{k(r)}^{r}\}_{r \in\oi}$. Notre deuxi{\`e}me remarque est la suivante : pour $r$
strictement plus petit que 1 ${\bar \S}_{k}$ est n{\'e}cessairement {\`a}
l'int{\'e}rieur de $\S^{r}_{k(r)}$. En effet, dans le cas contraire, nous
pourrions trouver une surface compacte $U$, $U\subset {\bar \S}_{k}$ dont
le bord serait dans $\S^{r}_{k(r)}$ (les deux surfaces s'intersectant
transversalement par le  principe de maximum g{\'e}om{\'e}trique \ref{max}) et en
particulier borderait $V\subset\S^{r}_{k}$. La surface $V$ serait alors
lentille pour $U$. Ceci est impossible : en effet la $k$-surface $W$
lentille pour $V$ serait {\'e}galement lentille pour $U$, et nous aurions
ainsi deux $k$-surfaces lentilles pour $U$, {\`a} savoir $U$ et $W$.

En faisant tendre $r$ vers 1, on obtient bien que  ${\bar \S}_{k}=\S_{k}$.
$\diamond$

\par
Nous allons par la suite avoir besoin d'un lemme plus explicite,
utilisant que $M$ est {\`a} g{\'e}om{\'e}trie born{\'e}e. Ce lemme est l'analogue de la
$\delta$-hyperbolicit{\'e} pour les g{\'e}od{\'e}siques, ou du caract{\`e}re fin des
triangles en courbure strictement n{\'e}gative. Pour cela, notons
$U_{\a}^{(z,u)}\subset M$ le  c{\^o}ne d'angle $\alpha$ autour du vecteur
$u$ en $z$.

Remarquons de plus que par compacit{\'e} et g{\'e}om{\'e}trie born{\'e}e, il existe une
constante $0<\alpha_{0}<\pi/2$ telle que $\alpha \leq \alpha_{0}$
entra{\^\i}ne 
que, pour tout $z$, $u$, le point $z$ 
appartient {\`a} l'enveloppe convexe de $\partial_{\infty}M\setminus \bo_{\infty}U_{\alpha}^{(z,u)}$.

\lem{deltahyper}{ Soit $0\leq \beta <\alpha < \alpha_{0}$, il existe une
constante $\delta (\alpha, \beta)$
ne d{\'e}pendant que de $\alpha$ et $\beta$ telle que si $z$ est un point
quelconque de $M$,  et $S$ la solution du probl{\`e}me de
plateau asymptotique pour $\bo_{\infty}U_{\alpha}^{(z,u)}$, alors elle
est incluse dans $U_{\a}^{(z,u)}$ et toute g{\'e}od{\'e}sique joignant
$z$ {\`a} un point de $\bo_{\infty}U^{(y,u)}_{\beta}$ intersecte $S$ {\`a} une distance plus
petite que $\delta (\alpha, \beta )$ de $y$}

\preu ce lemme d{\'e}coule imm{\'e}diatement du r{\'e}sultat pr{\'e}c{\'e}dent et d'un
argument de compacit{\'e} $\diamond$

\subsection{Unicit{\'e}} Nous voulons montrer 

\theo{uniasym}{Il existe au plus une solution du probl{\`e}me de Plateau
asymptotique. De plus si $(f,S)$ est une $k$-surface (pas n{\'e}cessairement
compl{\`e}te)  et si $\S\subset S$, alors il existe une
solution du probl{\`e}me de Plateau asymptotique pour $(N,\S )$, o{\`u} $N$ est
l'application de Gauss Minkowski, et cette
solution est un graphe au dessus de $(f,\S)$}

Ce th{\'e}or{\`e}me d{\'e}coule imm{\'e}diatement des trois  propositions suivantes :

\pro{uniasym1}{ Supposons que $S$ est un disque,  et soient $(f,S)$ et $(g,S)$ deux solutions  du
probl{\`e}me de plateau asymptotique, il existe alors une troisi{\`e}me solution du
probl{\`e}me de Plateau $(h,S)$ qui est un graphe {\`a} la fois au dessus de
$f(S)$ et $g(S)$.}

\pro{uniasym2}{Soit $\S$ une $k$-surface {\'e}ventuellement d{\'e}g{\'e}n{\'e}r{\'e}e et soit $\S_{1}$ une
$k$-surface qui est un graphe au dessus de $\S$, alors $\S =\S_{1}$}

L'unicit{\'e} d{\'e}coule imm{\'e}diatement de ces deux premi{\`e}res  propositions ( en
passant au besoin au rev{\^e}tement universel ), et la derni{\`e}re
partie du th{\'e}or{\`e}me du r{\'e}sultat suivant, qui va {\^e}tre l'outil fondamental
dans la d{\'e}monstration des deux premi{\`e}res propositions.

\pro{uniasym3}{ Soit  $(f,S)$ une $k$-surface  {\'e}ventuellement d{\'e}g{\'e}n{\'e}r{\'e}e et soit  $\S$ un ouvert de $S$. Il existe alors une
solution du probl{\`e}me de Plateau asymptotique pour $(N,\S)$, o{\`u} $N$ est
l'application de Gauss Minkowski, telle que cette solution est un graphe
au dessus de $(f,\S)$.}

\subsubsection{D{\'e}monstration de la proposition \ref{uniasym3}}\label{uniasym3preu}

\preu examinons tout d'abord le cas o{\`u} $\S$ est un ouvert {\`a} bord lisse
et relativement compact de $S$.

Posons $\S_{0}=f(\S)$, et consid{\'e}rons $\S^{0}_{r}$ la surface ``{\`a}
distance $r$'' de $\S_{0}$, c'est {\`a} dire, en notant  $n$ le champ de vecteur
normal ext{\'e}rieur {\`a} $f(S)$,
$$
\S^{0}_{r}=\{exp(rn(s)),s\in \S_{0}\}.
$$  

Soit ensuite $\S_{r}$ la solution du probl{\`e}me de Plateau pour
$D_{r}$. La famille $\{\S_{r}\}$ forme alors une famille continue, notre but
est bien s{\^u}r de montrer que  $\{\S_{r}\}$ converge, quand $r$ tend vers
$+\infty$ vers une  solution du
probl{\`e}me de Plateau asymptotique que nous recherchons.

(i) Montrons tout d'abord que pour tout $r$, $\S_{r}$ est un graphe au
dessus de $\S_{0}$. 

Remarquons tout d'abord que ceci est vrai pour  $r$ petit. En effet, si $g$ est la fonction correspondant {\`a} la variation
infinit{\'e}simale de  $\S_{r}$ en $r=0$, cette fonction satisfait
$L(g)=0$ et $g\vert_{\bo \S_{0}}>0$, o{\`u} $L$ est l'op{\'e}rateur elliptique
d{\'e}fini en \ref{defoloc}. Il suit de \ref{sign} et du principe de
maximum que $g>0$.

Pour  montrer que pour tout $r$, $\S_{r}$ est un graphe au
dessus de $\S_{0}$, consid{\'e}rons l'ensemble $O$ des $r$ tels que $\S_{r}$
soit un graphe. L'ensemble $O$ est  {\'e}vi\-dem\-ment un ouvert de ${\mathbb R}^{+}$. Soit alors $[0,r_{0}[$ la
composante connexe de $0$. Raisonnons par l'absurde et supposons que
$r_{0}$ est finie.

Pour $r$ plus petit que $r_{0}$, la surface $\S_{r}$ s{\'e}pare le bout $B$
de $\S_{0}$ en deux composantes connexes l'une born{\'e}e et contenant
$\S_{0}$, et l'autre non born{\'e}e. 

Par continuit{\'e}, la  surface  $\S_{r_{0}}$
{\'e}tant limite de graphes, la normale ext{\'e}rieure {\`a}
$\S_{r_{0}}$ va toujours pointer vers une composante connexe de
$B\setminus \S_{r_{0}}$ ne contenant pas $\S_{0}$. 
Mais par
ailleurs, cette  surface n'est pas un graphe
et il existe donc un point $x$ de  $\S_{0}$ en
tel que la normale issue de $x$ est tangent int{\'e}rieurement {\`a}
$\S_{r_{0}}$ en un point $y$. Par convexit{\'e} locale de $\S_{r_{0}}$, la
normale ext{\'e}rieur en $y$ serait alors dirig{\'e}e vers une composante connexe
de $B\setminus \S_{r_{0}}$ contenant  $x$ d'o{\`u} la contradiction.

(ii) Associons donc {\`a} $\S_{r}$, la fonction $f_{r}$ dont elle est le
graphe. Montrons maintenant que pour tout $y$, $f_{r}(y)$ est born{\'e}e.

Soit $u$ le vecteur normal en $y$, il existe alors $\alpha >0 $ petit
tel que $U_{\a}^{(y,u)} $ soit plong{\'e} dans le bout de $\S$. Si $D$ est la
solution du probl{\`e}me de Plateau asymptotique pour
$\bo_{\infty}U_{\a}^{y,u}$, elle va pouvoir {\^e}tre plong{\'e}e dans le bout de
$\S$, nous savons alors que la g{\'e}od{\'e}sique issue de $y$ dans la direction
de $u$ intersecte $D$ au temps $K$. Par le principe du maximum
g{\'e}om{\'e}trique la surface $D$ va alors servir de
barri{\`e}re {\`a} la famille $\{\S_{r}\}$, et nous en d{\'e}duisons que $\forall
r$, $f_{r}(y)\leq K$.

(iv) Nous pouvons maintenant montrer maintenant montrer que la suite
$\{\S_{r}\}$ converge apr{\`e}s extraction {\'e}ventuelle d'une sous-suite vers
une surface $\S_{\infty}$ qui est le graphe d'une fonction $\l$ d{\'e}finie
sur $\S_{0}$ et propre ({\it i.e.} tendant vers l'infini lorsque l'on
tend vers le bord de  $\S_{0}$.

Consid{\'e}rons la suite de fonctions $\{f_{r}\}$ d{\'e}finies au paragraphe
pr{\'e}c{\'e}dent. Pour chaque $y$ la suite
$\{f_{r(y)}\}$ est croissante et born{\'e}e. Nous en d{\'e}duisons
im\-m{\'e}\-dia\-tement que cette suite de fonctions converge vers une fonction
$\l$ propre et de graphe convexe. Pour cela, il faut bien s{\^u}r utiliser
la classique compacit{\'e} des convexes en courbure n{\'e}gative. Il nous reste
{\`a} montrer que le graphe de la limite est une $k$-surface. Ceci d{\'e}coule
imm{\'e}diatement du th{\'e}or{\`e}me \ref{comrap} appliqu{\'e} au probl{\`e}me de
Dirichlet {\`a} bord vide d{\'e}fini sur le bout de $\S_{0}$ et de la remarque
suivante : la suite de surfaces $\{\S_{r}\}$ est toujours transverse au
champ de vecteur normal {\`a}  $\S_{0}$ (ce sont des graphes), et en
particulier on ne peut pas voir appara{\^\i}tre de tubes {\`a} la limite.

(v) Montrons enfin que la surface limite $\S_{\infty}$ est une solution
du probl{\`e}me de Plateau asymptotique pour $(N,\S_{0})$. Notons $B$ le
bout de $S$. La surface $\S_{\infty}$ va s{\'e}parer $B$ en deux composantes
connexes. Soit $C$ la composante connexe de l'ext{\'e}rieur de
$\S_{\infty}$. Par construction $\bo_{\infty}C=N(\S_{0})$, en
particulier, nous en d{\'e}duisons que l'application de Gauss-Minkowski de
$\S_{\infty}$ est injective et {\`a} valeurs dans $\S_{0}$ vu comme
sous-ensemble du bord {\`a} l'infini de $B$ gr{\^a}ce {\`a} l'application de
Gauss-Minkowski de $S$. Il nous reste {\`a} montrer que cette application
est surjective. Pour cela, il suffit de remarquer  que si $v\in N(\S_{0})$ et
si $h$ est la fonction horosph{\'e}rique associ{\'e}e {\`a} $v$, alors $h$ est
propre sur  $\S_{\infty}$, elle admet donc un mimimum c'est {\`a} dire un
point dont la normale pointe vers $v$.

(vi) Il nous reste {\`a} examiner le cas o{\`u}  $\S_{0}$ n'est pas {\`a} bord lisse
et relativement compact. On construit alors une suite $\{\S_{r}^{0}\}$
d'ouverts {\`a} bord lisse et relativement compacts, emboit{\'e}s et dont la
r{\'e}union est $\S_{0}$. D'apr{\`e}s ce que nous venons de voir, nous pouvons
alors construire une suite $\{\S_{r}\}$ de $k$-surfaces dans $B$, le
bout de $\S_{0}$, dont les bouts
sont embo{\^\i}t{\'e}s, et fournissant des solutions des probl{\`e}mes de Plateau
asymptotiques associ{\'e}s {\`a} $\{\S_{r}^{0}\}$. Une adaptation sans douleur
de la d{\'e}monstration pr{\'e}c{\'e}dente montre alors que la suite $\{\S_{r}\}$ va
converger vers une $k$-surface solution de notre probl{\`e}me. $\diamond$

\subsubsection{D{\'e}monstration de la proposition \ref{uniasym2}}

\preu Soit $\S$ une $k$-surface {\'e}ventuellement d{\'e}g{\'e}n{\'e}r{\'e}e et soit $\S_{1}$ une
$k$-surface qui est un graphe au dessus de $\S$. Notons $f$, la fonction
associ{\'e}e. Pour montrer que $\S=\S_{1}$, il nous suffit de montrer que
$f$ est born{\'e}e par une constante $a$. En effet, $\S$ et $\S_{1}$ seront
alors toutes les deux lentilles pour la surface $S_{a}$ {\`a} distance $a$ de
$\S$, et nous aurons l'unicit{\'e} gr{\^a}ce {\`a} la proposition \ref{unilen}.

Montrons donc que $f$ est born{\'e}e. Consid{\'e}rons comme d'habitude le bout
$B$ de $\S$, et commen{\c c}ons par remarquer que dans le cas o{\`u} $\S$ est
compl{\`e}te, alors pour tout $y\in\S$ de vecteur normal ext{\'e}rieur  $n$ et
tout vecteur $u$ tel que $\langle u,n\rangle >0$, la g{\'e}od{\'e}sique issu de
$u$ reste trac{\'e}e dans $B$. Cela reste vrai dans le cas d'une
$k$-surface d{\'e}g{\'e}n{\'e}r{\'e}e : en effet une telle g{\'e}od{\'e}sique heurte transversalement en un temps
fini une surface {\'e}quidistante, et une telle surface est toujours compl{\`e}te.

En particulier,  le c{\^o}ne $U_{\a_{0}}^{y,n}$ est inclus
dans $B$,  o{\`u} $\a_{0}$ est
la constante de \ref{deltahyper} et  $n$ d{\'e}signe 
vecteur normal {\`a} $y$. Soit $D_{y}$ l'unique solution du probl{\`e}me de
Plateau asymptotique pour $\bo_{\infty}U_{\a_{0}}^{(y,n)}$ ({\it cf} \ref{diskplon}), ce disque $D_{y}$
est inclus dans $U_{\a_{0}}^{y,n}$, d'apr{\`e}s \ref{deltahyper} et en
particulier dans $B$ ; par \ref{uniasym3}, c'est donc  le graphe d'une fonction $f_{y}$ au
dessus d'un ouvert $U_{y}$ relativement compact de $\S$. 

Montrons que $f_{y}\geq f$. Dans le cas contraire, soit
$V_{y}=\{z/f_{y}\geq f$. D'apr{\`e}s la propret{\'e}
de $f_{y}$ sur $U_{y}$, $V_{y}$ est compact. De plus les deux surfaces
s'intersectant transversalement ({\`a} cause du principe du maximum
g{\'e}om{\'e}trique) $V_{y}$ est {\`a} bord lisse. Il est alors facile de voir que
le graphe de $f_{y}$ restreint {\`a} $V_{y}$ est lentille pour celui de $f$,
ce qui est impossible d'apr{\`e}s \ref{unilen}.

Ainsi $f(y)\leq f_{y}(y)$, or $f_{y}(y)\leq\delta(\a, 0)$, o{\`u} cette
constante {\`a} {\'e}t{\'e} introduite en \ref{deltahyper}.$\diamond$

\subsubsection{D{\'e}monstration de la proposition \ref{uniasym1}}

\preu soient $(f,S)$ et $(g,S)$ deux solutions du probl{\`e}me de Plateau
asymptotique pour $(i,S)$. Pour tout $y\in S$, soit $D_{y}$ un disque ouvert
tel que 
$$\
i(D_{y})\subset \bo_{\infty}U_{\a_{0}}^{(f(y),n_{f}(y))}\cap \bo_{\infty}U_{\a_{0}}^{(g(y),n_{g}(y))}
$$
o{\`u} $n_{f}$ et $n_{g}$ d{\'e}signent les champ de vecteur normaux ext{\'e}rieurs
associ{\'e}s {\`a} $f$ et $g$ respectivement.

Soit alors $S_{y}$ l'unique solution du probl{\`e}me de Plateau asymptotique
pour $(i,D_{y})$ produite par le lemme \ref{diskplon}. D'apr{\`e}s la proposition \ref{uniasym3}, $S_{y}$ est le
graphe de fonctions $f_{y}$ et $g_{y}$ au dessus de $(f,D_{y})$ et
$(g,D_{y})$ respectivement.

Consid{\'e}rons un recouvrement localement fini de $S$ par des $D_{y}$ o{\`u}
$y$ d{\'e}crit un ensemble $Y$. Et soit $f_{0}=inf (f_{y}),~~y\in Y$ et
$g_{0}=inf (g_{y}],~~ y\in Y$. Par construction le graphe de $f_{0}$ et
celui de $g_{0}$ sont les m{\^e}mes. Nous avons ainsi construit une surface
interm{\'e}diaire $S_{2}$ qui est un graphe {\`a} la fois au dessus de $(f,S)$ et
au dessus de $(g,S)$. Cette surface convexe n'a pas de raison d'{\^e}tre lisse. 

Rempla{\c c}ons cette surface $S_{2}$ par la surface $S_{3}$ {\`a} distance $r$
de de $S_{2}$, pour $r$ bien choisi, ind{\'e}pendamment de $S_{2}$, $S_{3}$
est {\`a} courbure strictement sup{\'e}rieure {\`a} $k$. A nouveau, $S_{3}$ est une
graphe au dessus de $(f,S)$ et $(g,S)$. 

Identifions comme d'habitude $S_{3}$ au disque  $D_{1}$ de rayon 1 de ${\mathbb
C}$. Soit $D_{r}$ les disques de rayon $r$, et $\S_{r}$ les solutions du
probl{\`e}mes de Plateau pour $D_{r}$. Pour $r$ petit, $\S_{r}$ est un
graphe {\`a} la fois au dessus de $(f,S)$ et $(g,S)$. Ceci est {\'e}galement
vrai par continuit{\'e} pour tout $r$, les k-surfaces  $(f,S)$ et $(g,S)$
faisant office de barri{\`e}re.

Quand $r$ tend vers l'infini la suite de surfaces $\{\S_{r}\}$ tend vers
une k-surface $\S_{0}$. Pour voir cela il suffit d'invoquer le m{\^e}me
raisonnement qu'en \ref{uniasym3preu} . .1.(iv).$\diamond$

%%%%%%%%%%%%%%%%%%%%%%%%%%%%%%%%11%%%%%%%%%%%%%%%%%%%%%%%%%%%%%%%%%%%
\subsection{Existence} Nous voulons d{\'e}montrer nos th{\'e}or{\`e}mes d'existence de
solutions du prob\-l{\`e}\-me de Plateau asymptotique
\theo{equib}{Soit $S$ une surface compacte de genre plus grand que $2$,
$\rho$ une repr{\'e}sentation de $\G$ dans le groupe des isom{\'e}tries de $M$,
$f$ un hom{\'e}omorphisme local de ${\bar S}$, le rev{\^e}tement universel de
$S$ dans $\bo_{\infty}M$ {\'e}quivariante sous $\r$, alors il existe une unique
solution non d{\'e}g{\'e}n{\'e}r{\'e}e du probl{\`e}me de Plateau asymptotique pour $(f,{\bar S})$ et cette
solution est {\'e}quivariante sous $\r$.}

\theo{pasequi}{Soit $f$ un hom{\'e}omorphisme local de $U$ dans
$\bo_{\infty}M$ et soit $S$ un ouvert relativement compact de $U$, alors
il existe une unique solution au probl{\`e}me de Plateau asymptotique pour $(f,S)$}

Nous montrerons {\'e}galement ce qui peut-{\^e}tre con{\c c}u comme une
g{\'e}\-n{\'e}\-ra\-li\-sa\-tion faible du petit th{\'e}or{\`e}me de Picard.

\theo{picard}{ Soit $f$ un hom{\'e}omorphisme local de $S$ dans
$\bo_{\infty}M$ {\'e}vitant trois po\-ints. Alors, il existe une solution non
d{\'e}g{\'e}n{\'e}r{\'e}e  du
probl{\`e}me de Plateau asymptotique pour $(f,S)$.}

Pour d{\'e}montrer nos deux th{\'e}or{\`e}mes d'existence, nous allons {\`a} chaque fois
construire une surface convexe immerg{\'e}e dans $M$ et dont l'application
de Gauss-Minkowski d{\'e}finit le un probl{\`e}me de Plateau
asymptotique. Nous construirons cette surface  en recollant des surfaces
obtenues comme solution des probl{\`e}mes de Plateau asymptotiques pour des
disques plong{\'e}s.

Les deux premiers th{\'e}or{\`e}mes vont chacun utiliser des lemmes ayant leur
int{\'e}r{\^e}t propre. Pour le premier, il nous faudra le 

\lem{lentasym}{ Soit $S$ une surface compl{\`e}te, qui n'est ni tubulaire,
ni ho\-ros\-ph{\'e}\-ri\-que {\`a} l'infini, et soit $N$ son application de Gauss
Minkowski {\`a} valeurs dans $\bo_{\infty}M$, alors le probl{\`e}me de
Plateau asymptotique $(N,S)$ admet une solution. Si de plus, $S$ n'est
pas tubulaire {\`a} l'infini, cette $k$-surface n'est pas d{\'e}g{\'e}n{\'e}r{\'e}e}

Pour le deuxi{\`e}me, nous utiliserons le
\lem{convasym}{ Soit $V$ une surface sans bord, pas n{\'e}cessairement
compl{\`e}te, immerg{\'e}e de fa{\c c}on localement convexe dans $M$. Soit $N$ son
application de Gauss-Minkowski {\`a} valeurs dans $\bo_{\infty}$, soit
$S\subset U$ un ouvert relativement compact. alors le probl{\`e}me de
Plateau asymptotique $(N,S)$ admet une solution.}

\subsubsection{D{\'e}monstration du lemme \ref{lentasym}}
\preu soit $S_{1}$ la surface {\`a} distance $r$ de $S$, cette surface n'est
{\'e}videmment
ni horosph{\'e}rique {\`a} l'infini, ni tubulaire. Pour $r$ bien choisi, cette
surface est {\`a} courbure strictement plus grande que $k$. De plus, toute
solution de $(N,S_{1})$ est une solution de $(N,S)$. Il est maintenant
facile de v{\'e}rifier que la $k$-surface lentille pour $S_{1}$, produite
par le lemme de Morse est une solution du probl{\`e}me de Plateau
asymptotique.$\diamond$
\bigskip

\subsubsection{D{\'e}monstration du lemme  \ref{convasym}}
\preu nous allons chercher {\`a} produire la solution du probl{\`e}me de Plateau
asymptotique comme limite de probl{\`e}mes de Plateau, comme dans la preuve
de la proposition \ref{uniasym3}.

Examinons tout d'abord le cas o{\`u} $S$ est un ouvert relativement compact {\`a} bord lisse. Consid{\'e}rons $S_{r}$ les surfaces {\`a} distance $r$ de
$S$, et $\S_{r}$ les surfaces lentilles solutions du probl{\`e}me de Plateau
pour $S_{r}$. Pour adapter la d{\'e}monstration de
\ref{uniasym3} il nous suffit de faire la remarque suivante. Nous
savons que, par d{\'e}finition des surfaces lentilles $S=S_{0}$ est un graphe au
dessus d'un ouvert $O$ du fibr{\'e} normal {\'e}tendu $N_{\S_{0}}$. Si $\S_{r}$ n'a pas de raisons d'{\^e}tre un graphe au dessus de
$\S_{0}$, les arguments de \ref{uniasym3preu}.1.(i), montre
qu'elle sera un graphe au dessus de ce m{\^e}me ouvert $O$. 

Le reste de la preuve, c'est-{\`a}-dire les points (ii), (iii), (iv), et (v)
s'adapte sans probl{\`e}me et nous permettent de construire une solution du
probl{\`e}me de Plateau pour $(N,S)$.

Enfin dans le cas o{\`u} $S$ n'est pas {\`a} bord lisse, il nous suffit de
construire ${\bar S}$ {\`a} bord lisse relativement compact et contenant
$S$. Ce que nous venons de dire nous permet de construire une solution
du probl{\`e}me de Plateau asymptotique pour $(N,{\bar S})$ et l'application
directe de la proposition \ref{uniasym3} permet de conclure. $\diamond$

\subsubsection{D{\'e}monstration du th{\'e}or{\`e}me \ref{equib}}
\preu par compacit{\'e}, nous  pouvons trouver des disques
ferm{\'e}s $D_{i}$,$~i\in J$, dont les int{\'e}rieurs recouvrent ${\bar S}$ le rev{\^e}tement
universel de $S$,  et v{\'e}rifiant, si nous notons 
$$
I_{i}=\{j\in I, D_{i}\cap D_{j}\not=\emptyset\}
$$ 
les hypoth{\`e}ses suivantes
\begitem
\ita Pour tout $i$, le cardinal de $I_{i}$ est fini,
\itb $D_{i}\not\subset \bigcup_{j\in I_{i}\setminus \{i\}}D_{i_{j}}$,
\itc Pour tout $i$, il existe un disque $\D_{i}$ contenant
$\bigcup_{j\in I_{i}}D_{i_{j}}$, tel que $f$ restreinte {\`a} $\D_{i}$ soit
un plongement,
\itd la famille $\{D_{i}\}$ est $\G$-{\'e}quivariante.
\enditem
D'apr{\`e}s le lemme \ref{diskplon}, nous pouvons construire des solutions
$S_{j}$ de chacun des prob\-l{\`e}\-mes de Plateau asymptotiques d{\'e}fini par
$(f,D^{j})$. Chacune de ces surfaces $S_{j}$ d{\'e}coupe $M$ en deux parties
dont l'une est convexe. Nous  appellerons cette composante convexe
$C_{j}$.

Remarquons maintenant que les hypoth{\`e}ses (a), (b) et (c) entra{\^\i}nent l'as\-ser\-tion suivante :
\begitem
\item[(e)] si $j,~k\in I_{i}$, alors $S_{j}\cap
S_{k}\not=\emptyset$ entra{\^\i}ne $D_{j}\cap
D_{k}\not=\emptyset$. 
\enditem
En effet, soit $\S_{i}$ la solution du probl{\`e}me de
Plateau asymptotique pour $\D_{i}$ d{\'e}finie par (d). D'apr{\`e}s la proposition
\ref{uniasym3}, $S_{j}$ est alors un graphe au dessus de
$N^{-1}(D_{j})\subset \S_{i}$ d'o{\`u} (e).
\par

Soit maintenant  
$$
U_{i}=\bigcap_{j\in I_{i}}C_{j}.
$$

Cet ensemble $U_{i}$ est convexe. Son bord est la r{\'e}union de
facettes. Parmi celles ci, nous distinguerons $F_{i}=S_{i}\cap\bo
U_{i}$. La facette  $F_{i}$ est une surface, non vide d'apr{\`e}s (b) dont le bord est une
r{\'e}union d'arcs $C^{\infty}$ pas n{\'e}cessairement connexes,
$\g_{(i,j)}$, $~j\in I_{i}$, portions de $S_{j}\cap
S_{i}$.

Plus pr{\'e}cis{\'e}ment encore 
$$
\g_{(i,j)}=S_{j}\cap
S_{i}\cap U_{i}, 
$$

nous avons m{\^e}me, d'apr{\`e}s (e)

$$
\g_{(i,j)}=S_{j}\cap
S_{i}\bigcap_{k\in I_{i}\cap I_{j}} C_{k}. 
$$

En particulier $\g_{(i,j)}=\g_{(j,i)}$.

Nous pouvons maintenant recoller $F_{j}$ avec $F_{i}$ le long de
$\g_{(i,j)}$, et obtenir ainsi de proche en proche une surface $F$
localement convexe, que nous nous empressons de lisser, de mani{\`e}re $C^{1}$, en prenant la
surface  $\S$ {\`a} distance $r$, et que nous appelons
$\S$. Il est maintenant clair que le probl{\`e}me de Plateau asymptotique
d{\'e}fini par $(N,\S)$, o{\`u} $N$ est l'application de Gauss-Minkowski, est
{\'e}quivalent {\`a} notre probl{\`e}me de d{\'e}part.

Nous voulons appliquer notre lemme \ref{lentasym}. Nous allons montrer
que $\S$ est {\`a} g{\'e}om{\'e}trie born{\'e}e et n'est ni tubulaire {\`a} l'infini, ni horosph{\'e}rique {\`a} l'infini.

Comme l'action de $\G$ sur $\S$ est cocompacte, il suffit bien
{\'e}videmment de montrer que $\S$ n'est ni tubulaire, ni horosph{\'e}rique.

Si $\S$ {\'e}tait horosph{\'e}rique, alors par d{\'e}finition $\G$ agirait de
mani{\`e}re cocompacte sur une horosph{\`e}re. Mais ceci est impossible : avec
nos hypoth{\`e}ses de g{\'e}om{\'e}trie born{\'e}e, la croissance des hororosph{\`e}res est
polynomiales \ref{poly} (ii).

De m{\^e}me, si $\S$ {\'e}tait tubulaire, $\G$ agirait de mani{\`e}re cocompacte sur
le fibr{\'e} unitaire normal {\`a} une g{\'e}od{\'e}sique. Ceci est impossible : ce
fibr{\'e} est conforme au plan priv{\'e} d'un point.

Le lemme \ref{lentasym} permet alors de conclure.$\diamond$

\subsubsection{D{\'e}monstration du th{\'e}or{\`e}me \ref{pasequi}}
\preu nous allons proc{\'e}der comme dans le paragraphe pr{\'e}c{\'e}dent. En effet,
nos hypoth{\`e}ses nous permettent de contruire une famille finie de disques
ouverts $D_{i}$ v{\'e}rifiant $~i\in J$, dont les int{\'e}rieurs recouvrent
l'adh{\'e}rence de  $S$,  et v{\'e}rifiant, si nous notons 
$$
I_{i}=\{j\in I, D_{i}\cap D_{j}\not=\emptyset\}
$$ 
les hypoth{\`e}ses suivantes
\begitem
\item[(b)] $D_{i}\not\subset \bigcup_{j\in I_{i}\setminus \{i\}}D_{i_{j}}$,
\item[(c)] Pour tout $i$, il existe un disque $\D_{i}$ contenant
$\bigcup_{j\in I_{i}}D_{i_{j}}$, tel que $f$ restreinte {\`a} $\D_{i}$ soit
un plongement,
\enditem
D'apr{\`e}s le lemme \ref{diskplon}, nous pouvons construire des solutions
$S_{j}$ de chacun des prob\-l{\`e}\-mes de Plateau asymptotiques d{\'e}fini par
$(f,D^{j})$. Chacune de ces surfaces $S_{j}$ d{\'e}coupe $M$ en deux parties
dont l'une est convexe. Nous  appelerons cette composante convexe
$C_{j}$.

Remarquons maintenant que les hypoth{\`e}ses  (b) et (c) entra{\^\i}nent, comme
dans le paragraphe pr{\'e}c{\'e}dent, l'assertion suivante :
\begitem
\item[e] si $j,~k\in I_{i}$, alors $S_{j}\cap
S_{k}\not=\emptyset$ entra{\^\i}ne $D_{j}\cap
D_{k}\not=\emptyset$. 
\enditem

En proc{\'e}dant exactement comme dans le
paragraphe pr{\'e}c{\'e}dent, nous ob\-te\-nons une surface convexe $\S$ telle que 
le probl{\`e}me de Plateau asymptotique
d{\'e}fini par $(N,\S)$, o{\`u} $N$ est l'application de Gauss-Minkowski, est
{\'e}quivalent au probl{\`e}me de Plateau asymptotique $(N,\bigcup_{i\in
J}D_{i})$. L'ouvert $S$ {\'e}tant par hypoth{\`e}se relativement compact dans $\bigcup_{i\in
J}D_{i})$, le lemme \ref{convasym} permet de conclure.$\diamond$

\subsubsection{D{\'e}monstration du th{\'e}or{\`e}me \ref{picard}}
\preu soit donc $i$ un hom{\'e}omorphisme local  de $S$ dans $\bo_{\infty}M$ {\'e}vitant
trois points $x_{1}$, $x_{2}$ et $x_{3}$. 

Pour d{\'e}montrer le resultat, il nous suffit que le probl{\`e}me de Plateau asymptotique $(i,S_{0})$
a une solution, o{\`u} $S_{0}=\bo_{\infty}\setminus\{x_{1},x_{2},x_{3}\}$ et
$i$ est l'injection. En effet, si cela est vrai,  si $(f,S)$ est un hom{\'e}omorphisme local de
$S$ {\'e}vitant  $x_{1}$, $x_{2}$ et $x_{3}$, nous pouvons relever $f$ en
une application ${\bar f}$ de ${\bar S}$ dans ${\bar S_{0}}$, o{\`u} ${\bar
S}$ et ${\bar S_{0}}$ d{\'e}signent les rev{\^e}tements universels de $S$ et
$S_{0}$ respectivement. La proposition \ref{uniasym3} nous permet de
construire une solution du probl{\`e}me de Plateau asymptotique pour
$(\pi\circ {\bar f}, {\bar S})$ o{\`u} $\pi$ d{\'e}signe la projection canonique
de ${\bar S_{0}}$ dans $S_{0}$. L'unicit{\'e} \ref{uniasym} nous permet
enfin de montrer que la solution ainsi obtenue est bien une solution de
$(f,S)$.
 
Concentrons nous donc sur le cas de $(i,S_{0})$. Distinguons deux cas.

(i) Premi{\`e}rement supposons qu'il existe un triangle id{\'e}al $T$ totalement
 g{\'e}o\-d{\'e}\-si\-que dont les sommets sont nos trois points. La surface
 $S_{\epsilon}$, bord de la boule de centre $T$ et de rayon $\epsilon$ est
 alors {\`a} courbure plus grande que $k$ pour un certain $\epsilon$.
 Elle n'est de plus ni horosph{\'e}rique {\`a} l'infini, ni tubulaire. Le
 lemme \ref{lentasym} nous affirme alors que la $k$-surface lentille
 pour $S_{\epsilon}$ produite par le lemme de Morse fournit une solution
 de notre probl{\`e}me de Plateau. Cette solution borde un convexe, elle va
 donc {\^e}tre compl{\`e}te, autrement dit la $k$-surface ainsi obtenue n'est
 pas d{\'e}g{\'e}n{\'e}r{\'e}e.

(ii) Pla{\c c}ons nous dans le cas g{\'e}n{\'e}ral maintenant. Soit $\{S_{n}\inn$
 une exhaustion de $S_{0}$ par des ouverts relativement compacts. Le
 th{\'e}or{\`e}me \ref{pasequi} nous permet de construire pour tout $n$ une
 solution $(f_{n},S_{n})$ du probl{\`e}me de Plateau asymptotique pour
 $(i,S_{n})$.  Notre but est bien s{\^u}r de montrer que la suite
 $\{f_{n}\inn$
converge apr{\`e}s extraction d'une sous suite.

Soit $x$ un point appartenant {\`a} l'enveloppe convexe de
$\{x_{1},x_{2},x_{3}\}$. Les surfaces $f_{n}(S_{n})$ bordent des
convexes $C_{n}$ qui contiennent tous $x$ et qui de plus sont tels que
$$
n\geq p\Rightarrow C_{p}\subset C_{n}~~~~~(*).
$$

Consid{\'e}rons 
$$
C=\bigcap_{n\in {\mathbb N}}C_{n}.
$$

 Dans le cas dans lequel nous nous
sommes plac{\'e}, l'int{\'e}rieur de $C$ est non vide. Soit donc $x$ tel que la
sph{\`e}re $S_{a}$, de centre $x$ et de rayon $a$, soit incluse dans $C$.
Identifions cette sph{\`e}re  canoniquement {\`a} $\bo_{\infty}M$. 

Les surfaces $f_{n}(S_{n})$
sont maintenant des graphes au-dessus de $S_n$, vus commes
sous-ensembles de $S_{a}$ de  fonctions $\l_{n}$. Par $(*)$,
pour tout $y$, la suite $\{\l_{n}(y)\inn$ est d{\'e}croissante.

Nous pouvons en conclure, en utilisant le corollaire \ref{graphe},
que la suite de fonctions $\{\l_{n}\inn$ converge $C_{\infty}$ vers une
fonction $\l$ dont le graphe va nous fournir la solution de notre
probl{\`e}me de Plateau $(i,S_{0})$. $\diamond$

\subsection{Non existence} Le th{\'e}or{\`e}me suivant montre qu{\'\i}l n'y a pas toujours
des solutions au prob\-l{\`e}\-me de Plateau asymptotique 

\theo{inex}{Si $U$ est $\bo_{\infty}M$ priv{\'e}e de 0, 1 ou 2 points. Le
probl{\`e}me de Plateau asymptotique d{\'e}fini par $(i,U)$, o{\`u} $i$ est
l'injection canonique, n'a pas de solutions.}

\preu raisonnons par l'absurde. Dans les trois cas, la solution serait
une surface $\S$ globalement convexe.

Dans le premier cas, une telle surface serait une sph{\`e}re compacte ce qui est
clairement impossible puisque par l'{\'e}quation de Gauss, la m{\'e}trique
induite sur une $k$-surface est {\`a} courbure n{\'e}gative.

Dans le deuxi{\`e}me cas, $\S$ est une pseudo-horosph{\`e}re, d'apr{\`e}s \ref{psudo}
elle est alors horosph{\'e}rique {\`a} l'infini. Il existe donc une suite de
points $\{x_{n}\inn$ de $\S$ telle que $\{(\S,x_{n})\inn$ converge vers une
surface horosph{\'e}rique $\S_{0}$. Mais ceci est impossible, car
 par compacit{\'e} $\S_{0}$ est {\'e}galement une $k$-surface et une
$k$-surface
ne peut-{\^e}tre horosph{\'e}rique {\`a} cause du principe du maximum g{\'e}om{\'e}trique.

Dans le dernier cas, soit $\gamma$ la
g{\'e}od{\'e}sique joignant les deux points du bord {\`a} l'infini et soit $\mu$ la
fonction distance {\`a} cette g{\'e}od{\'e}sique. Soit ensuite $x_{n}$ une suite de
points de $\S$ tendant vers le maximum $\mu_{0}$ de $\mu$. Nous pouvons
extraire une sous suite telle que $(\S,x_{n})$ converge vers $\S_{0}$.

Mais $\mu_{0}$ ne peut-{\^e}tre born{\'e}, puisqu'alors $\S_{0}$ serait tangente
int{\'e}rieurement {\`a} une {\'e}quidistante {\`a} une g{\'e}od{\'e}sique, ce qui est
impossible par le principe du maximum g{\'e}om{\'e}trique.

De m{\^e}me, $\mu_{0}$ ne peut-{\^e}tre infini, car dans ce cas la surface
$\S_{0}$ serait une pseudo-horosph{\`e}re, et nous venons de voir que c'est impossible.$\diamond$

%%%%%%%%%%%%%%%%%%%%%%%%%%18%%%%%%%%%%%%%%%%%%%%%%%
\section{Espace des $k$-surfaces}\label{laminations}
Abusivement, nous consid{\'e}rerons une $k$-surface comme une surface immerg{\'e}e
dans le fibr{\'e} unitaire de la vari{\'e}t{\'e} M, par son relev{\'e} de Gauss, \cad l'application qui 
a un point de la surface associe sa normale ext{\'e}rieure.  

Soyons  un  
peu plus pr{\'e}cis dans la d{\'e}finition d'une $k$-surface pour {\'e}viter le probl{\`e}me des rev{\^e}tements 
multiples.

Les $k$-surfaces {\'e}tant les solutions d'un probl{\`e}me elliptique, toute telle surface est 
totalement 
d{\'e}termin{\'e}e par son jet d'ordre infini en un point. On en d{\'e}duit qu'{\'e}tant donn{\'e}e une $k$-surface $(f,S)$ 
sans bord, 
$f$ d{\'e}signe une immersion de $S$ dans le fibr{\'e} unitaire,
il existe un {\it repr{\'e}sentant minimal} de la $k$-surface $(g,\S)$, 
tel que pour tout $k$-surface
$({\bar f},{\bar S})$ ayant la m{\^e}me image que $f$ ( \cad $f(S)={\bar f}({\bar S})$), 
alors il
existe un rev{\^e}tement $p$ de ${\bar S}$ sur $\S$, tel que ${\bar f}=g\circ p$.  

le repr{\'e}sentant minimal d'une $k$-surface est totalement d{\'e}termin{\'e} par son image, et
nous ommetrons souvent de parler de l'immersion sous-jacente.

L'{\it espace des $k$-surfaces} ${\cal N}$ d'une vari{\'e}t{\'e} $N$ est l'espace des paires $(S,x)$, o{\`u} $ x\in S$, et o{\`u} $S$ est 
soit 
un tube, soit le repr{\'e}sentant minimal $k$-surface sans bord {\'e}ventuellement d{\'e}gen{\'e}r{\'e}e.

Cet espace ${\cal N}$ est munie d'une topologie d{\'e}crite dans \cite{L1}, pour laquelle 
la notion de convergence est celle des sous-vari{\'e}t{\'e}s immerg{\'e}es point{\'e}es d{\'e}crites en \ref{}.

Nous avons montr{\'e} que dans la section 8.1 de \cite{L1} que l'espace ${\cal N}$ {\'e}tait compact et 
muni d'une structure 
de lamination, c'est {\`a} dire de produit local. Une feuille ${\cal L}_S$ de cette lamination est d{\'e}crite par $S$ qui est soit
le repr{\'e}sentant minimal d'une $k$-surface, soit un tube ; les points de cette feuille 
${\cal L}_S$    sont alors tous les points de la forme $(S,x)$ o{\`u} $x$ 
d{\'e}crit les points de $S$.
%%%%%%%%%%%%%%%%%%%%%%%%%%%%%%12%%%%%%%%%%%%%%%%%%%%%%
\section{Densit{\'e} des feuilles p{\'e}riodiques}

A partir de maintenant, $N$ d{\'e}signera une vari{\'e}t{\'e} compacte {\`a} courbure
strictement n{\'e}gative et $M$ sera son rev{\^e}tement universel.

Soit ${\cal N}$ l'espace compact lamin{\'e} associ{\'e} d{\'e}crit en \ref{laminations}. Nous voulons montrer 

\theo{densper}{L'ensemble de feuilles compactes est dense dans ${\cal N}$}

Les feuilles compactes de ${\cal N}$ sont par d{\'e}finition les immersions
localement convexes d'une surface compacte dans $N$, et  dont les images
sont des $k$-surfaces. Nous pouvons
interpr{\'e}ter ces surfaces compactes dans le rev{\^e}tement universel de la
mani{\`e}re suivante. Elles s'identifient aux quadruplets $(f,\S, \G,\r )$, o{\`u} 
\begitem
\iti $\G$ est un groupe discret agissant de mani{\`e}re cocompacte sur
$\S$ ;
\itii $\r$ est une repr{\'e}sentation de $\G$ dans $\pi_{1}(N)$ ;
\itiii $f$ est une immersion $\r$-{\'e}quivariante, localement convexe du disque $\S$ dans
$M$,  et  dont l'image est une $k$-surface.

\enditem
Introduisons une autre d{\'e}finition, un quadruplet $(f,\S ,\G,\r )$
est une {\it immersion {\'e}quivariante} {\`a} valeurs dans $\bo_{\infty} M$, si
$\S$ est une surface, {\'e}ventuellement {\`a} bord, $\G$ est un groupe agissant sur $\S$ de facon {\`a} ce
que le quotient soit une surface, $f$ est un
hom{\'e}omorphisme local de $\S$ dans $\bo_{\infty}M$, $\r$ est une
repr{\'e}sentation de $\G$ dans $iso(M)$ le groupe des isom{\'e}tries de
$M$. L'ensemble devant v{\'e}rifier la relation d'{\'e}quivariance \cad
$$
\forall s\in \S , \forall g\in\G , f(g.s)=\r (g).f(s).
$$
Une immersion {\'e}quivariante est {\it cocompacte} si $\S/\G$ est une
surface compacte.

D'apr{\`e}s le th{\'e}or{\`e}me \ref{equib}, les feuilles compactes de notre
lamination s'i\-den\-ti\-fient aux immersions {\'e}quivariantes cocompactes.

Il n'est peut-{\^e}tre pas inutile de remarquer que dans la d{\'e}monstration que nous allons donner, $\r (\G )$
sera toujours un groupe libre.

\bigskip 
La structure de cette section est la suivante.

Dans un premier
paragraphe ( \ref{parinclu} .  ~),  nous introduisons une d{\'e}finition, l'inclusion des donn{\'e}es
asymptotiques, et d{\'e}montrons un lemme de compacit{\'e} s'y rapportant.

Nous expliquons ensuite  ( en \ref{densper1} .  ) comment associer {\`a} un plongement du disque une
suite d'immersions {\'e}quivariantes l'approximant en un certain
sens.

Le paragraphe suivant ( \ref{densper2} . ~) montre comment construire de nouvelles
immersions {\'e}quivariantes {\`a} partir d'anciennes par un proc{\'e}d{\'e} de fusion.

Ensuite ( en \ref{densper3} . ~), nous introduisons une classe d'immersions topologiques du disque
dans les sph{\`e}res que nous appelons {\it quasi-plongement} et commen{\c c}ons
par montrer que toute donn{\'e}e asymptotique est limite de quasi-plongements.

En \ref{densper4} .  , nous montrons que tout quasi-plongement est limite de 
d'im\-mer\-sions {\'e}quivariantes cocompactes : celles-ci sont obtenues par
fusion d'immersions {\'e}quivariantes construites en \ref{densper1} .

Enfin, nous concluons la d{\'e}monstration dans le dernier paragraphe.

\subsection{Inclusion, lemme de compacit{\'e} pour les inclusions}\label{parinclu} Dans ce qui suit,
$D_{i}$ d{\'e}signera toujours une vari{\'e}t{\'e} connexe de dimension 2 sans bord. 
\subsubsection{Inclusion}
Introduisons encore une d{\'e}finition :
si $(f,D_{1})$ et $(g,D_{2})$ sont deux donn{\'e}es de probl{\`e}mes de Plateau asymptotiques, nous
dirons que $(f,D_{1})$ est {\it  inclus} dans $(g,D_{2})$ par $i$ , et
noterons
$$
(f,D_{1})\subset
(g,D_{2})~~[     i]
$$
 s'il existe un plongement $i$ de $D_{1}$ dans $D_{2}$, tel que $f=g\circ i$.
\par
La petite observation suivante\label{obsi} nous sera utile : 
\begitem
\iti soit $x_{1}$ et $x_{2}$
des points de $D_{1}$ et $D_{2}$ respectivement, si pour tout ouvert $U$
relativement compact de $D_{1}$, nous pouvons trouver $i_{U}$ telle que  $(f,U)$ est inclus dans $(g,D_{2})$ par
$i_{U}$ de telle sorte que $i_{U}(x_{1})=x_{2}$, alors
$(f,D_{1})$ est inclus dans $(g,D_{2})$.
\enditem
\par
\subsubsection{Intersection} 
Plus g{\'e}n{\'e}ralement, si $\{(f_{n},D_{n})\inn$ est  une suite de probl{\`e}mes de Plateau
asymptotiques, si $(g,D_{\infty})$ est un autre probl{\`e}me de Plateau asymptotique,
 inclus dans
tous les  $(f_{n},D_{n})$ par $i_{n}$, nous dirons que $(g,D_{\infty})$ est {\it l'intersection
des } $(f_{n},D_{n})$ {\it le long de} $i_{n}$, et nous noterons
$$
(g,D)=\bigcap_{n\in {\mathbb N}}(f_{n},D)~~ [ i_{n}]
$$
si toute donn{\'e}e $(h,U)$ ( o{\`u} $U$ est connexe), 
 incluse pour tout $n$ dans $(f_{n},D_{n})$ pour tout $n$ par $j_{n}$ et
 v{\'e}rifiant de plus 
$$
(g,D_{\infty})\subset (h,U)~~[j]
$$
o{\`u} $j$ v{\'e}rifie $i_{n}=j_{n}\circ j$, alors $(h,U)=(g,D)$.
On v{\'e}rifie  la proposition suivante

\pro{obsii}{Supposons que 
$$
(g,D)=\bigcap_{n\in {\mathbb N}}(f_{n},D)~~ [ i_{n}].
$$
Soit par ailleurs, $(h,U)$ telle que, pour tout $n$ 
$$
(h,U)\subset (f_{n},D) ~~[ j_{n}].
$$
Supposons qu'il existe $x_{0}\in D$ et $x_{1}\in U$ tels que pour tout
$n$,  $i_{n}(x_{0})=j_{n}(x_{1})$. 
Alors
$$
(h,U)\subset (g,D) ~~[ j],
$$
de telle sorte que $j_{n}=i_{n}\circ j$.}

\subsubsection{Limite inf{\'e}rieure}
Enfin, nous dirons que  $(g,D_{\infty})$ est {\it la limite inf{\'e}rieure
des } $(f_{n},D_{n})$ par $i_{n}$ et nous noterons
$$
(g,D_{\infty})=\liminf_{n\in {\mathbb N}}(f_{n},D_{n})~~ [ i_{n}],
$$
si pour toutes les sous-suites $s(n)$, il existe $p$ tel que
$$
(g,D_{\infty})=\bigcap_{n>p}(f_{s(n)},D_{s(n)})~~ [ i_{s(n)}].
$$

Il faut faire un petit peu attention avec ces notions : si en effet
nous con\-si\-d{\'e}\-rons la suite d{\'e}croissante d'ouvert $\{U_{n}\inn$ obtenue en enlevant les
uns apr{\`e}s les autres les points d'un ensemble d{\'e}nombrable dont
l'adh{\'e}rence est un {\'e}quateur, et si $f_{n}$ est la repr{\'e}sentation
conforme envoyant z{\'e}ro sur le p{\^o}le nord, nous avons, suivant le choix de
nos intersections deux intersections des
$\{(f_{n},U_{n})\inn$ : la repr{\'e}sentation de l'h{\'e}misph{\`e}re nord et celle
de l'h{\'e}misph{\`e}re sud.

\subsubsection{Encore un lemme de compacit{\'e}}
 
Nous voulons montrer le

\lem{inclu}{Soit $\{(f_{n},D)\inn$ une suite de probl{\`e}mes de Plateau
asymptotiques, admettant des solutions $\{\phi_{n}\inn$ et $x_{0}$ un
point de $D$; soit
{\'e}galement une donn{\'e}e asymptotique 
$$
(g,D)=\liminf_{n\in\mathbb N}
(f_{n},D_{n})~~ [ i_{n}].
$$
Alors, la
suite de $k$-surfaces $\{(\phi_{n},D)\inn$, point{\'e}es en $\{i_{n}(x_{0}) \inn$, converge
vers une solution du probl{\`e}me de Plateau asymptotique pour  $(g,D)$.}

\preu pour simplifier, nous nous noterons $0=x_{0}=i_{n}(x_{0})$.
 Soit  $\psi$ la solution du probl{\`e}me de Plateau asymptotique
d{\'e}fini par $(g,D)$. Cette solution existe bien  d'apr{\`e}s le th{\'e}or{\`e}me \ref{pasequi}.

Nous savons  ({\it cf} \ref{uniasym3}) que, pour tout $n$, $\psi(D)$ est un
graphe au dessus d'un ouvert de $\phi_{n}(D)$.

Notre premi{\`e}re
{\'e}tape consiste {\`a} d{\'e}montrer l'assertion :
\smallskip
\begitem
\item[(1)]la suite $\{\phi_{n}(0)\inn$ reste {\`a} distance born{\'e}e de
$\psi(0)$.
\enditem
\smallskip
Raisonnons par l'absurde. 
Dans le cas contraire, il existe une sous-suite $s(n)$ telle que
$\phi_{s(n)}(0)$ converge vers un point $m$ de $\bo_{\infty}M$.

Alors pour tout voisinage $U$ de $m$, pour $n$ suffisamment grand, si $N_{n}$
d{\'e}signe l'application de Gauss-Minkowski de $\phi_{s(n)}$ et $B_{n}$ la
boule de rayon 1 de centre $\phi_{s(n)}(0)$ trac{\'e}e sur $\phi_{s(n)}(D)$,
nous avons 
$$
\bo_{\infty}M\setminus U\subset N_{n}(B_{n}).
$$

En particulier,  $(i,\bo_{\infty}M\setminus U)$  est 
inclus dans les $(f_{s(n)},D)$ pour tous les $n$ suffisamment grand et
ainsi dans $(g,D)$, par \ref{obsii}.

Maintenant, d'apr{\`e}s notre observation pr{\'e}liminaire (\ref{obsi}), nous pouvons
con\-clu\-re 
que $(i,\bo_{\infty}M\setminus \{m\})$ est
inclus dans $(g,D)$. Nous avons l{\`a} notre contradiction, car le th{\'e}or{\`e}me
 \ref{pasequi}
fournirait alors une solution du probl{\`e}me de Plateau asymptotique pour
$(i,\bo_{\infty}M\setminus \{m\})$, ce qui est impossible d'apr{\`e}s
\ref{inex}.

Nous avons fini de d{\'e}montrer l'assertion (1).

Nous sommes maintenant en mesure d'appliquer notre th{\'e}or{\`e}me de compacit{\'e}
\ref{comrap} dans notre cas particulier o{\`u} il n'y a pas de condition
sur le bord. Extrayons  donc une sous-suite $s(n)$ telle que nous ayions
l'alternative :
\begitem
\item[(a)] soit $\{\phi_{s(n)}(D),\phi_{s(n)}(0)\inn$ converge vers une
$k$-surface (peut-{\^e}tre d{\'e}\-g{\'e}\-n{\'e}\-r{\'e}e) $S$,

\item[(b)] soit  $\{n\circ\phi_{s(n)}(D),n\circ\phi_{s(n)}(0)\inn$
converge vers un tube (n{\'e}cessairement complet) autour d'une g{\'e}od{\'e}sique
$\gamma$ joignant deux points {\`a} l'infini $\g_{+}$ et $\g_{-}$. Ici, bien s{\^u}r $n$
d{\'e}signe, l'application de Gauss de la surface (ici, $\phi_{n}(D)$)  {\`a} valeurs dans le fibr{\'e}
unitaire de $M$.
\enditem
Dans les deux cas de cette alternative, nous obtenons naturellement
associ{\'e}e {\`a} la limite une donn{\'e}e de probl{\`e}me de Plateau asymptotique $(f,D)$ :
l'application de Gauss-Minkowski de  la surface limite dans le premier cas, le
rev{\^e}tement universel de la sph{\`e}re moins deux points dans l'autre cas.

Par convergence uniforme sur tout compact, nous en d{\'e}duisons que pour
tout ouvert relativement compact $U$ de $D$, nous obtenons que $(f,U)$ est
inclus dans $(f_{s(n)},D)$, pour tout $n$ suffisamment grand et donc
dans $(g,D)$, par d{\'e}finition de limite inf{\'e}rieure.

Gr{\^a}ce {\`a} notre observation pr{\'e}liminaire, nous en d{\'e}duisons que $(f,D)$ est
inclus dans $(g,D)$.

En particulier, le deuxi{\`e}me cas de l'alternative
est exclus : en effet, le th{\'e}or{\`e}me \ref{pasequi} entrainerait alors une
solution du probl{\`e}me de Plateau d{\'e}fini par la sph{\`e}re moins deux points,
ce qui est impossible par \ref{inex}.

Enfin, prenons un ouvert relativement compact $U$ de $D$. Par
convergence uniforme sur tout compact, nous en d{\'e}duisons que $(g,U)$
est inclus dans $(f,D)$. Ainsi, $(g,D)$ est inclus dans $(f,D)$.

En conclusion, $(g,D)=(f,D)$ et, par unicit{\'e} de la solution d'un probl{\`e}me
de Plateau asymptotique \ref{uniasym}, nous en d{\'e}duisons que
$\{\phi_{s(n)}\inn$ converge vers $\psi$. $\diamond$

\subsection{Construction de groupes}\label{densper1} Nous allons dans un premier
temps exhiber des groupes associ{\'e}s {\`a} des plongements du disque.
 
Pour tout $\g\in\pi_{1}(N)$ diff{\'e}rent de l'identit{\'e}, $\g^{+}$ et $\g^{-}$ d{\'e}signeront respectivement les points
attracteurs et r{\'e}pulseurs de $\g$ sur $\bo_{\infty
}M$.

Si $F$ est un sous-groupe de $\pi_{1}(N)$, nous noterons $\L(F)$ son
ensemble limite et $\Omega (F)=\bo_{\infty}\setminus\L (F)$ son ensemble de
discontinuit{\'e}. Nous dirons que $F$ est {\it convexe cocompact} si
$F$ agit de mani{\`e}re cocompacte sur $\Omega(F)$.

La proposition que nous allons {\'e}noncer va {\^e}tre assez technique.

\pro{ping}{Donnons nous $(f,D)$ un plongement du disque dans  
$\bo_{\infty}M$ et $x_{i}^{\pm}$ des points de $f(\partial D)$. Soit maintenant de plus,   $\{\g_{i,n}\inn$,   $i\in\{1,\ldots,
  p\}$,
  $p$ suites d'{\'e}l{\'e}ments de $\pi_{1}(N)$,
telles que 
\begitem
\iti les suites $\{\g^{\pm}_{i,n}\inn$ convergent vers $x_{i}^{\pm}$,
\itii $\forall i,n,~~\g_{i,n}^{\pm}\notin f(D)$.
\enditem
Il existe alors une suite de groupes libres convexes cocompacts
$\{F_{n}\inn$
v{\'e}rifiant
$$
\forall n , ~~\exists q,~~\forall i ,~~ \g^{q}_{i,n}\in F_{n},\leqno(1)
$$
$$
\forall \g\in F_{n}\setminus\{{\rm id}\},~~ \g (f(D))\cap f(D) =\emptyset\leqno(2)
$$
$$
(f,D)= \liminf_{n\in{\mathbb N}} (i_{n},\Omega (F_{n}) )~~[f\circ i_{n}^{-1}],\leqno(3)
$$
o{\`u} $i_{n}$ d{\'e}signe l'injection canonique  de $\Omega (F_{n})$ dans $\bo_{\infty}(M)$}

\preu munissons $\bo_{\infty}M$ d'une m{\'e}trique arbitraire, fixons un
entier $n$ et soit $\e$ le r{\'e}el strictement positif ( {\`a} cause de (ii) ) d{\'e}fini par
$$
\e=Sup\{ d(\g_{i,n}^{\pm},x_{i}^{\pm})\}.
$$
Par hypoth{\`e}se $\e$ tend vers 0 quand $n$ tend vers l'infini.

Par densit{\'e} des g{\'e}od{\'e}siques
p{\'e}riodiques, il existe des {\'e}l{\'e}ments de $\pi_{1}(N)$,
$\{
\l_{1},\ldots\l_{k}\}$, o{\`u}  $k$ est un nombre fini d{\'e}pendant de $n$,
tels que l'ensemble 
$$
L=\{\l_{i}^{+},\l_{i}^{-} ; 1\leq i\leq k\}
$$
soit $\e$-proche pour la distance de Haussdorff de $f(\bo D)$ et {\`a}
l'ext{\'e}rieur de $f(D)$. On peut
enfin s'arranger pour que tous les points de $L$ soient distincts et que
 $l_{i}=\g_{i,n}$, $\forall i\in \{1,\ldots,p\}$.

Par l'argument ping-pong classique, il existe $q$ suffisamment grand tel
que le groupe 
$$
F_{n}=\langle \l_{1}^{q},\ldots,\l_{k}^{q}\rangle
$$
soit libre et agisse de mani{\`e}re cocompacte sur son ensemble de discontinuit{\'e}.

On peut de plus s'arranger en prenant $q$ suffisamment grand  pour que son
ensemble  limite $\L_{n}$ soit dans un $\e$
voisinage de $L$, et {\`a} l'ext{\'e}rieur de $f(D)$. 

Toujours en prenant $q$ suffisamment grand, on obtient la condition (2).

Pour achever la d{\'e}monstration de la proposition, il nous suffit de
 remarquer que $\L_{n}$ est une suite de compacts convergeant pour la
 m{\'e}trique de Haussdorff vers le bord de $f(D)$, tout en restant {\`a}
 l'ext{\'e}rieur de $f(D)$, et qu'ainsi nous avons bien 
$$
(f,D)=\liminf_{n\in {\mathbb N}} (i_{n},\Omega (F_{n}) )~~[f\circ i_{n}^{-1}],
$$
ce que nous voulions d{\'e}montrer
$\diamond$
\par

Nous avons d{\'e}j{\`a} tout en main pour montrer que chaque plongement du
disque est limite d'immersions {\'e}quivariantes cocompactes : il nous
suffit de mettre ensemble la proposition pr{\'e}c{\'e}dente et le lemme de
compacit{\'e}
\ref{inclu}

Malheureusement, toute immersion n'est pas limite de plongement. Il nous
faut donc raffiner notre construction et donner une proc{\'e}dure, que nous appelerons
{\it fusion},  permettant
de construire de nouvelles immersions {\'e}quivariantes cocompactes, et
montrer que celle ci permet d'approximer ce que nous appelerons des
{\it quasi-plongements} qui eux vont s'av{\'e}rer  dense dans les
immersions.

%%%%%%%%%%%%%%%%%%%%%%%%%%%%%%%%%%%13%%%%%%%%%%%%%%%%%%%%%%%%%%%%

\subsection{Fusion de groupes et de surfaces}\label{densper2} Nous allons
expliquer dans ce paragraphe une construction
permettant de construire de nouveaux probl{\`e}mes de Plateau {\'e}quivariants {\`a}
partir d'anciens.

La proc{\'e}dure que nous allons d{\'e}crire est donn{\'e}e par la construction
suivante. Introduisons tout d'abord nos notations et hypoth{\`e}ses.

Soient $(f_{i},S_{i},\G_{i},\r_{i})$, o{\`u} $i\in\{1,2\}$, deux
immersions {\'e}quivariantes. Nous supposerons ici que  les $S_{i}/\G_{i}$
sont des surfaces {\`a} bord. Soient  $c_{1}$ et  $c_{2}$  deux composantes connexes
du bord respectivement de  $S_{1}/\G_{1}$ et
$S_{2}/\G_{2}$ . On suppose qu'il existe deux relev{\'e}s $v_{i}$ de ces
courbes dans $S_{i}$, associ{\'e}s {\`a} des {\'e}l{\'e}ments $\g_{i}\in\G_{i}$ tels que
$\r_{1}(\g_{1})=\r_{2}(\g_{2})$  ainsi qu'un hom{\'e}omorphisme  $j$ entre $v_{1}$ et $v_{2}$, v{\'e}rifiant

$$
\forall s\in v_{2},~~f_{1}\circ j (s)=f_{2} (s)
$$
ainsi que

$$
\g_{2}\circ j=j\circ \g_{1}.
$$
En particulier $j$ descend en un hom{\'e}omorphisme $h$ entre $c_{1}$ 
et $c_{2}$.
Notre construction d{\'e}coule alors de l'imm{\'e}diate proposition suivante :

\pro{prefusion}{Avec les hypoth{\`e}ses et notations pr{\'e}c{\'e}dentes, il existe 
 une unique immersion {\'e}quivariante $(f, \S,\G,\r)$ telle
que $\S/\G=\S_{1}/\G_{1}\cup_{h}\S/\G_{2}$ et, si $\iota_{i}$ est l'injection
de $\S_{i}$ dans $\S$ qui se d{\'e}duit de cette identification, alors
$f\circ \iota_{i}=f_{i}$. }

\preu il s'agit d'une construction standard $\diamond$

L'immersion {\'e}quivariante ainsi construite sera appel{\'e}e {\it fusion} des
deux immersions {\'e}quivariantes pr{\'e}c{\'e}dentes.

On peut remarquer alors que $\G$ est un produit amalgamm{\'e} de $\G_{1}$ et
$\G_{2}$

Nous verrons plus tard comment la fusion se  comporte vis-a-vis de
l'o\-p{\'e}\-ra\-tion limite inf{\'e}rieure. 

\subsection{Quasi-plongements}\label{densper3} Introduisons une d{\'e}finition interm{\'e}diaire.

Un hom{\'e}omorphisme local $f$ du disque ouvert $D$ dans $\bo_{\infty}M$ sera
appel{\'e} un {\it quasi-plongement}, si $f$ se prolonge continuement en un
hom{\'e}omorphisme local de l'adh{\'e}rence de $D$, et s'il existe une famille
finie d'arcs plong{\'e}s $\{c_{i}\}$, $1\leq i\leq q $ deux {\`a} deux disjoints, dont les extr{\'e}mit{\'e}s
sont dans $\bo D$, et telle que $f$ s'{\'e}tend en un plongement de l'adh{\'e}rence
de chaque composante connexe de 

$$
D\setminus\bigcup_{1\leq i\leq q}c_{i}.
$$
Nous appelerons la famille d'arcs de la d{\'e}finition, {\it d{\'e}coupe} du quasi-plongement.

En un certain sens, le r{\'e}sultat suivant affirme que tout hom{\'e}omorphisme local est limite de quasi-plongements.

\pro{quasiplong}{Soit $f$ un hom{\'e}omorphisme local de $D$ dans
$\bo_{\infty}M$, alors pour tout ouvert relativement compact $U$ de $D$
il existe un ouvert relativement compact $V$ inclus dans $D$ contenant $U$ tel que
$(f,V)$ soit un quasi-plongement}

\preu nous allons utiliser des id{\'e}es contenues dans la param{\'e}trisation
faite par  Thuston de l'espace des ${\mathbb CP}^{1}$-structures par les laminations g{\'e}od{\'e}siques
mesur{\'e}es. Notre m{\'e}thode va {\^e}tre laborieuse, et nous aimerions savoir s'il
y a plus simple.

Identifions $\bo_{\infty}M$ {\`a} la sph{\`e}re $S^{2}$. Munissons $S^{2}$ d'une
m{\'e}trique {\`a} courbure constante  et $D$ de la
m{\'e}trique induite par $f$.  Nous noterons ${\bar D}$ sa compl{\'e}tion
m{\'e}trique et $Fr(D)={\bar D}\setminus D$. Au besoin, en restreignant un
peu $D$, nous pouvons supposer que $\bar D$ est hom{\'e}omorphe au disque
ferm{\'e} et que $f$ s'{\'e}tend en un hom{\'e}omorphisme local de ${\bar D}$ dans $S^2$.

Une boule m{\'e}trique ouverte  de $D$ sera 
appel{\'e}e  une {\it bonne boule}, si $f$ est une isom{\'e}trie de cette boule
sur une boule de $S^{2}$. Une boule sera {\it maximale}, si c'est une
bonne boule et si elle n'est incluse dans aucune autre bonne boule qu'elle
m{\^e}me.

Si $B$ est une boule maximale, nous noterons ${\bar B}$ son adh{\'e}rence
dans ${ D}$ (et non dans ${\bar D}$) et $Fr(B)$ l'ensemble ${\bar B}\setminus B$. Cet
ensemble $Fr(B)$ est une r{\'e}union 
d'intervalles ouverts. Un peu de g{\'e}om{\'e}trie sph{\'e}rique montre que si deux
points de $Fr(B)$ sont dans une m{\^e}me composante connexe de $Fr(B)$,
alors 
ils sont
inclus dans une m{\^e}me bonne boule.

Nous dirons enfin que deux points $x$ et $y$ de $Fr(D)$ sont {\it joignables},
s'il existe une boule maximale $B$ et une composante connexe $A$ de $Fr(B)$
telle que $x$ et $y$ appartiennent tous deux {\`a} l'adh{\'e}rence de $A$
dans ${\bar D}$.

Reamarquons que  si $x_{1}$ et $x_{2}$ sont
joignables de m{\^e}me que $y_{1}$ et $y_{2}$ avec $y_{i}\not= x_{j}$, alors
$y_{1}$ et $y_{2}$ sont tous les deux inclus dans la m{\^e}me composante
connexe de $Fr(D)\setminus \{x_{1},x_{2}\}$. Explicitons ce
raisonnement. Notons $B_{x}$ et $B_{y}$ les bonnes boules respectives
pour les paires $(x_{1},x_{2})$ et $(y_{1},y_{2})$, ainsi que $A_{x}$ et
$A_{y}$ les composantes connexes de $Fr(B_{x})$ et $Fr(B_{y})$ qui s'en 
d{\'e}duisent. Il est calir que $A_{x}$ et $A_{y}$ {\'e}tant des portions de
cercles  ne peuvent s'intersecter en une infinit{\'e} de points ar lles
seraient alors confondues. De plus, $y_{1}$ et $y_{2}$ sont dans des  composantes
connexes diff{\'e}rentes de $Fr(D)\setminus \{x_{1},x_{2}\}$, nous en
d{\'e}duisons que $A_{x}$ et $A_{y}$ se rencontrent transversalement en
exactement un point, puisque deux cercles ne peuvent se renconter en
plus de deux points. Mais alors, l'une des extr{\'e}mit{\'e}s de $A_{x}$ (
c'est-{\`a}-dire $x_{1}$ et $x_{2}$ ) serait
inclus dans l'int{\'e}rieur de $B_{y}$ d'o{\`u} la contradiction.

Si maintenant $x$ et $y$ sont deux points joignables correspondant {\`a} une
boule maximale $B$ (pas n{\'e}cessairement unique d'ailleurs), tra{\c
c}ons entre $x$ et $y$ la g{\'e}od{\'e}sique pour la m{\'e}trique hyperbolique
conforme de $B$. L'arc correspondant sera appel{\'e}e {\it bon arc}. Un
raisonnement g{\'e}om{\'e}trique {\'e}l{\'e}mentaire montre que deux bons arcs ne
peuvent s'intersecter qu'en leurs extr{\'e}mit{\'e}s.

La r{\'e}union ${\cal L}$ des bons arcs forme un ensemble ferm{\'e} et nous avons construit
ainsi une lamination de $D$.

Cette construction n'est bien s{\^u}r pas originale, c'est celle qui associe
{\`a} toute ${\mathbb CP}^1$-structure une lamination comme dans la construction de
Thurston d{\'e}j{\`a} cit{\'e}e.

On remarque ais{\'e}ment que si $A$ est le compl{\'e}mentaire de ${\cal L}$,
toute composante connexe de $A$ est inclus dans une bonne boule.

Nous pouvons maintenant construire une exhaustion de $D$ par des disques
compacts $D_{n}$ telle que toute intersection non vide d'un bon arc avec $D_{n}$ soit connexe
et de longueur minor{\'e}e par une constante ne d{\'e}pendant que de $n$ : par exemple, on peut trouver une m{\'e}trique
hyperbolique telle que la lamination soit g{\'e}od{\'e}sique, puis prendre une
exhaustion par des convexes, et enfin {\'e}ventuellement d{\'e}couper ces
convexes en enlevant la ``petite''composante connexe du compl{\'e}mentaire  des bons arcs de petite longueur.

Fixons ensuite l'un des $D_{p}$ ; nous voulons montrer que $f$ restreinte
{\`a} $D_{p}$ est un quasi-plongement. 

Pour cela, il nous suffit de contruire un recouvrement de $D_{p}$ tel
que tout ouvert de ce recouvrement 
est bord{\'e} par une courbe de Jordan contitu{\'e}e d'un nombre fini de bons arcs et d'arcs
inclus dans le bord de $D_{p}$, et tel qu'il existe une d{\'e}coupe de cet
ouvert par des bons arcs pour laquelle la restriction de $f$ {\`a} cet
ouvert soit un quasi-plongement. En effet, lorsque l'on aura extrait un
recouvrement fini par de tels ouverts, la famille de bons arcs inclus
dans le bord de ces ouverts et provenant des d{\'e}coupes de ces m{\^e}mes ouverts, nous donnera une d{\'e}coupe de quasi-plongement.

Ce recouvrement est ais{\'e} {\`a} construire : nous avons quatre cas {\`a} consid{\'e}rer.
 
 (1) Si $x\in A=D\setminus{\cal L}$, nous prenons
simplement la composante connexe de  $x$ dans $A\cap D_{p}$. Le longueur
de la trace des bons arcs sur $D_{p}$ {\'e}tant minor{\'e}e, il n'y a qu'un
nombre fini de bons arcs dans le bord de cette composante connexe.

 (2) Si  $x$
se trouve sur un bon arc isol{\'e}. Ce bon arc s{\'e}pare deux composantes connexes de $A$. Nous
prenons alors comme ouvert la r{\'e}union de ces composantes connexes, avec
la d{\'e}coupe donn{\'e} par le bon arc passant par $x$.

 (3) Si $x$ se trouve sur un bon arc $c$ associ{\'e} {\`a} une boule maximale $B$,  isol{\'e} d'un c{\^o}t{\'e} mais pas de
l'autre, c'est {\`a} dire si ce bon arc
borde une composante $O$ connexe de $A$,
et s'il existe une suite $\{c_{n}\inn$ d'arcs distincts de $c$ tendant vers
$c$. Dans ce cas,  nous pouvons remarquer que pour $n$ suffisamment grand, la r{\'e}gion
$U_{n}$ bord{\'e} par $c$ et $c_{n}$ est telle que $U_{n}\cap D_{p}$ est incluse dans
$B$. Nous pouvons alors prendre comme ouvert $(O\cup U_{n}\cup c)\cap D_{p}$
qui est inclus dans $B$.

 (4)  Le dernier cas correspond au cas $x$ appartient {\`a} un bon arc $c$ et o{\`u} nous pouvons trouver deux suites
de bons arcs tendant vers $c$,
$\{c_{n}\inn$ et $\{g_{n}\inn$, telles que la r{\'e}gion $U_{n}$
comprise entre $c_{n}$ et $g_{n}$ contienne $c$. Il nous suffit de
remarquer que pour $n$ suffisamment grand $U_{n}\cap D_{p}$ est inclus
dans $B$, et de prendre enfin comme ouvert $U_{n}\cap D_{p}$.

Notre laborieuse construction est termin{\'e}e. $\diamond$

%%%%%%%%%%%%%%%%%%%%%%%%%%%%%%%14%%%%%%%%%%%%%%%%%%%%%%%%%%
\subsection{Quasi-plongements et immersions {\'e}quivariantes
cocompactes}\label{densper4} 
Le r{\'e}sultat essentiel de cette section est une g{\'e}n{\'e}ralisation de
\ref{ping} au cas des quasi-plongements. 

\pro{densqp}{Soit $(f,D)$ un quasi-plongement, il existe alors une suite d'immersions
{\'e}quivariantes cocompactes $\{(f_{n},S_{n})\inn$ telles que
$$
(f,D)=\liminf_{n\in {\mathbb N}}(f_{n},S_{n})~~[i_{n}].
$$}

\preu soit donc $(f,D)$ un quasi-plongement. D'apr{\`e}s la d{\'e}finition, il existe une famille
finie d'arcs plong{\'e}s $\{c_{i}\}$, $1\leq i\leq q $ deux {\`a} deux disjoints, dont les extr{\'e}mit{\'e}s
sont dans $\bo D$, et telle que $f$ s'{\'e}tend en un plongement de l'adh{\'e}rence
de chaque composante connexe $D_{j}$ de 
$$
D_{c}=D\setminus\bigcup_{1\leq i\leq q}c_{i}.
$$ 

Orientons les arcs $c_{i}$, notons $x_{i}^{\pm}$ leurs extr{\'e}mit{\'e}s et
$D_{i^{\pm}}$ les composantes connexes de $D_{c}$ bord{\'e}es par $c$.

Par densit{\'e} des g{\'e}od{\'e}siques p{\'e}riodiques, nous pouvons trouver des suites
$\{\g_{i,n}\inn$ d'{\'e}l{\'e}ments de $\pi_{1}(N)$ telles que

\begitem
\iti les suites $\{\g^{\pm}_{i,n}\inn$ convergent vers
$f( x_{i}^{\pm})$
\itii $\forall~ i,n~~\g_{i,n}^{\pm}\notin f(D_{i}^{\pm}).$
\enditem
D'apr{\`e}s le lemme \ref{ping}, nous pouvons alors trouver pour toute
composante $D_{j}$ de $D_{c}$ une suite de  groupes libres convexes
cocompacts $\{F_{j,n}\inn$ 
v{\'e}rifiant, si $c_{i}$ est dans le bord de $D_{j}$,
$$
\forall n , ~~\exists q,~~ \g^{q}_{i,n}\in F_{j,n},\leqno{\rm~~~~(iii)}
$$
ainsi que
$$
\forall\g\in F_{j,n}\setminus\{id\},~~ \g(f(D_{j}))\cap
f(D_{j})=\emptyset~~\leqno{\rm~~~~ (iv)}
$$
et 
$$
(f,D_{j})= \liminf_{n\in{\mathbb N}} (p_{j,n}, D).\leqno{\rm~~~~ (v)}
$$

Ici $p_{j,n}$ d{\'e}signe la projection du rev{\^e}tement universel de $\Omega
(F_{n})$ identifi{\'e} {\`a} $D$ et v{\'e}rifiant $p_{n}(0)=y_{j}$, o{\`u} $y_{j}$ est
un point  arbitraire de $D_{j}$.

Pour tout $j,~n$,   notons alors  ${\bar
S}_{j,n}$  la surface compacte
$$
{\bar S}_{j,n}=\Omega(F_{j,n})/F_{j,n}
$$
et $\pi_{j,n}$ la projection naturelle de $\Omega(F_{j,n})$ sur  ${\bar
S}_{j,n}$.

La condition (iv) nous assure que $\pi_{j,n}\circ f$ est un plongement
de $D_{j}$ dans ${\bar S}_{j,n}$.

Il nous est toujours loisible dans la contruction pr{\'e}c{\'e}dente de prendre des
sous-groupes d'indice finis des $F_{j,n}$. Nous avons par ailleurs une
certaine libert{\'e} concernant le choix des g{\'e}n{\'e}rateurs des groupes
$F_{j,n}$ autres que les $\g_{i}$. 

Pour continuer notre d{\'e}monstration, nous allons raffiner notre
construction en utilisant la latitude que nous venons de d{\'e}crire de
telle sorte que les groupes 
$F_{j,n}$ v{\'e}rifient les propri{\'e}t{\'e}s suppl{\'e}mentaires suivantes :
\begitem
\itvi le groupe
$F_{n}$ engendr{\'e} par la r{\'e}union des $F_{j,n}$ est lui-m{\^e}me convexe
cocompact, 
\enditem
$$
f(c_{i})\subset \Omega (F_{n}).\leqno{\rm~~~~(vii)}
$$

\smallskip

Choisissons maintenant, pour tout $i$, $n$, des voisinages
$B^{\pm}_{i,n}$ de $\g^{\pm}_{i,n}$, hom{\'e}\-omor\-phes  {\`a} la boule ferm{\'e}e,
disjoints deux {\`a} deux, 
tels que
$$
\g_{i,n}(\bo_{\infty}M \setminus B_{i,n}^{-})\subset
B_{i,n}^{+}, \leqno{\rm ~~~~(viii)}
$$

$$
\lim_{n\rightarrow\infty}( diam (B^{\pm}_{i,n}))=0,\leqno{\rm ~~~~(ix)}
$$
et enfin, ces voisinages intersectent les images des $c_{i}$ exactement
en leurs extr{\'e}mit{\'e}s
$$
f(c_{i})\cap B_{i,n}^{\pm}=f(x_{i}^{\pm}).\leqno{\rm ~~~~(x)}
$$

A nouveau, si besoin {\'e}tait, nous avons pris des sous-groupes d'indices
finis des  $F_{j,n}$.

Nous avons presque fini  nos pr{\'e}liminaires. Pour  les terminer, nous
construisons  pour tout entier $n$, et tout $j\in\{1,\ldots , q\}$, une courbe $c_{i,n}$ trac{\'e}e sur
$\Omega (F_{n})$ telle que
$$
f(c_{i})\subset c_{n,i},\leqno{\rm~~~~(xi)}
$$

$$
\g_{i,n}(c_{i,n})=c_{i,n},\leqno{\rm ~~~~(xii)}
$$ 
et
$$ 
c_{n,i}\setminus f(c_{i})\subset B_{i,n}^{+}\cup B_{i,n}^{-}.\leqno{\rm
~~~~(xiii)}
$$

Pour ce faire, il faut proc{\'e}der de la mani{\`e}re suivante : nous savons par
(viii) que
$\g_{i,n}(f(c_{i}))\subset B_{i,n}^{+}$. Par ailleurs,
$B_{i,n}^{+}\setminus B_{i,n}^{-}$ est connexe de m{\^e}me que son adh{\'e}rence
$A$. Par (x), cette adh{\'e}rence $A$ intersecte $f(c_{i})$ exactement en
$f(x_{i}^{+})$ et $\g_{i,n}(f(c_{i}))$ en
$\g_{i,n}(f(x_{i}^{-}))$. Choisissons donc un chemin plong{\'e} $\l_{i,n}$
joignant $f(x_{i}^{+})$ et $\g_{i,n}(f(x_{i}^{-}))$ dans $A$. Nous pouvons
alors prendre comme courbe $c_{i,n}$, la r{\'e}union des images it{\'e}r{\'e}es de
l'arc $f(c_{i})\cup \l_{i,n}$ :
$$
c_{i,n}=\bigcup_{p\in{\mathbb Z}} \g_{i,n}^{p}(f(c_{i})\cup \l_{i,n}).
$$

\begitem

\item[(xiv)] Enfin,
en prenant {\'e}ventuellement des sous-groupes d'indice fini des
$F_{j,n}$, la condition (xiii) et le fait que les $B_{i,n}^{\pm}$
sont deux {\`a} deux disjoints,  nous pemettent d'assurer que,
pour tout $j$ et $n$,  si $c_{i_{1}}, \ldots
c_{i_{m}}$, sont trac{\'e}es dans le bord de $D_{j}$, alors 
les  courbes $\pi_{j,n}(c_{i_{k},n})$ sont des courbes plong{\'e}es et
d'intersection vide dans  ${\bar
S}_{j,n}=\Omega (F_{n})/F_{n}$.

\enditem
Rappelons qu'ici  $D_{j}$ est une composante connexe de
$$
D\setminus\bigcup_{1\leq i\leq q}c_{i}.
$$
\par
Nous avons maintenant fini d'imposer des conditions suppl{\'e}mentaires aux
groupes $F_{j,n}$ et pouvons continuer notre d{\'e}monstration.
\smallskip
Par la condition (iv), $\pi_{j,n}\circ f$ est un plongement de $D_{j}$
dans
$$
U_{j,n}={\bar S}_{j,n}\setminus \bigcup_{1\leq k\leq m}\pi_{j,n}(c_{i_{k},n}).
$$

Notons donc $\S_{j,n}$ la composante connexe de $U_{j,n}$ dans
laquelle s'envoie $D_{j}$. Par construction,  $\S_{j,n}$ est une
surface dont le bord s'identifie {\`a} la r{\'e}union des courbes
$\pi_{j,n}(c_{i_{k},n})$. Soit enfin $\G_{j,n}$ le groupe fondamental de
$\S_{j,n}$, et $S_{j,n}$ son rev{\^e}tement universel.

Notre {\'e}tape suivante va {\^e}tre de fusionner toutes les surfaces
$S_{j,n}$. La proc{\'e}dure est la suivante : pour tout $i$ et $n$, la
construction d{\'e}crite en \ref{densper2}, nous permet pr{\'e}\-ci\-sem\-ment de
fusionner $S_{n}^{i^{+}}$ et $S_{n}^{i^{-}}$ le long de
$c_{i,n}$. En proc{\'e}dant de proche en proche, nous
produisons
ainsi une immersion {\'e}quivariante cocompacte $(f_{n},S_{n},\G_{n})$ ayant
la propri{\'e}t{\'e} suivante : il existe un plongement $i_{n}$ de $D$ dans $S_{n}$
tel que $f=f_{n}\circ i_{n}$ \cad
$$
(f,D)\subset (f_{n},S_{n})~~[i_{n}].
$$

Par ailleurs les courbes $c_{i,n}$ donnent naissance {\`a} des courbes
$q_{i,n}$ trac{\'e}es sur $S_{n}$ et telles que $i_{n}(c_{i})\subset
q_{i,n}$ ainsi que $f_{n}(q_{i,n})=c_{i,n}$.

Pour conclure, nous voulons montrer

$$
(f,D)= \liminf_{n\in {\mathbb N}}(f_{n},S_{n})~~[i_{n}].
$$

Soit donc $(g,U)$ un plongement d'un ouvert $U$ connexe, tel que 
$$
(g,U)\subset (f_{s(n)},S_{s(n)})~~[j_{s(n)}],
$$
et 
$$
(f,D)\subset (h,U) ~~[j],
$$
avec $i_{n}=j_{n}\circ j$.

Pour simplifier les notations d{\'e}j{\`a} lourdes, nous allons supposer $s(n)=n$
ce qui est indolore.

Nous supposerons, en consid{\'e}rant  $j$ comme une inclusion, que $D$ est un ouvert de $U$
et notons $Z(D)$ la fronti{\`e}re de $D$ dans $U$. Nous voulons montrer que
$Z(D)$ est vide. Raisonnons par l'absurde et supposons que cet ensemble
est non vide.

Munissons $\bo_{\infty}M$ d'une m{\'e}trique annexe, $U$, $S_{n}$ et $D$ des
m{\'e}triques induites.

Par notre hypoth{\`e}se (ix),
$$
\lim_{n\rightarrow\infty}( diam (B^{\pm}_{i,n}))=0,
$$
ce qui nous permet de trouver un point $x$ de $Z(D)$ et un nombre $\e>0$
tel que  
$$
\forall i,n ~~d(g(x),B^{\pm}_{i,n})\geq\e.
$$
Nous pouvons {\'e}galement nous assurer que
$$
\forall i, ~d (x,c_{i})_geq\e
$$

et en particulier, {\`a} cause de (xiii)
$$
\forall i,n ~~d(i_{n}(x),q{i,n})\geq\e. \leqno{~~~~\rm (xv)}.
$$

En prenant $\e$ suffisamment petit, nous pouvons de plus assurer que $g$
est une bijection de la boule $B_{\e}^{U}$ de rayon $\e$ de centre $x$ dans $U$ et la
boule $B_{\e}$ de centre $g(x)$ de rayon $\e$ dans $\bo_{\infty}M$.

Le point $x$ appartient n{\'e}cessairement {\`a} la fronti{\`e}re d'une des
composantes
$D_{j}$ de
$$
D\setminus\bigcup_{1\leq i\leq q}c_{i}.
$$

Rappelons maintenant que $S_{n}$ est la r{\'e}union de copies isom{\'e}triques
des surfaces $S_{j,n}$, bord{\'e}es par les courbes $q_{i,n}$ qui s'envoient par
$f_{n}$ sur les courbes $c_{i,n}$.

La condition (xv) nous permet d'affirmer que $j_{n}(B_{\e}^{U})$ est
inclus dans l'int{\'e}rieur d'une de ces copies de $S_{j,n}$. En particulier,
puisque
$$
B_{\e}=g(B_{e}^{U})=f_{n}\circ j_{n}(B_{e}^{U}),
$$
nous en d{\'e}duisons que, pour tout $n$, $g(B_{\e})$ est inclus dans
$\Omega (F_{j,n})$, ce qui est impossible car la limite de Haussdorff des
$\L (F_{j,n})$ contient $Fr(D_{j})$, et en particulier $g(x)$, par
construction.

Ce dernier point est la contradiction recherch{\'e}e.$\diamond$

\subsection{D{\'e}monstration du th{\'e}or{\`e}me \ref{densper}} Nous pouvons maintenant
d{\'e}\-mon\-trer la densit{\'e} des feuilles p{\'e}riodiques dans l'espace lamin{\'e}
${\cal M}$.

Tout d'abord la proposition \ref{densqp} et le lemme \ref{inclu}
entra{\^\i}ne que l'a\-dh{\'e}\-ren\-ce de l'ensemble des feuilles p{\'e}riodiques contient
l'ensemble des quasi-plongements.

Le th{\'e}or{\`e}me \ref{densper} suit alors de la proposition
\ref{quasiplong} et de la proposition suivante

\pro{toto}{Soit $(f,D)$ une donn{\'e}e asymptotique, soit $\{U\inn$ une suite
em\-bo{\^\i}\-t{\'e}e d'ouvert
relativement compacts de $D$ tel que
$$
\bigcup_{n\in{\mathbb N}}U_{n}=D.
$$ Notons $\{(\phi_{n},U_{n})\inn$ la suite de solutions des probl{\`e}mes
de Plateau asymptotiques d{\'e}finis par  $\{(f,U_{n})\}_{n\in{\mathbb N}}$ et
obtenues par le th{\'e}or{\`e}me \ref{pasequi}. Nous avons alors les deux
r{\'e}sultats suivants :
\begitem
\iti si $(f,D)$ admet une solution $(\phi ,D)$ alors la suite $\{(\phi_{n},U_{n})\inn$ 
point{\'e}es en l'origine de $D$ converge vers  $(\phi ,D)$.
\itii si $D$ est le rev{\^e}tement universel $\bo_{\infty}M$ auquel on a {\^o}t{\'e}
les deux ex\-tr{\'e}\-mi\-t{\'e}s d'une g{\'e}od{\'e}sique $\g$, et si $f$ est la projection
canonique, alors la suite de surfaces $\{S_{n}\inn$, immerg{\'e}es dans $UM$
-- o{\`u} $S_{n}$ est l'ensemble des vecteurs normaux ext{\'e}rieurs {\`a} la surface
 $\phi_{n}( U_{n})$ -- point{\'e}es en l'origine de $D$, converge vers le
tube de $\g$.
\enditem
 }

\preu d{\'e}montrons (i) tout d'abord. Nous allons ressasser  nos
id{\'e}es habituelles. D'apr{\`e}s \ref{uniasym}, pour tout $n$,  $\phi_{n}( U_{n})$
est un graphe d'une fonction $f_{n}$ au dessus de $\phi(U_{n})$ . D'apr{\'e}s le lemme
\ref{deltahyper}
$$
\forall x\in D,~\exists K, ~s.t.~\forall n,~~~f_{n}(x)\leq K.
$$ 

Nos arguments de compacit{\'e} montre que la suite de fonctions $f_{n}$
converge alors $\ci$ sur tout compact. Le graphe de la limite est alors
n{\'e}cessairement une solution du probl{\`e}me de Plateau asymptotique d{\'e}fini
par $(f,D)$. L'unicit{\'e} de la solution du probl{\`e}me de Plateau
asymptotique, montre que cette fonction limite est nulle et donc (i).

En ce qui concerne (ii), il nous suffit de montrer que quelle que soit
la sous-suite  $s(n)$, $\{(\phi_{s(n)},U_{s(n)})\inn$ 
point{\'e}es en l'origine de $D$ ne converge pas vers une $k$-surface. En
effet dans ce cas, la suite  $\{S_{n}\inn$ converge vers le tube d'une
g{\'e}od{\'e}sique d'apr{\`e}s \ref{comrap}, et par construction de $D$, cette
g{\'e}od{\'e}sique ne peut {\^e}tre que $\g$.
Or si  $\{(\phi_{s(n)},U_{s(n)})\inn$ 
point{\'e}es en l'origine de $D$ ne convergeait vers une $k$-surface, celle
ci seriat une solution du probleme de Plateau asymptotique pour
$\bo_{\infty}M$ auquel on {\^o}t{\'e} un ou deux points  et d'apr{\`e}s \ref{inex}
une telle solution n'existe pas.$\diamond$

\subsubsection{Remarque}{\label{tutu}} Si $(f,D)$ est une donn{\'e}e
asymptotique qui s'{\'e}tend en un homeomorphisme  d'un voisinage du disque ferm{\'e}
dans $\bo_{\infty}M$, et $(\phi ,D)$ est une $k$-surface solution du probl{\`e}me de Plateau
asymptotique pour $(f,D)$, alors $(\phi,D)$ n'est pas
tubulaire {\`a} l'infini. En effet, si $(\phi,D)$ {\'e}tait tubulaire {\`a} l'infini,
alors il existerait une suite $\xnn$ de points de $\bo_{\infty}M$, telle
que le cardinal de $f^{-1} (x_{n})$ tende vers l'infini.

%%%%%%%%%%%%%%%%%%%%%%%%%%%%%%%%%15%%%%%%%%%%%%%%%%%%%%%%%%%%%%%%%
\subsection{Feuilles p{\'e}riodiques de m{\^e}me genre}

Nous voulons d{\'e}montrer le 

\pro{memegenre}{ Soit $\{\Theta_{n}\inn =\{(f_{n},S_{n},\Gamma_{n},\rho_{n})\inn$ une suite
d'im\-mer\-sions {\'e}\-qui\-va\-riantes cocompactes telle que le genre de
$S_{n}/\Gamma_{n}$ soit born{\'e}. Soit $\xnn$ une suite de points telle
que $x_{n}\in S_{n}$ et  que $\{f_{n}(x_{n})\inn$ reste dans un compact.

Supposons que  $\{(f_{n},S_{n},x_{n})\inn$ converge vers
$(f_{\infty},S_{\infty},x_{\infty})$, nous avons alors les possibilit{\'e}s suivantes :
\begitem
\iti soit il existe $\Gamma_{\infty}$ et $\rho_{\infty}$ telle que
$(f_{\infty},S_{\infty},\Gamma_{\infty},\rho_{\infty})$ est cocompacte et
auquel cas $\{\Theta_{n}\inn$ est constante {\`a} partir d'un certain rang,

\itii soit $(f_{\infty},S_{\infty})$ est d{\'e}g{\'e}n{\'e}r{\'e}e ou tubulaire {\`a}
l'infini.
\enditem
}

Remarquons qu'il d{\'e}coule de \ref{toto} et \ref{tutu} que l'ensemble des $k$-surfaces
qui ne sont ni d{\'e}g{\'e}n{\'e}r{\'e}es, ni tubulaires {\`a} l'infini est dense
dans l'espace ${\cal N}$. Par ailleurs, a cause de l'unicite dans le
lemme de Morse, l'ensemble des $k$-surfaces
compactes de m{\^e}me genre est discret. Nous en d{\'e}duisons un raffinement de
\ref{densper}:

\theo{superdensper}{Soit $g\in{\mathbb N}$, alors l'ensemble des feuilles
compactes de genre plus grand que $g$ est dense dans ${\cal N}$}

\subsubsection{D{\'e}monstration de \ref{memegenre}}

Montrons tout d'abord la

\pro{coucou}{Soit $S$ une $k$-surface {\`a} courbure moyenne born{\'e}e, alors il existe une
constante $\e$ strictement positive telle que, pour tout $x\in S$, le
rayon d'injectivit{\'e} de $S$ en $x$ pour la m{\'e}trique induite de celle de $N$ est
minor{\'e} par $\e$}

\preu la d{\'e}monstration est imm{\'e}diate. La courbure de $S$ {\'e}tant
uniform{\'e}ment born{\'e}e, le  th{\'e}or{\`e}me de compacit{\'e} \ref{comrap} nous assure
que $S$ est {\`a} g{\'e}om{\'e}trie born{\'e}e. Le
rayon d'injectivit{\'e} est donc bien uniform{\'e}ment minor{\'e}.$\diamond$

\par
Nous pouvons maintenant d{\'e}montrer le th{\'e}or{\`e}me \ref{superdensper}. D'apr{\`e}s
l'{\'e}\-qua\-tion de Gauss, la courbure des m{\'e}triques sur $S_{n}$ est coinc{\'e}e
entre deux constantes strictement n{\'e}gatives.

Si le diam{\`e}tre de $S_{n}/\Gamma_{n}$ est unform{\'e}ment born{\'e}, nous en
d{\'e}duisons que $S_{n}/\Gamma_{n}$ converge vers une surface compacte, et
que n{\'e}cessairement la limite de
$$
\{(f_{n},S_{n},\Gamma_{n},\rho_{n})\inn
$$  
est une immersion
{\'e}quivariante cocompacte. En particulier, nous en d{\'e}duisons que  $\{(\Gamma_{n},\rho_{n})\inn$
devient constante {\`a} partir d'un certain rang, et ceci entra{\^\i}ne 
$\{\Theta_{n}\inn$ elle m{\^e}me est constante {\`a} partir d'une certain rang,
par le lemme de Morse.

Si le diam{\`e}tre de $S_{n}/\Gamma_{n}$ tend vers l'infini,  nous en d{\'e}duisons 
que le rayon d'injectivit{\'e} de $S_{\infty}$ est nul. Par \ref{coucou},
la courbure moyenne de $f_{\infty}$ n'est pas born{\'e}e. $(f_{\infty},S_{\infty})$
est donc d{\'e}g{\'e}n{\'e}r{\'e}e ou tubulaire {\`a} l'infini.$\diamond$

%%%%%%%%%%%%%%%%%%%%%%%%%%%%%%%%%%%%%%%%16%%%%%%%%%%%%%%%%%%%%%
\section{G{\'e}n{\'e}ricit{\'e} des feuilles denses}
Nous nous proposons de d{\'e}montrer

\theo{gendens}{L'ensemble des points de ${\cal N}$ par lesquel passent des
feuilles denses est une intersection d{\'e}nombrable d'ouverts denses}

\preu nous allons en fait d{\'e}montrer une r{\'e}sultat plus fort. Notons 
${\cal M}$ l'espace lamin{\'e} associ{\'e} {\`a} $M$, c'est-{\`a}-dire tel que 
${\cal N}={\cal M}/\pi_{1}(N)$, alosr, d{\'e}j{\`a} dans ${\cal M}$ une feuille g{\'e}n{\'e}rique est dense.

Remarquons tout d'abord que l'espace ${\cal N}$ poss{\'e}de une base
d{\'e}nombrable d'ouverts. Pour cela, il suffit  de trouver une application
continue injective de  ${\cal N}$ dans un espace {\`a} base
d{\'e}nombrable. Notons donc  $N(k)$, l'espace des $k$-jets de plans
de $N$, et  $N({\infty})$ la limite projective des $N(k)$. Cet espace
est {\`a} base d{\'e}nombrable et l'application naturelle de  ${\cal N}$ dans
$N(\infty )$ est injective, par ellipticit{\'e}. Nous en d{\'e}duisons bien s{\^u}r
que ${\cal M}$ est lui aussi {\`a} base d{\'e}nombrable.

Pour conclure, il nous faut donc montrer que si $x$ et $y$ sont deux
points de ${\cal M}$, il existe deux suites de points $\xnn$ et $\ynn$,
tendant respectivement vers $x$ et $y$, telles
que  pour tout $n$, $x_{n}$ et $y_{n}$ soient sur la m{\^e}me feuille ${\cal L}_{n}$.

Nous nous donnons donc deux $k$-surfaces point{\'e}es $(f_{\infty},S,x)$ et $(\phi_{\infty},\S,y)$,
dans $M$; notons $f$, et $\phi$ respectivement, leurs applications de
Gauss-Minkowski {\`a} valeurs dans $\bo_{\infty}M$. 

Soit maintenant $\Snn$, et $\Sinn$ des exhaustions de $S$
et $\S$ par des ouverts relativement compacts, contenant $x$ et $y$ respectivement.

Soit alors $(f_{n},\S_{n})$, et $(\phi_{n},\Sigma_{n})$,
les solutions des probl{\`e}mes de Plateau asymptotiques d{\'e}finis
respectivement par $(f,\S_{n})$, et $(\phi,\Sigma_{n})$.

D'apr{\`e}s la proposition \ref{toto}, $\{(f_{n},\S_{n},x)\inn$, et
$\{(\phi_{n},\Sigma_{n},y)\inn$, convergent respectivement vers
$(f_{\infty},\S,x)$ et
$(\phi_{\infty},\Sigma,y)$.

Par ailleurs, pour tout entier $p$, nous pouvons construire une famille de donn{\'e}es
asymptotiques $\{(h_{(p,n)},U_{(p,n)}\inn$ et des injections $i_{(p,n)}$ et
$j_{(p ,n)}$ de $S_{p}$ et $\S_{p}$ dans $U_{(p,n)}$ respectivement, telles
que 
l'on ait
$$
(f_{p},S_{p})=\liminf_{n\in {\mathbb N}}(h_{(p,n)},U_{(p,n)})~~ [ i_{(p,n)}],
$$
ainsi que   $$
(\phi_{p},S_{p})=\liminf_{n\in {\mathbb N}}(h_{(p,n)},U_{(p,n)})~~ [ j_{(p,n)}].
$$

Donnons l'esquisse de cette construction, il suffit de rejoindre
$f(S_{p}$ et $\phi(S_{p})$ respectivement par des rubans dont
l'{\'e}paisseur tend vers 0 quand $n$ tend vers l'infini.

Nous pouvons enfin nous d{\'e}brouiller pour que, pour tout $n$ et $p$,  le probl{\`e}me de Plateau
asymptotique $(h_{(p,n)},U_{(p,n)})$ poss{\`e}de une solution
$(H_{(p,n)},U_{(p,n)})$.

D'apr{\`e}s le lemme \ref{inclu}, la suite de $k$-surfaces point{\'e}es\-
$$
\{(H_{(p,n)},U_{(p,n)},i_{(p,n)}(x))\inn
$$ 
converge vers  $(f_{p},S_{p},x)$, et, respectivement, la suite
$$
\{(H_{(p,n)},U_{(p,n)},j_{(p,n)}(y))\inn
$$
converge vers $(\phi_{p},S_{p},y)$

Nous pouvons donc trouver une sous-suite $s(n)$ telle que les suites de
$k$-surfaces point{\'e}es

$$
\{(H_{(s(n),n)},U_{(s(n),n)},i_{(s(n),n)}(x))\inn
$$ 
et
$$
\{(H_{(s(n),n)},U_{(s(n),n)},j_{(s(n),n)}(y))\inn
$$
convergent respectivement vers $(f_{\infty},S,x)$ et $(\phi_{\infty},\S,
y)$.

En enon{\c c}ant ceci dans le cadre de ${\cal N}$, les suites de points de
${\cal N}$ d{\'e}finies par

$$
\{X_{n}\inn=\{(H_{(s(n),n)},U_{(s(n),n)},i_{(s(n),n)}(x))\inn
$$
et
$$
\{Y_{n}\inn=\{(H_{(s(n),n)},U_{(s(n),n)},j_{(s(n),n)}(y))\inn ,
$$
convergent vers $X=(f_{\infty},S,x)$ et $Y=(\phi_{\infty},\S,
y)$. Par ailleurs, pour tout $n$, $X_{n}$ et $Y_{n}$ appartiennent
{\`a} la m{\^e}me feuille ${\cal L}_{n}$, d{\'e}finie par
$$
{\cal L}_{n}=(H_{(s(n),n)},U_{(s(n),n)}).
$$
\par

C'est ce que nous voulions d{\'e}montrer.$\diamond$
%%%%%%%%%%%%%%%%%%%%%%%%%%%%%%%%%%%%%%%%17%%%%%%%%%%%%%%%%%%%%%%%%%
\section{Stabilit{\'e}}
Notre but dans cette section est de d{\'e}montrer un analogue du th{\'e}or{\`e}me de stabilit{\'e} pour
le flot g{\'e}od{\'e}sique des vari{\'e}t{\'e}s {\`a} courbure strictement n{\'e}gative~:

\theo{conjug}{Soit $N$ une vari{\'e}t{\'e} compacte de dimension 3, si $g$ est une m{\'e}trique {\`a} courbure 
strictement plus petite que -1, et $k\in]0,1[$, nous noterons ${\cal N}^{k}_g$, son espace 
des $k$-surfaces convexes.

Si $k$ et $l$ appartiennent {\`a} $]0,1[$, si $g$ et $h$ sont deux m{\'e}triques appartenant {\`a} la 
m{\^e}me composante connexe de l'espace 
des m{\'e}triques {\`a} courbure plus petite que -1, il existe alors un hom{\'e}omorphisme
$\Phi$ de ${\cal N}^{k}_g $ dans ${\cal N}^{l}_h$ envoyant feuille sur feuille.}

\subsection{Un lemme}

Nous voulons montrer le 

\lem{stabhoros}{Soit $g$ une m{\'e}trique {\`a} courbure strictement n{\'e}gative sur $N$ et $c$ une 
constante strictement positive, il existe un voisinage $U$ de $g$ pour la topologie $\ci$ tel
que, $\forall h, ~ {\bar h}\in U$, pour toute surface immerg{\'e}e compl{\`e}te {\`a} g{\'e}om{\'e}trie born{\'e}e $(f,S)$ 
dont les courbures principales (pour la m{\'e}trique $g$) sont plus grandes que $c$, si $(f,S)$ 
est horosph{\'e}rique {\`a} l'infini pour $h$, alors $(f,S)$ est horosph{\'e}rique
{\`a} l'infini  pour ${\bar h}$.}

\preu on peut choisir $U$ pour que $(f,S)$ ait une courbure plus grande que $c^2 /4$ pour 
toutes
les m{\'e}triques de $U$, et, 
en particulier, soit une surface convexe. 

Pour montrer le lemme, remarquons qu'il nous suffit de montrer que si $(f,S)$ est 
horosph{\'e}rique pour $h$,
 alors
$(f,S)$ est horosph{\'e}rique {\`a} l'infini pour ${\bar h}$. Or si, $(f,S)$ est horosph{\'e}rique, 
elle (ou plus exactement son image inverse dans le rev{\^e}tement universel) borde un ensemble
convexe ayant exactement un point {\`a} l'infini.  Le bord {\`a} l'infini {\'e}tant un invariant de 
quasi-isom{\'e}trie,
nous en d{\'e}duisons que le convexe bord{\'e} par $(f,S)$ {\`a} lui aussi un seul point {\`a} l'infini 
pour ${\bar h}$.

La surface $(f,S)$ est donc une pseudo-horosph{\`e}re pour la m{\'e}trique $g$ et par le lemme \ref{psudo},
elle est horosph{\'e}rique {\`a} l'infini pour ${\bar h}$.$\diamond$

\subsubsection{D{\'e}monstration du th{\'e}or{\`e}me \ref{conjug}}.

Il nous suffit bien s{\^u}r de montrer le th{\'e}or{\`e}me pour deux m{\'e}triques ${\bar h}$ et $h$ 
suffisamment proches au sens $C^{\infty}$. Nous allons les prendre dans un ouvert $U$
donn{\'e} par lemme pr{\'e}c{\'e}dent pour une constante $c$ ad{\'e}quate.

Soit tout ${\cal L}=(f,S)$ une feuille de ${\cal N}^{k}_{\bar h}$. Nous voyons ici $(f,S)$ comme une
immersion de $S$ dans $UN$ le fibr{\'e} unitaire de $N$. Notons alors $f_R$ l'application 
de $S$ dans $N$, donn{\'e}e par
$$
f_R ~:~s\mapsto exp(Rf(s)).
$$  

Nous savons ({\it cf} corollaire \ref{equic}), que, pour $R$
choisi suffisamment grand in\-d{\'e}\-pen\-dam\-ment de $S$, chacune des 
courbures principales de $f_R$ est plus grande qu'une constante $c$,  
avec $c>l^{1/2}$,

Nous en d{\'e}duisons que si $h$ est une m{\'e}trique suffisamment proche de ${\bar h}$, l'immersion 
$(f_R , S)$ est pour cette m{\'e}trique $h$ {\`a} courbure plus grande que $l$.

Nous sommes en mesure de construire notre hom{\'e}omorphisme $\Phi$. S{\'e}parons en deux cas,

(i) Si $(f,S)$ est un tube alors, par le lemme de Morse pour les g{\'e}od{\'e}siques,
$(f_R ,S)$ est tubulaire pour la m{\'e}trique $h$, soit $\gamma$ la g{\'e}od{\'e}sique correspondante et 
$T$ son tube. Nous avons donc une projection radiale naturelle 
$\pi$ de $f_R (S)$ sur $T$. Nous posons alors si $x\in S$
$$
\Phi (f,S,x) =(T,\pi\circ f_R (x)),
$$
o{\`u} nous voyons maintenant $(f,S,x)$, {\it resp.} $(T,\pi\circ f_R (x))$, comme un point
de  ${\cal N}^{k}_{\bar h} $, {\it resp.}  de ${\cal N}^{l}_h$.

(ii) Si $(f,S)$ n'est pas un tube, $(f_R ,S)$ n'est pas tubulaire pour $h$
( par le lemme de Morse pour les g{\'e}od{\'e}siques ). De plus, par le lemme \ref{stabhoros}
$(f_R ,S)$ n'est pas horosph{\'e}rique {\`a} l'infini. Notre lemme de Morse pour les surfaces convexes
\ref{morse}, nous permet alors de construire une $l$-surface lentille pour la m{\'e}trique 
$h$, $({\bar f},{\bar S})$
pour $(f_R ,S)$. Enfin, notons $\pi$ la projection qui a tout point de $(f_R ,S)$ associe 
son pied sur  $({\bar f},{\bar S})$.
Nous pouvons alors d{\'e}finir
$$
\Phi (f,S,x) =({\bar f},{\bar S},\pi\circ f_R (x)).
$$
A nouveau, nous voyons $(f,S,x)$, {\it resp.}$ ({\bar f},{\bar S},\pi\circ f_R (x))$, 
comme un point
de  ${\cal N}^{k}_{\bar h} $, {\it resp.}  de ${\cal N}^{l}_h$.

Par construction, $\Phi$ envoie feuille sur feuille. La continuit{\'e} de $\Phi$ provient de 
\ref{compalent}. 

Pour d{\'e}montrer, que $\Phi$ est une bijection, remarquons tout d'abord que $\Phi$
envoie bijectivement chaque feuille sur chaque feuille. Soit maintenant ${\bar \Phi}$ 
l'ap\-pli\-ca\-tion obtenue en inversant les r{\^o}les de $h$ et ${\bar h}$.
L'unicit{\'e} dans le lemme de Morse permet de montrer $\Phi\circ {\bar \Phi}$ envoie chaque 
feuille dans elle m{\^e}me.  L'application  $\Phi\circ {\bar \Phi}$ est donc une bijection, 
ce qui entra{\^\i}ne que ${\Phi}$ elle-m{\^e}me est une bijection.$\diamond$

%%%%%%%%%%%%%%%%%%%%%%%%%%%%%%%%%%%%%%%

\end{document}